\documentclass[11pt, a4paper]{article}

\usepackage[margin=2cm, bottom=3cm]{geometry}

\usepackage[utf8]{inputenc}    
\usepackage[T1]{fontenc}       
\usepackage{times}             

\usepackage{amsmath, amssymb, amsthm, mathtools}
\usepackage{stmaryrd}          
\usepackage{wasysym}           
\usepackage{cases}             
\usepackage{bigints}           

\usepackage{graphicx}          
\usepackage{float}             
\usepackage[labelfont=bf]{caption}  
\usepackage{placeins}          

\usepackage[colorlinks=true,
linkcolor=blue,
citecolor=blue,
urlcolor=blue]{hyperref}   
\usepackage{cleveref}
\usepackage{xcolor}

\usepackage{array, multirow, makecell}
\usepackage[shortlabels]{enumitem}
\usepackage[math]{cellspace}          
\usepackage[round, authoryear]{natbib}  

\usepackage{setspace}
\usepackage{empheq}
\usepackage{enumitem}

\newcommand*\dif{\mathop{}\!\mathrm{d}}

\renewcommand{\eqref}[1]{Equation~(\ref{#1})}

\theoremstyle{definition}

\newtheorem{lem}{Lemma}

\title{ 
	\textbf{Second-order prestress stability and third-order rigidity of polyhedral surfaces}
} 

\author{Zeyuan He and Daniel Robertz
\\~
	\\ \small Chair of Algebra and Number Theory, Faculty of Mathematics and Natural Sciences, RWTH Aachen University
	\\ \small zeyuan.he@rwth-aachen.de, daniel.robertz@rwth-aachen.de}
	
\date{}

\begin{document}
	
\maketitle



\begin{abstract}
There has been a longstanding confusion on the proper definition of higher order rigidity and flexibility in geometric constraint systems. Recently, an energy-based formulation of higher-order rigidity was introduced by Steven Gortler, Miranda Holmes-Cerfon, and Louis Theran (2025). In this article, we apply the framework to polyhedral surfaces and introduce new criteria for testing second-order prestress stability and third-order rigidity. Furthermore, we present comprehensive case studies of polyhedral surfaces exhibiting different levels of higher-order shakiness. These results advance the understanding of higher-order rigidity and flexibility in origami-inspired structures, which are also applicable to a broad class of near-mechanisms. Clarifying the transition from higher-order flexibility to finite mechanisms opens new directions for both theoretical investigation and mechanism design.
\end{abstract}



\maketitle

\section{Introduction} \label{section: introduction}



\begin{figure}
	\includegraphics[width=1\linewidth]{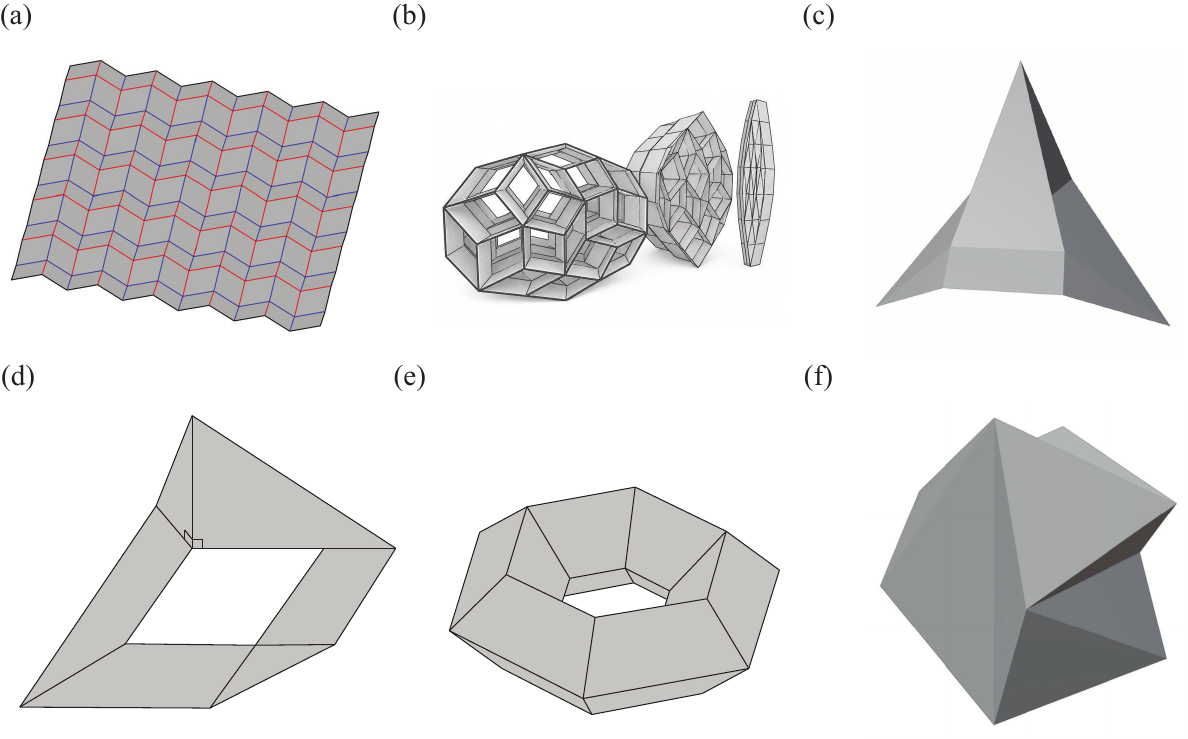}
	\caption{\label{fig: introduction}A gallery of polyhedral surfaces with diverse topologies, exhibiting flexibility or different degrees of rigidity. (a) The classical Miura-ori \citep{miura_method_1985}, notable for its flexibility despite being overconstrained, with all quadrilateral panels remaining rigid. (b) A tubular origami `honeycomb' \citep{tachi_rigid-foldable_2012} satisfying both rigid-panel constraints and global closure conditions; it remains flexible due to symmetry. (c) A newly discovered non-triangulated flexible polyhedron \citep{he_new_2025}, subject to stronger closure constraints in addition to panel rigidity, yet still flexible owing to its special geometric structure.
(d) A rigid polyhedral surface formed by a cycle of panels. It admits first-order flexes, i.e., non-trivial infinitesimal motions that preserve constraints to first order; such motions can be stabilized by appropriate self-stress.
(e) A rigid polyhedral torus that is first-order rigid, admitting no non-trivial infinitesimal flex.
(f) Jessen’s orthogonal icosahedron, which, like (d), is prestress stable.}
\end{figure}

Polyhedral surfaces in $\mathbb{R}^3$, composed of planar polygonal faces joined along their edges (Figure \ref{fig: introduction}), form a fundamental class of discrete geometric structures. Such geometry is essential to a wide range of applications, including architecture, robotics, metamaterials, deployable structures, biological assemblies, and computer graphics \citep{meloni_engineering_2021, misseroni_origami_2024}. Understanding the rigidity and flexibility of polyhedral surfaces has long been a central topic in mechanics and rigidity theory. A polyhedral surface is \textit{flexible} (Figure~\ref{fig: introduction}(a)–(c)) if it admits a continuous deformation in $\mathbb{R}^3$ that preserves the distances between any two vertices belonging to the same face and is not induced by a rigid motion of $\mathbb{R}^3$. Otherwise, it is said to be \textit{rigid} (Figure~\ref{fig: introduction}(d)–(f)). Classical rigidity theory focuses primarily on infinitesimal (i.e., first-order) rigidity and flexibility \citep{connelly_frameworks_2022}, which are determined by the existence of infinitesimal motions preserving the geometric constraints to first order. However, more subtle phenomena arise when a structure admits higher-order flexes that preserve these constraints to higher orders. In this work, we investigate the higher-order rigidity and prestress stability of polyhedral surfaces and develop methods for detecting and characterizing such higher-order flexes.

Roughly speaking, a polyhedral surface is first-order flexible if it admits a `speed' that maintains its geometric constraint. It is second-order flexible if, in addition to the speed, it admits an `acceleration' that maintains the geometric constraint, and third-order flexible if it admits an `acceleration of acceleration' (i.e. jerk), and so on. This interpretation aligns with the `classical' definition of $n$-th ($n \in \mathbb{Z}_+$) order flexibility for mechanisms \citep{rembs_verbiegungen_1933, efimov_theorems_1952, tarnai_higher-order_1989, sabitov_local_1992, salerno_how_1992, chen_order_2011, kuznetsov_underconstrained_2012}. A polyhedral surface with higher order flexibility becomes more and more `shaky' under perturbation, yet still has no actual motion as long as it is not flexible. This perspective naturally raises the question of how a highly shaky polyhedral surface approaches the transition from rigidity to flexibility. One way to approach this question is to introduce the idea of \textbf{critical order} for flexibility. The $n$-th order flexibility serves as a polynomial approximation to the flexibility. By the Artin approximation theorem, for each polyhedral surface, there is a sufficiently large $n^*$ such that the $n^*$-th order flexibility is equivalent to flexibility \citep{sabitov_local_1992, alexandrov_sufficient_1998}. This observation highlights the importance of developing criteria for higher-order rigidity and flexibility, since the iterative computation of higher-order flexibility is guaranteed to eventually determine flexibility.

However, the precise definition of $n$-th order flexibility, particularly for $n \ge 3$, in geometric constraint systems has long been a subject of ambiguity. Physically, one expects that structures with second-order flexibility form a subset exhibiting increased shakiness of those with first-order flexibility. Similarly, structures with third-order flexibility are anticipated to be a more shaky subset of those with second-order flexibility, and so on, until a sufficiently large $n^*$-th order that is equivalent to flexibility. This was later challenged by  \citet{connelly_higher-order_1994}, who demonstrated the existence of a flexible framework --- known as the double-Watt linkage --- that is classified as third-order rigid under the classical definition, despite being actually flexible. This `paradox' arises because the double-Watt linkage cannot initiate motion with nonzero speed. A flexible framework exhibiting this behaviour is also referred to as a \textit{cusp mechanism} in other literature. To address this issue, various studies have explored the limitation of the classical definition and proposed modifications \citep{gaspar_finite_1994, stachel_proposal_2007, nawratil_global_2025, tachi_proper_2024, nawratil_flexes_2025, gortler_higher_2025}. For a comprehensive overview of the historical context, see \citet{nawratil_global_2025, gortler_higher_2025}.

Prestress stability complements higher-order rigidity concepts that do not involve self-stress, by showing how internal forces can stabilize shaky structures. In its classical form, prestress stability describes how an appropriate self-stress can render a first-order flexible framework first-order rigid. This idea was formalized in rigidity theory for tensegrity structures, where appropriate self-stress ensures a positive-definite Hessian of the total energy, i.e. tangent stiffness \citep{connelly_second-order_1996}. Since then, the concept has found broad applications across structural engineering \citep{guest_stiffness_2006}, polymer and cytoskeletal networks \citep{kroy_glassy_2007}, soft matter \citep{marchetti_hydrodynamics_2013}, and mechanobiology \citep{ingber_tensegrity_1997}, providing a unified framework for tailoring structural stiffness to control static response and vibrational behavior. The breadth of applications naturally motivates the study of higher-order analogues of prestress stability, more specifically, whether approximate self-stresses can stabilize higher-order flexible frameworks to the corresponding level of rigidity; however, a systematic formulation of this concept remains absent.

Building on our previous work on local rigidity theory for polyhedral surfaces \citep{he_rigid_2022, hayakawa_panel-point_2024} up to second-order rigidity, this article develops criteria for determining second-order prestress stability and third-order rigidity. We adopt the energy-based definition recently introduced by \citet{gortler_higher_2025}, which overcomes several limitations of the classical definition. This formulation is particularly suited to our structural engineering perspective, as it quantifies how `floppy' a higher-order rigid polyhedral surface is through the growth rate of an associated energy function. It thereby provides a unified hierarchical framework for developing and comparing higher-order notions of rigidity and prestress stability. We further present a computational workflow for testing rigidity and prestress stability up to third order, as summarized in Figure~\ref{fig: flowchart}. The framework is illustrated through three polyhedral surfaces with distinct topologies and different forms of higher-order shakiness in Section \ref{section: examples}. Finally, Maple scripts are provided as supplementary material to facilitate customized symbolic computation, parametric design, and visualization.

\section{Modelling} \label{section: modelling}

\begin{figure}
	\includegraphics[width=1\linewidth]{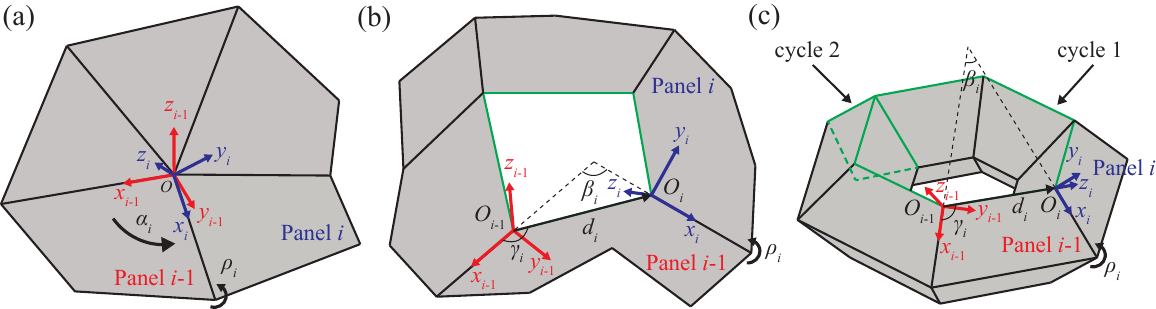}
	\caption{\label{fig: vertex and cycle}Graphical explanation on the measurement of `angular and distance misfits'. (a), (b), and (c) illustrate the rotation and possible translation of local coordinate systems around three polyhedral surfaces: a degree-5 vertex, a degree-5 hole, and a representative cycle on a $6 \times 4$ toroidal polyhedral surface. The local coordinate systems for panels $P_{i-1}$ and $P_i$ are depicted in red and blue, respectively, while the representative cycles in (b) and (c) are shown in green. Specifically, the closure constraint in (a) corresponds to a vertex and is expressed in the form of \eqref{eq: closure 1}; in (b), it applies to a hole and follows the form of \eqref{eq: closure 2}; in (c), it involves 24 vertices and 2 cycles, combining the forms of \eqref{eq: closure 1} and \eqref{eq: closure 2}. This figure is reproduced from the author's previous work \citet{he_rigid_2022}.}
\end{figure}

In this article, we use the \textbf{folding-angle model} to describe the configuration space of a polyhedral surface. In this model, panels are treated as rigid bodies connected by hinges, and the closure of hinge cycles serves as the geometric constraint. At present, this modelling framework is applicable to orientable surfaces. The dihedral angle at a hinge, measured with respect to a prescribed orientation and offset by $\pi$, is defined as the \textit{folding angle}. The interior angles of the polyhedral panels are referred to as \textit{sector angles}. The configuration of the entire 
polyhedral surface is described by the collection of folding angles, while the closure constraints ensure that the angular and distance misfits around each vertex or representative cycle vanish. More precisely, if a local orthogonal frame is attached to a panel and propagated through a sequence of rotations and translations around a cycle, the frame must return to its original position (Figure~\ref{fig: vertex and cycle}).

Each closure constraint associated with a vertex corresponds to three components of a moment, representing the work conjugates of the angular misfit. Similarly, each closure constraint associated with a representative cycle corresponds to three components of a moment and three components of a force. These quantities represent the work conjugates of the angular misfit and the translational distance misfit, respectively. The corresponding stresses satisfy equilibrium conditions at every hinge. Further details of this modelling framework can be found in the origami literature \citep{tachi_design_2012, tachi_rigid_2015, he_rigid_2022}.

Two other commonly used modelling frameworks are the panel-hinge model and the point-panel model. In the \textbf{panel-hinge model}, rigid planes (panels) are treated as bodies, while the shared lines where adjacent panels meet (hinges) are treated as constraints. Each panel undergoes a rigid-body motion consisting of a rotation and a translation. The hinge constraint requires that the motions of the two adjacent panels map their common hinge to the same spatial line. Each panel is therefore described by six variables (three for rotation and three for translation), while each hinge imposes five independent constraints. The space of overall rigid-body (i.e. trivial) motions in $\mathbb{R}^3$ has dimension six. This formulation is standard in the geometric rigidity literature \citep{crapo_statics_1982, whiteley_rigidity_2002}.

In the \textbf{point-panel model}, vertices are treated as points in space, while panels impose geometric constraints on groups of vertices. Each vertex is described by three variables corresponding to its Cartesian coordinates. A rigid panel constrains its vertices to remain coplanar and fixes the in-plane shape of the polygon through length constraints on selected pairs of vertices. For a panel with $k$ sides ($k \ge 4$), coplanarity can be enforced using $k-3$ constraints, requiring the signed volume of tetrahedra formed by four vertices to vanish. In addition, $2k-3$ length constraints are used to fix the in-plane geometry. As in the previous model, the space of trivial motions has dimension six. Further details of this modelling framework can be found in \citet{hayakawa_panel-point_2024}.

\begin{figure}
	\noindent \centering	\includegraphics[width=1\linewidth]{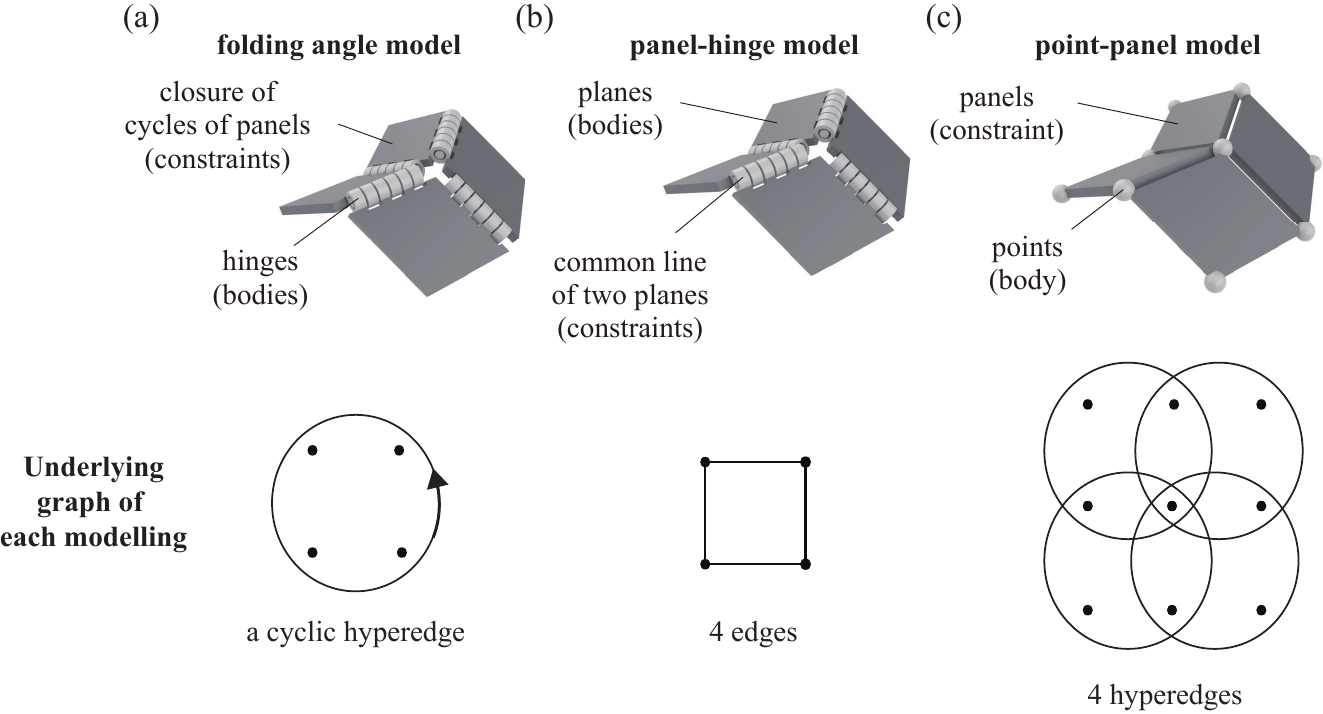}
	\caption{Three common modellings for polyhedral surfaces from the origami and geometric rigidity research community. The underlying graph for each model provides an abstract depiction of how bodies are connected with the geometric constraints. Here black dots represent bodies, while edges --- both standard edges connecting two bodies and hyperedges connecting multiple bodies --- are shown as lines and circles. A cyclic hyperedge, which constrain the sequence of body connections, is depicted with arrows. This figure is reproduced from the author's previous work \citet{hayakawa_panel-point_2024}.}
	\label{fig: description introduction}
\end{figure}

A graphical comparison of these three modelling frameworks is provided in Figure~\ref{fig: description introduction}.

We employ the folding-angle model in the subsequent analysis for several reasons: 
(1) the constraint associated with each hyperedge is symmetric with respect to all its variables, 
(2) the formulation is computationally convenient for evaluating derivatives of arbitrary order (Section \ref{section: n-th order derivatives}), 
(3) overall rigid-body motions are inherently excluded, and 
(4) since the constraint function possesses nonvanishing derivatives of arbitrarily high order, the theories of higher-order rigidity and prestress stability developed for this model are not tied to a particular finite-order constraint representation and can be extended naturally to geometric constraint systems with more general forms of constraint.

\clearpage

\clearpage

We now describe the geometric constraints of the folding-angle model in detail. 
Let $\rho \in \mathbb{R} \slash 2\pi\mathbb{Z}$ denote a folding angle and $\alpha \in (0,\pi)$ a sector angle. 
Here $\mathbb{R} \slash 2\pi\mathbb{Z}$ indicates that two angles $a,b\in\mathbb{R}$ are considered equivalent if 
$a=b+2k\pi$ for some $k\in\mathbb{Z}$. By convention, we often restrict $\rho$ to the interval $(-\pi,\pi]$. 
The closure constraint applies to every vertex and representative cycle, as illustrated in Figure~\ref{fig: vertex and cycle}. 
A vertex may be regarded as a degenerate cycle.

A representative cycle corresponds to a homology class in the first homology group, with the number of such independent cycles defined as the first Betti number. Homology theory was developed to analyze and classify manifolds based on their cycles (i.e., closed loops) that cannot be continuously deformed into each other. Informally, cycles that can be transformed into one another by continuous deformation belong to the same homology class in the first homology group. The first Betti number also represents the maximum number of cuts that can be made on a surface without dividing it into two separate pieces. For instance, the first Betti number is 0 for a sphere or a disk, 1 for a cylindrical surface, and 2 for a torus. 

The local coordinate systems are constructed as follows for deriving the closure constraint.
\begin{description}
	\item[\textbf{Around a vertex surrounded by $n$ hinges, Figure \ref{fig: vertex and cycle}(a)}] A local coordinate system is established on each panel $i$ ($i \in \mathbb{Z}_+,~ i \le n$), with the origin $O_i$ located at the central vertex. The $x$-axis lies along an interior crease, pointing outward from the origin, while the $z$-axis is normal to the panel. The directions of all $z$-axes are consistent with the surface orientation, aligning with the defined sign convention for folding angles. Specifically, the transition between the local coordinate systems of panels $i-1$ and $i$ consists of a rotation $\alpha_i$ about $z_{i-1}$, followed by a rotation $\rho_i$ about $x_i$. After a series of such rotations and translations (\eqref{eq: closure 1}), the transition matrix evaluates again to the identity.
\begin{align}
\label{eq: closure 1}
\overset{\curvearrowright}{\prod_{i=1}^{n}}
\left[
\begin{array}{ccc}
\cos \alpha_{i} & -\sin \alpha_{i} & 0\\
\sin \alpha_{i} &  \cos \alpha_{i} & 0\\
0               & 0                & 1
\end{array}
\right]
\left[
\begin{array}{ccc}
1 & 0 & 0\\
0 & \cos \rho_{i} & -\sin \rho_{i}\\
0 & \sin \rho_{i} & \cos \rho_{i}
\end{array}
\right]
=
\left[
\begin{array}{ccc}
1 & 0 & 0\\
0 & 1 & 0\\
0 & 0 & 1
\end{array}
\right]
\end{align}
where $\rho_i$ is the folding angle about axis $x_i$, $\alpha_i$ is the angle between axes $x_{i-1}$ and $x_i$ ($2 \le i \le n$), and $\alpha_1$ is the angle between axes $x_n$ and $x_1$. The symbol $\overset{\curvearrowright}{}$ indicates that the product is evaluated by multiplying matrices in an order where those with higher indices appear on the right.
	
	\item[\textbf{Around a cycle enclosed by $n$ hinges, Figure \ref{fig: vertex and cycle}(b) and (c)}] Local coordinate systems are constructed similarly. Each origin $O_i$ is located at a vertex of the cycle, with the 
$x$-axis aligned along a crease incident to the cycle, pointing outward from the origin, and the $z$-axis normal to the panel. The transition between the local coordinate systems of panels $i-1$ and $i$ includes a translation $[l_i \cos \gamma_i;~ l_i \sin \gamma_i;~ 0]$ measured in the coordinate system built on panel $i-1$, followed by a rotation $\beta_i$ about $z_{i-1}$, and a rotation $\rho_i$ about $x_i$. The angles $\beta_i$ and $\gamma_i$ can be expressed linearly in terms of all the sector angles. The corresponding transition matrix is provided in \eqref{eq: closure 2}.
\begin{equation}\label{eq: closure 2}
\begin{aligned}
\overset{\curvearrowright}{\prod_{i=1}^{n}}
\left[
\begin{array}{cccc}
\cos \beta_i & -\sin \beta_i & 0 & l_i \cos \gamma_i \\
\sin \beta_i & \cos \beta_i  & 0 & l_i \sin \gamma_i \\
0 & 0 & 1 & 0 \\
0 & 0 & 0 & 1
\end{array}
\right]
\left[
\begin{array}{cccc}
1 & 0 & 0 & 0 \\
0 & \cos \rho_i & -\sin \rho_i & 0 \\
0 & \sin \rho_i & \cos \rho_i  & 0 \\
0 & 0 & 0 & 1
\end{array}
\right]
=
\left[
\begin{array}{cccc}
1 & 0 & 0 & 0 \\
0 & 1 & 0 & 0 \\
0 & 0 & 1 & 0 \\
0 & 0 & 0 & 1
\end{array}
\right]
\end{aligned}
\end{equation}
where $\rho_i$ is the folding angle about axis $x_i$, $\beta_i$ is the angle between axes $x_{i-1}$ and $x_i$ ($2 \le i \le n$), and $\beta_1$ is the angle between axes $x_n$ and $x_1$. The vector $[l_i \cos \gamma_i;~ l_i \sin \gamma_i;~ 0]$ represents the position of $O_i$ measured in the local coordinate system of panel $i-1$. The symbol $\overset{\curvearrowright}{}$ indicates that the product is evaluated by multiplying matrices in an order where those with higher indices appear on the right. If the interior creases are concurrent, the cycle degenerates to a vertex.
\end{description} 

Following \citet{tachi_design_2012}, three independent constraints are selected from \eqref{eq: closure 1} and six from \eqref{eq: closure 2}. This choice is motivated by the fact that the derivatives (i.e.,  instant rotations) of the diagonal entries in \eqref{eq: closure 1} are consistently zero. Further, any instant rotation is skew-symmetric, and selecting only the off-diagonal entries ensures the formulation of equations that describe the (angular) velocity of the structure. Hence, we focus on the entries defined as follows:
\begin{equation}\label{eq: independent closure 1}
\begin{bmatrix}
\ast & \ast & f_2 \\
f_3 & \ast & \ast \\
\ast & f_1 & \ast
\end{bmatrix}
\coloneqq \overset {\curvearrowright} {\prod_{i=1}^{n}}
		\left[ \begin{array}{ccc}
			\cos \alpha_{i} & -\sin \alpha_{i} & 0\\
			\sin \alpha_{i} &  \cos \alpha_{i} & 0\\
			0                &                 0 & 1
		\end{array} \right]
		\left[ \begin{array}{ccc}
			1                &                 0 & 0             \\
			0                &     \cos \rho_{i} & -\sin \rho_{i}\\
			0                &     \sin \rho_{i} &  \cos \rho_{i}
		\end{array} \right]
\end{equation}

\begin{equation}\label{eq: independent closure 2}
\begin{bmatrix}
\ast & \ast & f_2 & f_4 \\
f_3 & \ast & \ast & f_5 \\
\ast & f_1 & \ast & f_6 \\
0 & 0 & 0 & 1
\end{bmatrix}
\coloneqq
		 \overset {\curvearrowright}{\prod_{i=1}^{n}}
			\left[ \begin{array}{cccc}
				\cos \beta_i & -\sin \beta_i & 0 & l_i \cos \gamma_i \\
				\sin \beta_i & \cos \beta_i & 0 & l_i \sin \gamma_i \\
				0 & 0 & 1 & 0   \\
				0 & 0 & 0 & 1
			\end{array} \right]
			\left[ \begin{array}{cccc}
				1 & 0 & 0 & 0 \\
				0 & \cos \rho_i & -\sin \rho_i & 0 \\
				0 & \sin \rho_i & \cos \rho_i & 0 \\
				0 & 0 & 0 & 1
			\end{array} \right]
\end{equation}
Here, $\ast$ denotes entries that are not relevant to the subsequent discussion. The constraint vector for the entire polyhedral surface $\boldsymbol{f} = \boldsymbol{0}$ consists of $3n_v + 6n_c$ conditions, where $n_v$ is the number of interior vertices, and $n_c$ is the number of representative cycles. Specifically, $\boldsymbol{f}$ is constructed by assembling three components for each vertex and six components for each cycle.


\section{$n$-th order derivative of the closure constraint} \label{section: n-th order derivatives}

In this section, we show that the $n$-th order derivative of the closure constraint $\boldsymbol{f}$ can be calculated explicitly from the coordinates of vertices in \textbf{one (global) coordinate system}, rather than from coordinate systems attached to every panel in the derivation of \eqref{eq: closure 1} and \eqref{eq: closure 2}. 

For a polyhedral surface with multiple vertices and cycles, the $i$-th order derivative $(i \in \mathbb{Z}_+)$ of $\boldsymbol{f}$ is an order $i+1$ tensor, which is assembled by 3 components of every vertex and 6 components of every cycle. Let $\boldsymbol{f}^v$ be the 3 components for a vertex:
\begin{equation*}
	\dfrac{\partial \boldsymbol{f}^v}{\partial \rho_j} = x_j = \begin{bmatrix}
			x_{j1} \\
			x_{j2} \\
			x_{j3} 
		\end{bmatrix},
\end{equation*}
where $x_j$ is the direction vector $\left(\|x_j\|=1\right)$ along the crease pointing outward the vertex, as illustrated in Figure \ref{fig: vertex and cycle}. Next, 
\begin{equation*}
	\begin{gathered}
	\dfrac{\partial^2 \boldsymbol{f}^v}{\partial \rho_k \partial \rho_j} = \begin{bmatrix}
			\left(x_k^\times x_j^\times\right)_{32} \\[4pt]
			\left(x_k^\times x_j^\times\right)_{13} \\[4pt]
			\left(x_k^\times x_j^\times\right)_{21}
		\end{bmatrix} = \begin{bmatrix}	x_{j3}x_{k2} \\[4pt]
	x_{j1}x_{k3} \\[4pt]
	x_{j2}x_{k1}
	\end{bmatrix} \mathrm{~for~} k \le j, ~~\left(\dfrac{\dif^2 \boldsymbol{f}^v}{\dif \boldsymbol{\rho}^2}\right)_{jk} = \left(\dfrac{\dif^2 \boldsymbol{f}^v}{\dif \boldsymbol{\rho}^2}\right)_{kj},
	\end{gathered}
\end{equation*}

where $_{32}, ~_{13}, ~_{21}$ denote the positions of entries in the matrix, and
\begin{equation*}
	\begin{gathered}
		x_j^\times = \begin{bmatrix}
			0 & -x_{j3} & x_{j2}\\
			x_{j3} & 0 & -x_{j1} \\
			-x_{j2} & x_{j1} & 0
		\end{bmatrix}.
	\end{gathered}
\end{equation*}

Furthermore,
\begin{equation*}
	\dfrac{\partial^3 \boldsymbol{f}^v}{\partial \rho_l \partial \rho_k \partial \rho_j} = \begin{bmatrix}
		\left(x_l^\times x_k^\times x_j^\times\right)_{32} \\[4pt]
		\left(x_l^\times x_k^\times x_j^\times\right)_{13} \\[4pt]
		\left(x_l^\times x_k^\times x_j^\times\right)_{21}
	\end{bmatrix} \mathrm{~~for~~} l \le k \le j,
\end{equation*}
\begin{equation*}
	\left(\dfrac{\dif^3 \boldsymbol{f}^v}{\dif \boldsymbol{\rho}^3}\right)_{lkj} = \left(\dfrac{\dif^3 \boldsymbol{f}^v}{\dif \boldsymbol{\rho}^3}\right)_{jlk} = \left(\dfrac{\dif^3 \boldsymbol{f}^v}{\dif \boldsymbol{\rho}^3}\right)_{kjl} = \left(\dfrac{\dif^3 \boldsymbol{f}^v}{\dif \boldsymbol{\rho}^3}\right)_{ljk} = \left(\dfrac{\dif^3 \boldsymbol{f}^v}{\dif \boldsymbol{\rho}^3}\right)_{klj} = \left(\dfrac{\dif^3 \boldsymbol{f}^v}{\dif \boldsymbol{\rho}^3}\right)_{jkl}.
\end{equation*}
For every cycle, let $\boldsymbol{f}^c$ be the 3 components other than $\boldsymbol{f}^v$:
\begin{equation*}
	\dfrac{\partial \boldsymbol{f}^c}{\partial \rho_j} = - x_j \times O_j,
\end{equation*}
where $O_j$ is the coordinate vector of the vertex on the cycle incident to crease $x_j$. The higher order derivatives of $\boldsymbol{f}^c$ are:
\begin{equation*}
	\dfrac{\partial^2 \boldsymbol{f}^c}{\partial \rho_k \partial \rho_j} = - x_k \times \left(x_j \times O_j\right) \mathrm{~~for~~} k \le j, ~~\left(\dfrac{\dif^2 \boldsymbol{f}^c}{\dif \boldsymbol{\rho}^2}\right)_{jk} = \left(\dfrac{\dif^2 \boldsymbol{f}^c}{\dif \boldsymbol{\rho}^2}\right)_{kj},
\end{equation*}
\begin{equation*}
	\dfrac{\partial^3 \boldsymbol{f}^c}{\partial \rho_l \partial \rho_k \partial \rho_j} = -x_l \times (x_k \times (x_j \times O_j)) \mathrm{~~for~~} l \le k \le j,
\end{equation*}
\begin{equation*}
	\left(\dfrac{\dif^3 \boldsymbol{f}^c}{\dif \boldsymbol{\rho}^3}\right)_{lkj} = \left(\dfrac{\dif^3 \boldsymbol{f}^c}{\dif \boldsymbol{\rho}^3}\right)_{jlk} = \left(\dfrac{\dif^3 \boldsymbol{f}^c}{\dif \boldsymbol{\rho}^3}\right)_{kjl} = \left(\dfrac{\dif^3 \boldsymbol{f}^c}{\dif \boldsymbol{\rho}^3}\right)_{ljk} = \left(\dfrac{\dif^3 \boldsymbol{f}^c}{\dif \boldsymbol{\rho}^3}\right)_{klj} = \left(\dfrac{\dif^3 \boldsymbol{f}^c}{\dif \boldsymbol{\rho}^3}\right)_{jkl}.
\end{equation*}
The $n$-th order derivatives are explicit by applying induction to $n$. The detailed derivation is presented in Appendix \ref{app: derivatives}.

\section{Second-order prestress stability and third-order rigidity} \label{section: second and third}

Instead of defining higher-order rigidity purely through derivatives of the constraint equations, \citet{gortler_higher_2025} characterize rigidity by measuring how rapidly physical energy increases when the framework is perturbed from its equilibrium configuration $\boldsymbol{\rho}$. The energy functional is assumed to be a \textit{stiff-bar energy} $E_b$, meaning that for every closure constraint energy, the self-stress is zero and the stiffness matrix (Hessian) is positive-definite, so that every closure constraint behaves mechanically like a standard unloaded spring. Within this framework, the equilibrium configuration $\boldsymbol{\rho}$ is a local minimum of $E_b$. Moreover, $\boldsymbol{\rho}$ is rigid if and only if this minimum is strict; that is, there exists a neighbourhood $B$ of $\boldsymbol{\rho}$ such that $E_b(\boldsymbol{\phi}) > E_b(\boldsymbol{\rho})$ for all $\boldsymbol{\phi} \in B \setminus \{\boldsymbol{\rho}\}$. Conversely, $\boldsymbol{\rho}$ is flexible if and only if the minimum is non-strict, meaning that for every neighbourhood $B$ of $\boldsymbol{\rho}$, there exists $\boldsymbol{\phi} \in B \setminus \{\boldsymbol{\rho}\}$ with $E_b(\boldsymbol{\phi}) = E_b(\boldsymbol{\rho})$.

Consider a \textit{trajectory} $\boldsymbol{\rho}(t)$, $t \in [0,1]$, starting from $\boldsymbol{\rho}(0)$ = $\boldsymbol{\rho}$, expressed as a power series
\begin{equation} \label{eq: trajectory}
\boldsymbol{\rho}(t)
=
\boldsymbol{\rho}
+ \boldsymbol{\rho}' t
+ \boldsymbol{\rho}'' t^2
+ \cdots
+ \boldsymbol{\rho}^{(k)} t^k
+ \cdots .
\end{equation}

The \textit{energy growth} along the trajectory at $\boldsymbol{\rho}$ is characterized, for $t \to 0$, by
\begin{equation}
E_b\big(\boldsymbol{\rho}(t)\big) - E_b(\boldsymbol{\rho})
\sim
\lVert \boldsymbol{\rho}(t) - \boldsymbol{\rho} \rVert^{s},
\quad \text{for some } s \in \mathbb{Q}_{+}.
\end{equation}

If the energy increases at least to order $s$ for all trajectories, and there exist trajectories along which the growth is exactly of order $s$, then the energy grows \textit{tightly} at order $s$ and the \textit{order of rigidity} is defined as $s/2$. It is shown that this rigidity order is independent of the particular choice of stiff-bar energy.

Moreover, the computation of the rigidity order can be reduced from arbitrary trajectories to \textit{flexes}. A $(j,k)$-\textit{flex} is a trajectory whose first nonzero term occurs at order $j$, while the geometric constraints are satisfied up to order $k$ ($j \le k$ from Taylor expansion). With this notion, for a rigid framework, the rigidity order satisfies:
\begin{equation} \label{eq: rigidity order and flex}
\dfrac{s}{2}
=
\max
\left\{
\frac{k+1}{j}
:
\text{there exists a } (j,k)\text{-flex} \text{ of the framework}
\right\},
\end{equation}
and for a flexible framework, $s = \infty$.

This energy-based definition resolves limitations of classical higher-order rigidity notions. In particular, the absence of a $(1,k)$-flex does not necessarily imply $k$-th order rigidity. For example, a cusp mechanism with no $(1, k)$-flex but admitting a $(2,\infty)$-flex is correctly classified as flexible within this definition.

From \eqref{eq: rigidity order and flex}, the classical criteria for higher-order rigidity can, in principle, be derived by analyzing admissible $(j,k)$-flexes. While the framework extends naturally to higher orders, the criteria and computation of $n$-th order flexes become increasingly intricate. For clarity and concreteness, we therefore focus on a detailed treatment up to third-order rigidity.  

\subsection{Rigidity tests}

We denote by $\mathcal{K}_b = \mathrm{Null}(\dif \boldsymbol{f}/\dif \boldsymbol{\rho})$ the right null space of the rigidity matrix.

\begin{enumerate}[label={[{\arabic*}]}] \item First-order rigidity is characterized by the absence of any $(j,j)$-flex for all $j$. It follows that this is equivalent to the absence of a $(1,1)$-flex, and hence to the condition $\mathrm{dim}(\mathcal{K}_b)=0$.

\item Second-order rigidity is characterized by the existence of a $(j,2j-1)$-flex for some $j$, together with the absence of any $(j,2j)$-flex. It follows that this condition reduces to the existence of a $(1,1)$-flex and the non-existence of a $(1,2)$-flex. The key reparametrization method underlying this reduction traces back to \citet{connelly_second-order_1996}. More precisely,
\begin{subequations}\label{eq:second-order-rigidity-test-1}
\begin{empheq}[left=\empheqlbrace]{align}
&\dim(\mathcal{K}_b) > 0. \label{eq:second-order-rigidity-test-1a}\\[6pt]
&\text{For any } \boldsymbol{\rho}' \in \mathcal{K}_b \setminus \{\boldsymbol{0}\},
\text{ there does not exist } \boldsymbol{\rho}'' \text{ satisfying } \notag\\
& \qquad\qquad \frac{\dif \boldsymbol{f}}{\dif \boldsymbol{\rho}} \cdot \boldsymbol{\rho}''
+ \frac{\dif^2 \boldsymbol{f}}{\dif \boldsymbol{\rho}^2}
 \cdot \left(\boldsymbol{\rho}' \otimes \boldsymbol{\rho}'\right)
= \boldsymbol{0}. \label{eq:second-order-rigidity-test-1b}
\end{empheq}
\end{subequations}

Here $\cdot$ denotes the contraction (inner product) between tensors. Those $\boldsymbol{\rho}'$ that meet \eqref{eq:second-order-rigidity-test-1b} are called \textit{extendable}. In index notation, the above equation reads
\begin{equation*}
\dfrac{\partial f_i}{\partial \rho_j} \rho_j'' 
+ \dfrac{\partial^2 f_i}{\partial \rho_j \partial \rho_k} \rho_j'\rho_k'  
= 0.
\end{equation*}

In practice, we often use the so-called \textit{stress test} \eqref{eq:second-order-rigidity-test-2b}, which is equivalent to \eqref{eq:second-order-rigidity-test-1b}. This equivalence follows from the fact that a linear system admits a solution if and only if its right-hand side is orthogonal to every vector in the left null space, which is also known as the the Fredholm alternative.

\begin{subequations}\label{eq:second-order-rigidity-test-2}
\begin{empheq}[left=\empheqlbrace]{align}
&\dim(\mathcal{K}_b) > 0.
\label{eq:second-order-rigidity-test-2a}\\[6pt]
&\begin{aligned}
&\text{For any } \boldsymbol{\rho}' \in \mathcal{K}_b \setminus \{\boldsymbol{0}\}, 
\text{ there exists self-stress } \boldsymbol{\omega} 
\text{ (as a function of } \boldsymbol{\rho}') \text{ such that}\\[4pt]
&\qquad\qquad\qquad\qquad
\boldsymbol{\omega} \cdot 
\dfrac{\dif^2 \boldsymbol{f}}{\dif \boldsymbol{\rho}^2} 
\cdot \left(\boldsymbol{\rho}' \otimes \boldsymbol{\rho}'\right) 
> 0.
\end{aligned}
\label{eq:second-order-rigidity-test-2b}
\end{empheq}
\end{subequations}

In index notation, \eqref{eq:second-order-rigidity-test-2b} reads
\begin{equation*}
 \omega_i \dfrac{\partial^2 f_i}{\partial \rho_j \partial \rho_k} \rho_j'\rho_k' > 0.
\end{equation*}

Following the terminology in geometric rigidity theory, we define the \textit{stress matrix}
\begin{equation*}
\Omega := \boldsymbol{\omega} \cdot \dfrac{\dif^2 \boldsymbol{f}}{\dif \boldsymbol{\rho}^2}.
\end{equation*}

\item Third-order rigidity is characterized by the existence of a $(j,3j-1)$-flex for some $j$, together with the absence of any $(j,3j)$-flex. When $\mathcal{K}_b$ is one-dimensional, this condition reduces to the existence of a $(1,2)$-flex and the non-existence of a $(1,3)$-flex. More precisely,  
\begin{subequations}\label{eq:third-order-rigidity-test-1}
\begin{empheq}[left=\empheqlbrace]{align}
&\dim(\mathcal{K}_b)=1. \label{eq:third-order-rigidity-test-1a}\\[6pt]
&\text{Let } \boldsymbol{\rho}' \in \mathcal{K}_b, \|\boldsymbol{\rho}'\|=1, 
\text{ there exists a solution } \boldsymbol{\rho}'' \text{ satisfying} \notag\\[3pt]
&\qquad\qquad\qquad\qquad
\dfrac{\dif \boldsymbol{f}}{\dif \boldsymbol{\rho}} \cdot \boldsymbol{\rho}''
+ \dfrac{\dif^2 \boldsymbol{f}}{\dif \boldsymbol{\rho}^2} 
\cdot \left(\boldsymbol{\rho}' \otimes \boldsymbol{\rho}'\right)
= \boldsymbol{0}.
\label{eq:third-order-rigidity-test-1b}\\[6pt]
&\text{For any } \boldsymbol{\rho}'',  \text{there does not exist } \boldsymbol{\rho}''' \text{ satisfying} \notag\\[3pt]
&\qquad\qquad
\dfrac{\dif \boldsymbol{f}}{\dif \boldsymbol{\rho}} \cdot \boldsymbol{\rho}'''
+ 3\,\dfrac{\dif^2 \boldsymbol{f}}{\dif \boldsymbol{\rho}^2} 
\cdot \left(\boldsymbol{\rho}'' \otimes \boldsymbol{\rho}'\right)
+ \dfrac{\dif^3 \boldsymbol{f}}{\dif \boldsymbol{\rho}^3} 
\cdot (\boldsymbol{\rho}')^{\otimes 3}
= \boldsymbol{0}.
\label{eq:third-order-rigidity-test-1c}
\end{empheq}
\end{subequations}
In index notation, \eqref{eq:third-order-rigidity-test-1c} reads:
\begin{equation*}
\dfrac{\partial f_i}{\partial \rho_j} \rho_j''' 
+ 3\dfrac{\partial^2 f_i}{\partial \rho_j \partial \rho_k} \rho_j''\rho_k' + \dfrac{\partial^3 f_i}{\partial \rho_j \partial \rho_k \partial \rho_l} \rho_j'\rho_k'\rho_l' 
= 0.
\end{equation*}

Similarly, we have the stress test \eqref{eq:third-order-rigidity-test-2} equivalent to \eqref{eq:third-order-rigidity-test-1}:

\begin{subequations}\label{eq:third-order-rigidity-test-2}
\begin{empheq}[left=\empheqlbrace]{align}
&\dim(\mathcal{K}_b)=1. \label{eq:third-order-rigidity-test-2a}\\[6pt]
&\text{Let } \boldsymbol{\rho}' \in \mathcal{K}_b, \|\boldsymbol{\rho}'\|=1, 
\text{ for every self-stress } \boldsymbol{\omega}, \text{we have~}  \notag\\[3pt]
&\qquad\qquad\qquad\qquad
\Omega \cdot \left(\boldsymbol{\rho}' \otimes \boldsymbol{\rho}'\right)
= \boldsymbol{0}.
\label{eq:third-order-rigidity-test-2b}\\[6pt]
&\text{Solve } \boldsymbol{\rho}'' \text{ from } \eqref{eq:third-order-rigidity-test-1b}, \text{ for every } \boldsymbol{\rho}'', \text{there exists self-stress~} \boldsymbol{\omega} \text{~such that~}
\notag\\[3pt]
&\qquad\qquad
 3\,\Omega \cdot \left(\boldsymbol{\rho}'' \otimes \boldsymbol{\rho}'\right)
+ \boldsymbol{\omega} \cdot \dfrac{\dif^3 \boldsymbol{f}}{\dif \boldsymbol{\rho}^3} \cdot (\boldsymbol{\rho}')^{\otimes 3}
> 0
\label{eq:third-order-rigidity-test-2c}
\end{empheq}
\end{subequations}

In index notation, \eqref{eq:third-order-rigidity-test-2c} reads:
\begin{equation*}
 3\omega_i\dfrac{\partial^2 f_i}{\partial \rho_j \partial \rho_k} \rho_j''\rho_k' + \omega_i\dfrac{\partial^3 f_i}{\partial \rho_j \partial \rho_k \partial \rho_l} \rho_j'\rho_k'\rho_l' > 0.
\end{equation*}
\end{enumerate}

The first- and second-order derivatives introduced in Section~\ref{section: n-th order derivatives}, together with the characterizations of first- and second-order flexibility developed here, suffice to recover classical graphical results for simply-connected polyhedral surfaces. First-order flexibility is equivalent to the existence of a reciprocal-parallel diagram of the fold lines, while second-order flexibility is further equivalent to the zero-area property of this reciprocal diagram \citep{watanabe_method_2009,tachi_design_2012,demaine_zero-area_2016}. The first-order equivalence follows from the fact that the instantaneous relative rotation of two adjacent panels must occur about their common fold line, an observation made much earlier by \citet{crapo_statics_1982}.

\subsection{Generalized $n$-th order flexes}

Suppose that a configuration $\boldsymbol{\rho}$ is first-order flexible. We use the tuple of coefficients below to represent a trajectory defined in \eqref{eq: trajectory}.
\begin{equation*}
    \left(\boldsymbol{\rho}',\, \boldsymbol{\rho}'',\, \cdots,\, \boldsymbol{\rho}^{(k)},\, \cdots \right)
\end{equation*}

In the classical kinematic definition, a first-order flex is a nonzero vector $\boldsymbol{a} \in \mathcal{K}_b \setminus \{\boldsymbol{0}\}$. The energy-based formulation generalizes this notion by allowing the first $j-1$ coefficients to vanish. In this setting, a \textit{$j$-active first-order flex} is defined as a $(j,j)$-flex. Specifically, for $\boldsymbol{a} \in \mathcal{K}_b \setminus \{\boldsymbol{0}\}$, such a flex can be expressed as
\begin{equation}\label{eq: first-order-flexes}
\left(\boldsymbol{0},\, \cdots,\, \boldsymbol{0},\, 
\underset{j\text{-th position}}{\boldsymbol{a}}\right).
\end{equation}

Suppose that a configuration $\boldsymbol{\rho}$ is second-order flexible, which means we can find $\boldsymbol{a} \in \mathcal{K}_b \setminus \{\boldsymbol{0}\}$ and $\boldsymbol{b}$ such that:
\begin{equation*}
    \dfrac{\dif \boldsymbol{f}}{\dif \boldsymbol{\rho}} \cdot \boldsymbol{b} 
+ \dfrac{\dif^2 \boldsymbol{f}}{\dif \boldsymbol{\rho}^2} \cdot \left(\boldsymbol{a} \otimes \boldsymbol{a}\right) 
= \boldsymbol{0}.
\end{equation*}

In the classical kinematic definition, a second-order flex is a tuple $(\boldsymbol{a}, ~\boldsymbol{b})$. Similarly, now in the energy-based definition, a $j$-active \textit{second-order flex} refers to a $(j,2j)$-flex. The first $j-1$ coefficients are allowed to vanish. Such flexes can be represented by the following ordered tuples, whose entries can be read directly from the higher-order derivatives $\dif^k \boldsymbol{f}/\dif t^k$ given in Appendix~\ref{app: Faa}.

\begin{equation} \label{eq: second-order flexes}
\left(\boldsymbol{0},\, \cdots,\, \boldsymbol{0},\, \underset{j\text{-th position}}{\boldsymbol{a}},\, \boldsymbol{\rho}^{(j+1)},\, \cdots,\, \boldsymbol{\rho}^{(2j-1)}, \, \underset{2j\text{-th position}}{\frac{(2j)!}{2 (j!)^2} \boldsymbol{b}} \right)
\end{equation}
where $\boldsymbol{\rho}^{(k)} \in \mathcal{K}_b$ for all $j+1 \le k \le 2j-1$.

Suppose that a configuration $\boldsymbol{\rho}$ is third-order flexible with $\dim(\mathcal{K}_b)=1$, which means we can find $\boldsymbol{a} \in \mathcal{K}_b \setminus \{\boldsymbol{0}\}$, $\boldsymbol{b}$, and $\boldsymbol{c}$ such that:
\begin{equation*}
\begin{cases}
    \dfrac{\dif \boldsymbol{f}}{\dif \boldsymbol{\rho}} \cdot \boldsymbol{b} 
+ \dfrac{\dif^2 \boldsymbol{f}}{\dif \boldsymbol{\rho}^2} \cdot \left(\boldsymbol{a} \otimes \boldsymbol{a}\right) 
= \boldsymbol{0}, \\
\dfrac{\dif \boldsymbol{f}}{\dif \boldsymbol{\rho}} \cdot \boldsymbol{c}
+ 3\,\dfrac{\dif^2 \boldsymbol{f}}{\dif \boldsymbol{\rho}^2} \cdot \left(\boldsymbol{b} \otimes \boldsymbol{a} \right)
+ \dfrac{\dif^3 \boldsymbol{f}}{\dif \boldsymbol{\rho}^3} \cdot \boldsymbol{a}^{\otimes 3}
= \boldsymbol{0}.
\end{cases}
\end{equation*}

A \textit{j-active third-order flex} can be represented by the following ordered tuple, whose entries can be read directly from the higher-order derivatives $\dif^k \boldsymbol{f}/\dif t^k$ given in Appendix~\ref{app: Faa}: 
\begin{equation} \label{eq: third-order flexes}
\left(\boldsymbol{0},\, \cdots,\, \boldsymbol{0},\, \underset{j\text{-th position}}{\boldsymbol{a}},\, \boldsymbol{a},\, \cdots,\, \boldsymbol{a}, \, \underset{2j\text{-th position}}{\frac{(2j)!}{2 (j!)^2} \boldsymbol{b}}, \, \frac{(2j)!}{2 (j!)^2} \boldsymbol{b} , \, \cdots, \, \frac{(2j)!}{2 (j!)^2} \boldsymbol{b} , \, \underset{3j\text{-th position}}{\dfrac{(3j) !}{6 (j !)^3}\boldsymbol{c}} \right)
\end{equation}

Note that, by \eqref{eq: second-order flexes} and \eqref{eq: third-order flexes}, a $(j,2j-1)$-flex can be constructed from any first-order flex, and a $(j,3j-1)$-flex can be constructed from any second-order flex. This recursive structure explains why, although fractional orders of rigidity are theoretically possible, they do not occur below third order. In contrast, no analogous construction is guaranteed from a third-order flex to a $(j,4j-1)$-flex, suggesting that fractional orders of rigidity may arise above third order.

\subsection{Prestress stability tests}

Next we discuss prestressed energy $E = E_\omega + E_b$, where $E_b$ is the aforementioned stiff-bar energy and $E_\omega$ is the additional energy term only dependent on self-stress $\boldsymbol{\omega}$. Unlike stiff-bar energies, which are guaranteed to attain a local minimum at an equilibrium configuration, a prestressed energy $E$ does not necessarily have this property when the self-stress is nonzero, even if the stiffness matrix is positive definite. 

Suppose that a prestressed energy $E$ attains a strict local minimum at $\boldsymbol{\rho}$ and grows tightly at order $s$ in a neighborhood of $\boldsymbol{\rho}$. Since the growth behavior depends on both $\boldsymbol{\omega}$ and $E_b$, different prestressed energies may lead to different values of $s$. This motivates defining the \textit{order of prestress stability} by taking the infimum of $s/2$ over all prestressed energies $E$ for which $\boldsymbol{\rho}$ is a strict local minimum. Physically, this corresponds to identifying the most effective prestressed energy for the geometric constraint system.

It is shown that the problem of identifying the most effective prestressed energy can be asymptotically reduced to selecting an appropriate self-stress $\boldsymbol{\omega}$. This follows from the fact that, when $\boldsymbol{\omega}$ has sufficiently small magnitude relative to the material stiffness (i.e., ~Hessian of $E_b$), the choice of the material stiffness only introduces a negligible perturbation to the resulting prestress stability order. Consequently, the dependence on $E_b$ becomes negligible, and the rigidity order is effectively governed by the prestress. The assumption of small prestress (or equivalently large material stiffness) is classical in the literature \citep{connelly_second-order_1996} and is also consistent with practical considerations in structural engineering.

To ensure that a prestressed energy $E$ attains a strict local minimum at an equilibrium configuration $\boldsymbol{\rho}$ when the Hessian is only positive semidefinite, \citet{gortler_higher_2025} develop a generalized `$2k$-th order test' based on so-called `indicative trajectories'. This framework yields conditions on the self-stress $\boldsymbol{\omega}$ under the aforementioned small prestress assumption. Their approach addresses a fundamental limitation of higher-order tests, as can be seen from the Taylor expansion at $\boldsymbol{\rho}$ below. Note that all energy functions are assumed to be analytic, i.e., smooth and locally expressible as a convergent power series.

\begin{equation}
\begin{aligned}
E(\boldsymbol{\rho}+\Delta \boldsymbol{\rho})
& = E(\boldsymbol{\rho})
+ \frac{\dif E}{\dif \boldsymbol{\rho}} \cdot \Delta \boldsymbol{\rho} + \frac{1}{2}\, \frac{\dif^2 E}{\dif \boldsymbol{\rho}^2} \cdot (\Delta \boldsymbol{\rho} \otimes \Delta \boldsymbol{\rho}) + \frac{1}{6}\, \frac{\dif^3 E}{\dif \boldsymbol{\rho}^3} \cdot (\Delta \boldsymbol{\rho})^{\otimes 3} + \frac{1}{24}\, \frac{\dif^4 E}{\dif \boldsymbol{\rho}^4} \cdot (\Delta \boldsymbol{\rho})^{\otimes 4}
+ \mathcal{O}(\|\Delta \boldsymbol{\rho}\|^5)
\end{aligned}
\end{equation}
First of all, equilibrium means:
\begin{equation} \label{eq: equilibrium}
    \dfrac{\dif E}{\dif \boldsymbol{\rho}} = \dfrac{\dif E}{\dif \boldsymbol{f}} \cdot  \dfrac{\dif \boldsymbol{f}}{\dif \boldsymbol{\rho}} = \boldsymbol{\omega} \cdot  \dfrac{\dif \boldsymbol{f}}{\dif \boldsymbol{\rho}} = \boldsymbol{0}.
\end{equation}
Consequently,
\begin{equation}
\begin{aligned}
E(\boldsymbol{\rho} + \Delta \boldsymbol{\rho})
& = E(\boldsymbol{\rho}) + \frac{1}{2}\, \frac{\dif^2 E}{\dif \boldsymbol{\rho}^2} \cdot (\Delta \boldsymbol{\rho} \otimes \Delta \boldsymbol{\rho}) + \frac{1}{6}\, \frac{\dif^3 E}{\dif \boldsymbol{\rho}^3} \cdot (\Delta \boldsymbol{\rho})^{\otimes 3} + \frac{1}{24}\, \frac{\dif^4 E}{\dif \boldsymbol{\rho}^4} \cdot (\Delta \boldsymbol{\rho})^{\otimes 4}
+ \mathcal{O}(\|\Delta \boldsymbol{\rho}\|^5).
\end{aligned}
\end{equation}
We analyze the case where the total stiffness matrix (i.e. tangent stiffness matrix) 
\begin{equation} \label{eq: total stiffness}
\dfrac{\dif^2 E}{\dif \boldsymbol{\rho}^2}
= \Omega + \dfrac{\dif^2 E}{\dif \boldsymbol{f}^2}
\cdot \left(\dfrac{\dif \boldsymbol{f}}{\dif \boldsymbol{\rho}} \otimes \dfrac{\dif \boldsymbol{f}}{\dif \boldsymbol{\rho}} \right)
\end{equation}
is positive semidefinite. 

Under the small prestress assumption, if $\Delta \boldsymbol{\rho} \notin \mathcal{K}_b$, the second term in \eqref{eq: total stiffness} dominates and is strictly positive. Consequently, 
$\dif^2 E/\dif \boldsymbol{\rho}^2 \cdot (\Delta \boldsymbol{\rho} \otimes \Delta \boldsymbol{\rho}) > 0$. It therefore remains to consider $\Delta \boldsymbol{\rho} \in \mathcal{K}_b \setminus \{\boldsymbol{0}\}$, where the second term vanishes and the behavior is governed by the stress matrix $\Omega$. However, if $\Delta \boldsymbol{\rho}$ is only a `weak zero' of $\Omega$:
\begin{equation}
\Omega \cdot (\Delta \boldsymbol{\rho} \otimes \Delta \boldsymbol{\rho}) = 0
\quad \text{but} \quad
\Omega \cdot \Delta \boldsymbol{\rho} \neq \boldsymbol{0},
\end{equation}
then $\dif^2 E/\dif \boldsymbol{\rho}^2$ becomes indefinite; a proof is provided in Appendix~\ref{app: weak zero}.

We thus conclude the following characterization: 
$\dif^2 E/\dif \boldsymbol{\rho}^2$ is positive definite if and only if 

\begin{subequations}\label{eq:prestress-stability-test}
\begin{empheq}[left=\empheqlbrace]{align}
&\dim(\mathcal{K}_b) > 0. \label{eq:prestress-stability-test-1a}\\[6pt]
&\text{For any } \boldsymbol{\rho}' \in \mathcal{K}_b \setminus \{\boldsymbol{0}\}, \text{there exists a constant $\boldsymbol{\omega}$~} (\text{not a function of } \boldsymbol{\rho}') \text{ such that~} \notag\\
& \qquad\qquad\qquad\qquad\qquad\qquad \Omega \cdot (\boldsymbol{\rho}' \otimes \boldsymbol{\rho}') > 0. \label{prestress-stability-test-1b}
\end{empheq}
\end{subequations}

Consequently, $\boldsymbol{\rho}$ is a strict local minimum of $E$, and the energy grows tightly at order $s=2$. We say the polyhedral surface is \textit{(first-order) prestress stable} at $\boldsymbol{\rho}$. 

Let $\mathcal{K}_\omega = \mathrm{Null}(\Omega)$ denote the right null space of the stress matrix, then $\dif^2 E/\dif \boldsymbol{\rho}^2$ is positive semidefinite but not positive definite if and only if, 
\begin{subequations}\label{eq:semi-definite-test}
\begin{empheq}[left=\empheqlbrace]{align}
&\mathrm{dim}{(\mathcal{K}_b \cap \mathcal{K}_\omega)} > 0. \label{eq:semi-definite-test-1a}\\[6pt]
& \text{For every~} \boldsymbol{\rho}' \in \mathcal{K}_b, \text{either~}
    \Omega \cdot (\boldsymbol{\rho}' \otimes \boldsymbol{\rho}') > 0
    ~~\text{or}~~
    \boldsymbol{\rho}' \in \mathcal{K}_\omega \label{semi-definite-test-1b}
\end{empheq}
\end{subequations}

In this case, the null space of $\dif^2 E/\dif \boldsymbol{\rho}^2$ is $\mathcal{K}_b \cap \mathcal{K}_\omega$. If $\dif^2 E/\dif \boldsymbol{\rho}^2$ is not even positive semi-definite, we say the equilibrium configuration $\boldsymbol{\rho}$ is unstable, or the prestressed energy $E$ is a saddle point at $\boldsymbol{\rho}$, or the configuration $\boldsymbol{\rho}$ cannot be (higher-order) stabilized by any prestress. 

The aforementioned limitation of the higher-order test is as follows: if the third-order term of the Taylor expansion vanishes, and the fourth-order term is positive on $\mathcal{K}_b \cap \mathcal{K}_\omega$, one still cannot conclude that $E$ attains a strict local minimum at $\boldsymbol{\rho}$. Indeed, $\boldsymbol{\rho}$ may still be a saddle point. This is because the fourth-order Taylor approximation may grow at a rate slower than fourth order along certain trajectories, so the fourth-order term does not dominate the higher-order remainder. This phenomenon does not arise at second order, since positive definiteness provides a uniform quadratic lower bound in every direction and guarantees quadratic growth of the energy.

Rather than considering only an arbitrary small perturbation  $\Delta \boldsymbol{\rho}$, the `4th order test' analyzes the growth of
the energy along all quadratic trajectories of the form $\boldsymbol{\rho}(t) = \boldsymbol{\rho} + \boldsymbol{a}t + \boldsymbol{b}t^2$ where \(\boldsymbol{a}\) and \(\boldsymbol{b}\) are undetermined coefficients.
The key idea is that, for testing second-order prestress stability, it is
sufficient to identify the quadratic trajectory along which the energy grows
most slowly. This slowest-growth trajectory is precisely of the form $\boldsymbol{\rho}(t) = \boldsymbol{\rho} + \boldsymbol{\rho}'t + \boldsymbol{\rho}''t^2$, where $(\boldsymbol{\rho}', ~\boldsymbol{\rho}'')$ is a $(1, 2)\text{-}$flex. In this
way, the test reduces the stability question to the behavior of the energy
along the most critical second-order admissible direction. We summarize the
resulting test for \textit{second-order prestress stability} below.

\begin{subequations}\label{eq:second-order-prestress-stability-test}
\begin{empheq}[left=\empheqlbrace]{align}
&\text{There exists a self-stress~} \boldsymbol{\omega} \text{~such that} \dif^2 E/\dif \boldsymbol{\rho}^2 \text{~is positive semidefinite}. \label{eq:second-order-prestress-stability-test-a}\\[6pt]
    &\text{There exists } \text{(1, 2)-flexes~} (\boldsymbol{\rho}', \boldsymbol{\rho}''). \text{ More specifically, there exist } \boldsymbol{\rho}' \in \mathcal{K}_b \setminus \boldsymbol{0} \text{ and } \boldsymbol{\rho}''  \text{ such that:} \notag\\[3pt]
&\qquad\qquad\qquad\qquad\qquad\qquad
\dfrac{\dif \boldsymbol{f}}{\dif \boldsymbol{\rho}} \cdot \boldsymbol{\rho}''
+ \dfrac{\dif^2 \boldsymbol{f}}{\dif \boldsymbol{\rho}^2} 
\cdot \left(\boldsymbol{\rho}' \otimes \boldsymbol{\rho}'\right)
= \boldsymbol{0}.
\label{eq:second-order-prestress-stability-test-b}\\[6pt]
& \boldsymbol{\omega} \cdot \dfrac{ \dif^3 \boldsymbol{f}}{\dif \boldsymbol{\rho}^3} = \boldsymbol{0}. \label{eq:second-order-prestress-stability-test-c}\\[6pt]
&\text{For all~(1, 2)-flexes~} (\boldsymbol{\rho}', \boldsymbol{\rho}''), 
\notag\\[3pt]
&\qquad\qquad\qquad\qquad\qquad\qquad
3\,\Omega \cdot (\boldsymbol{\rho}'' \otimes \boldsymbol{\rho}'') + \Omega^{\mathrm{II}} \cdot (\boldsymbol{\rho}')^{\otimes 4}>0,
\label{eq:second-order-prestress-stability-test-d}
\end{empheq}
\end{subequations}

where, by analogy with the stress matrix $\Omega$, we refer to $\Omega^{\mathrm{II}}$ as the \emph{second-order stress tensor}.
\begin{equation*}
    \Omega^{\mathrm{II}} = \boldsymbol{\omega} \cdot \dfrac{\dif^4 \boldsymbol{f}}{\dif \boldsymbol{\rho}^4}. 
\end{equation*}
The higher-order terms in \eqref{eq:second-order-prestress-stability-test-d} comes from higher-order multivariate Taylor expansion, see Appendix \label{appendix: derivatives} \ref{app: Faa}.

\begin{figure}
	\noindent \centering	\includegraphics[width=1\linewidth]{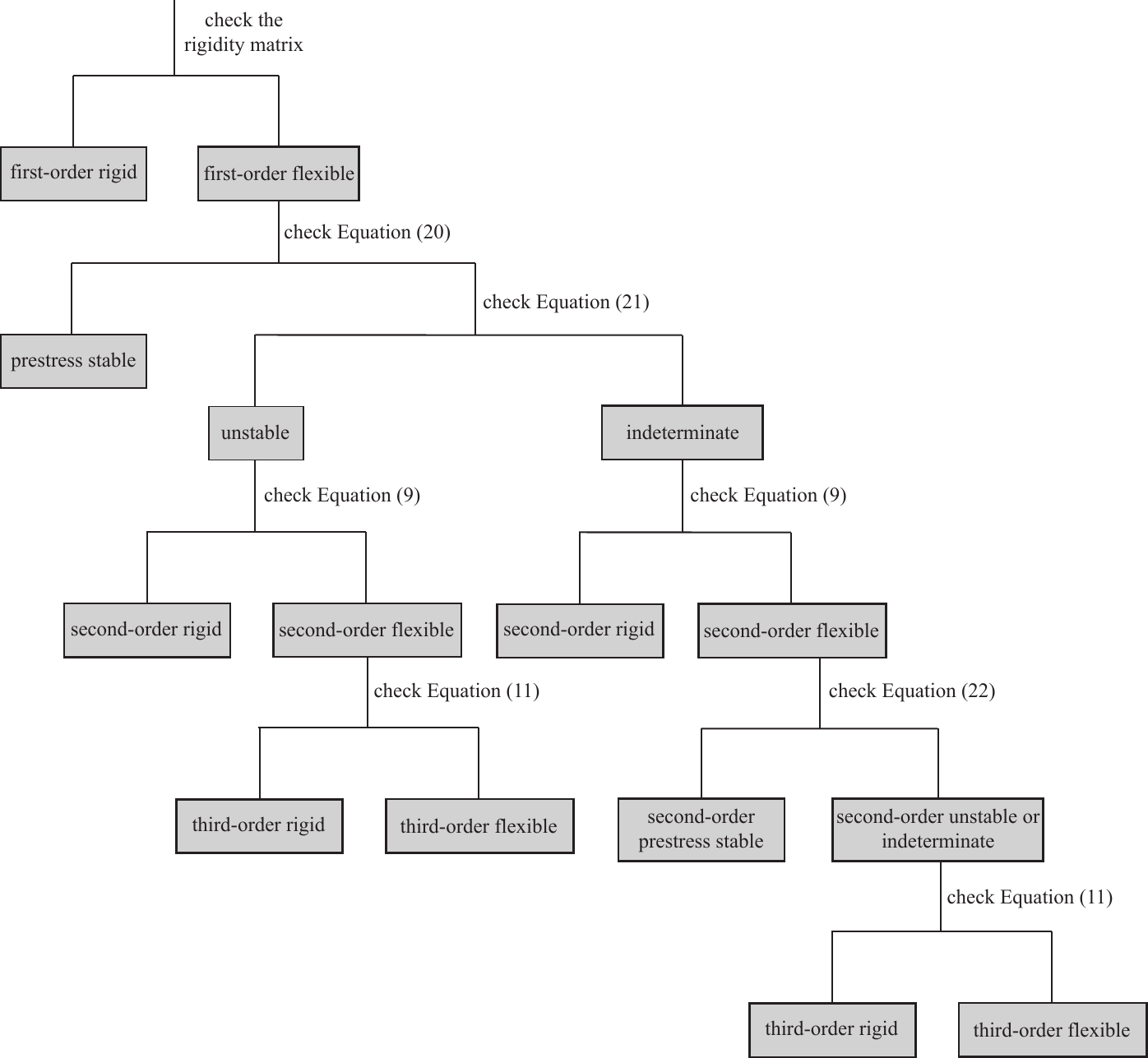}
	\caption{A flowchart illustrating the full procedure for rigidity and prestress stability tests, up to third-order rigidity.}
	\label{fig: flowchart}
\end{figure}

To conclude the rigidity and prestress stability tests, we present in Figure~\ref{fig: flowchart} a flowchart that illustrates the procedure for applying these tests to determine the exact order of rigidity of a polyhedral surface.

\clearpage

\clearpage

\section{Examples} \label{section: examples}

\subsection{Unstable and second-order rigid}

\begin{figure}
	\noindent \centering	\includegraphics[width=1\linewidth]{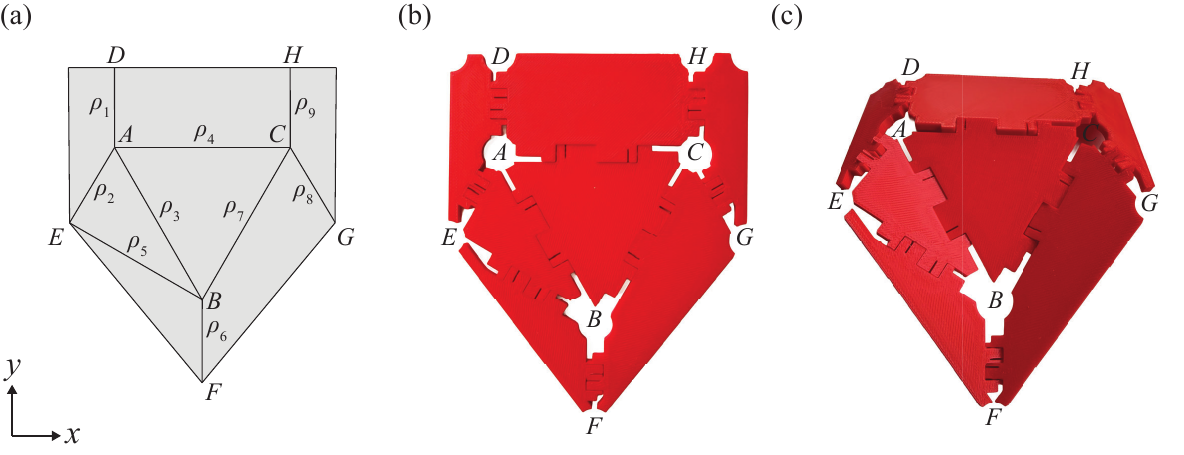}
	\caption{(a) A multi-vertex planar polyhedral surface that is unstable and second-order rigid in $\mathbb{R}^3$. (b) A 3D-printed model of the polyhedral surface and (c) its response to a small manual perturbation. The pip-and-socket hinges were independently reverse-engineered from a commercial Polydron element \citep{harvey_planar_1975}.}
	\label{fig: 2nd rigid and saddle}
\end{figure}

We will show that the planar polyhedral surface in Figure \ref{fig: 2nd rigid and saddle} is unstable and second-order rigid. This polyhedral surface is rigid in $\mathbb{R}^3$. Although each of the three degree-4 vertices is individually flexible, their folding motions cannot be combined into a compatible motion of the whole surface. The sector angles are chosen so that, on every possible folding branch, the motions prescribed by the three vertices are mutually incompatible. The coordinates for vertices $A$ to $H$ are in the following matrix
\begin{equation*}
\begin{bmatrix}
0 & 1 & 2 & 0 & -1/2 & 1 & 5/2 & 2 \\
0 & -\sqrt{3} & 0 & 1 & -\sqrt{3}/2 & -1-\sqrt{3} & -\sqrt{3}/2 & 1 \\
0 & 0 & 0 & 0 & 0 & 0 & 0 & 0
\end{bmatrix}.
\end{equation*}

The rigidity matrix is:
\begin{equation*}
\dfrac{\dif \boldsymbol{f}}{\dif \boldsymbol{\rho}} = \begin{bmatrix}
0 & -1/2 & 1/2 & 1 & 0 & 0 & 0 & 0 & 0 \\
1 & -\sqrt{3}/2 & -\sqrt{3}/2 & 0 & 0 & 0 & 0 & 0 & 0 \\
0 & 0 & 0 & 0 & 0 & 0 & 0 & 0 & 0 \\
0 & 0 & -1/2 & 0 & -\sqrt{3}/2 & 0 & 1/2 & 0 & 0 \\
0 & 0 & \sqrt{3}/2 & 0 & 1/2 & -1 & \sqrt{3}/2 & 0 & 0 \\
0 & 0 & 0 & 0 & 0 & 0 & 0 & 0 & 0 \\
0 & 0 & 0 & -1 & 0 & 0 & -1/2 & 1/2 & 0 \\
0 & 0 & 0 & 0 & 0 & 0 & -\sqrt{3}/2 & -\sqrt{3}/2 & 1 \\
0 & 0 & 0 & 0 & 0 & 0 & 0 & 0 & 0
\end{bmatrix}.
\end{equation*}

The right null space $\mathcal{K}_b$ of the rigidity matrix is three-dimensional. We denote by $K_b$ a matrix whose columns form a basis of $\mathcal{K}_b$:
\begin{equation*}
	K_b = \begin{bmatrix}
-5 & -2\sqrt{3} & -4\sqrt{3}/3 & -\sqrt{3}/3 & 2 & 0 & 2\sqrt{3}/3 & 0 & 1 \\
3\sqrt{3} & 4 & 2 & 1 & -\sqrt{3} & 0 & -1 & 1 & 0 \\
3 & \sqrt{3} & \sqrt{3} & 0 & -1 & 1 & 0 & 0 & 0
\end{bmatrix}^\mathrm{T}.
\end{equation*}
The self-stress $\boldsymbol{\omega}$ is:
\begin{equation*}
	\boldsymbol{\omega} = \begin{bmatrix}
		0 & 0 & w_1 & 0 & 0 & w_2 & 0 & 0 & w_3 
	\end{bmatrix}^{\mathrm{T}}, \quad w_1, ~w_2, ~w_3 \in \mathbb{R}.
\end{equation*}
The stress matrix $\boldsymbol{\Omega}= \boldsymbol{\omega} \cdot \dif^2 \boldsymbol{f}/\dif \boldsymbol{\rho}^2$ is:
\begin{equation*}
	\boldsymbol{\Omega} = \begin{bmatrix}
0 & 0 & 0 & 0 & 0 & 0 & 0 & 0 & 0 \\
0 & \frac{\sqrt{3}w_1}{4} & \frac{\sqrt{3}w_1}{4} & 0 & 0 & 0 & 0 & 0 & 0 \\
0 & \frac{\sqrt{3}w_1}{4} &
-\frac{\sqrt{3}w_1}{4}-\frac{\sqrt{3}w_2}{4} &
0 & -\frac{w_2}{4} & \frac{w_2}{2} &
-\frac{\sqrt{3}w_2}{4} & 0 & 0 \\
0 & 0 & 0 & 0 & 0 & 0 &
\frac{\sqrt{3}w_3}{2} & \frac{\sqrt{3}w_3}{2} & -w_3 \\
0 & 0 & -\frac{w_2}{4} & 0 &
-\frac{\sqrt{3}w_2}{4} & \frac{\sqrt{3}w_2}{2} &
-\frac{3w_2}{4} & 0 & 0 \\
0 & 0 & \frac{w_2}{2} & 0 &
\frac{\sqrt{3}w_2}{2} & 0 & 0 & 0 & 0 \\
0 & 0 & -\frac{\sqrt{3}w_2}{4} &
\frac{\sqrt{3}w_3}{2} & -\frac{3w_2}{4} & 0 &
\frac{\sqrt{3}(w_2+w_3)}{4} &
\frac{\sqrt{3}w_3}{4} & -\frac{w_3}{2} \\
0 & 0 & 0 & \frac{\sqrt{3}w_3}{2} & 0 & 0 &
\frac{\sqrt{3}w_3}{4} & -\frac{\sqrt{3}w_3}{4} &
\frac{w_3}{2} \\
0 & 0 & 0 & -w_3 & 0 & 0 &
-\frac{w_3}{2} & \frac{w_3}{2} & 0
\end{bmatrix}.
\end{equation*}
When checking prestress stability, the positive definiteness of $\Omega$ on
$\mathcal{K}_b$ is then equivalent to the positive definiteness
of the reduced matrix $K_b^T \Omega K_b$:
\begin{equation*}
K_b^T \Omega K_b =
\begin{bmatrix}
\frac{\sqrt{3}}{3}\left(17w_1-4w_2-w_3\right)
& -11w_1+2w_2+w_3
& \sqrt{3}\left(-3w_1+w_2\right) \\
-11w_1+2w_2+w_3
& \sqrt{3}\left(7w_1-w_2-\frac{w_3}{2}\right)
& 6w_1-\frac{3w_2}{2} \\
\sqrt{3}\left(-3w_1+w_2\right)
& 6w_1-\frac{3w_2}{2}
& \frac{\sqrt{3}}{2}\left(3w_1-w_2\right)
\end{bmatrix}.
\end{equation*}

Equivalently, we need to check whether all the 3 leading principal minors of $K_b^T \Omega K_b$ are positive:
\begin{subequations}\label{eq:sylvester-conditions}
\begin{empheq}[left=\empheqlbrace]{align}
& \dfrac{\sqrt{3}}{3}\left(17w_1-4w_2-w_3\right)>0, \label{eq:sylvester-d11} \\
& -2w_1^2+\dfrac{\left(-2w_2+13w_3\right)w_1}{2}-w_2w_3-\dfrac{w_3^2}{2}>0, \label{eq:sylvester-d21} \\
& -\dfrac{\sqrt{3}}{2}
\left[
\left(w_2+\dfrac{3w_3}{2}\right)w_1^2
+\dfrac{w_1}{2}\left(w_2^2-5w_2w_3+3w_3^2\right)
+\dfrac{w_2w_3}{2}\left(w_2-w_3\right)
\right]>0. \label{eq:sylvester-d31}
\end{empheq}
\end{subequations}

The polynomial inequalities \eqref{eq:sylvester-conditions} was tested using Maple’s RegularChains package. The commands SamplePoints and RealTriangularize employ regular-chain triangular decomposition, border-polynomial methods, and exact real-root isolation to determine the real semi-algebraic solution set. The computation shows that the system \eqref{eq:sylvester-conditions} has no solution. Therefore, there exists no
self-stress for which $K_b^T\Omega K_b$ is positive definite, and the
configuration $\boldsymbol{\rho}$ is not (first-order) prestress stable.

The next step is to determine whether the configuration $\boldsymbol{\rho}$
can be stabilized at higher order. We begin by analyzing the dimension of the
right null space $\mathcal{K}_\omega$, which is obtained from the rank of the
reduced stress matrix $\Omega K_b$:
\begin{equation*}
    \Omega K_b =
\begin{bmatrix}
0 & -5w_1/2 & -w_1/2 & 0 & -2\sqrt{3}w_2/3 & \sqrt{3}w_2/3 & -w_3/2 & w_3/2 & 0 \\
0 & 3\sqrt{3}w_1/2 & \sqrt{3}w_1/2 & 0 & w_2 & -w_2/2 & \sqrt{3}w_3/2 & 0 & 0 \\
0 & 3w_1/2 & 0 & 0 & \sqrt{3}w_2/2 & 0 & 0 & 0 & 0
\end{bmatrix}^\mathrm{T}.
\end{equation*}
The Gröbner basis of all the $2\times 2$ minors of $\Omega K_b$ is:
\begin{equation*}
    \left[
w_3^2,\;
w_2w_3,\;
w_1w_3,\;
w_2^2,\;
w_1w_2,\;
w_1^2
\right],
\end{equation*}
and they vanish
simultaneously only when
$w_1=w_2=w_3=0$.
Hence, there is no nonzero self-stress $\boldsymbol{\omega}$ for which
$\operatorname{rank}(\Omega K_b)=1$. Next, compute all $3\times 3$ minors of $\Omega K_b$, the Gröbner basis of all the $3\times 3$ minors of $\Omega K_b$ is:
\begin{equation*}
\begin{bmatrix}
w_2w_3^2 &
w_1w_3^2 &
w_2^2w_3 &
w_1w_2w_3 &
w_1^2w_3 &
w_1w_2^2 &
w_1^2w_2
\end{bmatrix}.
\end{equation*}
They vanish
simultaneously if and only if two of $w_1, ~w_2, ~w_3$ are zero, in which case $\operatorname{rank}(\Omega K_b)=2$. Let us do a case-by-case discussion:
\begin{enumerate} [label={[\arabic*]}]
\item Suppose that $w_1=w_2=0$. In this case, the right null space of
$\Omega K_b$ within $\mathcal{K}_b$ is spanned by
\[
K_\omega =
\begin{bmatrix}
0\\
0\\
1
\end{bmatrix}.
\]
A basis matrix for its orthogonal complement is given by
\[
K_\omega^\perp =
\begin{bmatrix}
1 & 0 \\
0 & 1 \\
0 & 0
\end{bmatrix}.
\]
It remains to check whether the restriction of $K_b^T\Omega K_b$ to this
orthogonal complement is positive definite. However,
\[
(K_\omega^\perp)^T K_b^T \Omega K_b K_\omega^\perp
=
\begin{bmatrix}
-\sqrt{3}w_3/3 & w_3 \\
w_3 & -\sqrt{3}w_3/2
\end{bmatrix}.
\]
This matrix is not positive definite, since its determinant $-w_3^2/2$ is not positive. Hence this case can be ruled out.
\item Suppose that $w_2=w_3=0$. In this case, the right null space of
$\Omega K_b$ within $\mathcal{K}_b$ is spanned by
\[
K_\omega =
\begin{bmatrix}
3/2 \\
\sqrt{3}/2 \\
1
\end{bmatrix}.
\]
A basis matrix for its orthogonal complement is given by
\[
K_\omega^\perp = \begin{bmatrix}
1 & 0 \\
0 & 1 \\
-3/2 & -\sqrt{3}/2
\end{bmatrix}\].
It remains to check whether the restriction of $K_b^T\Omega K_b$ to this
orthogonal complement is positive definite. However,
\[
(K_\omega^\perp)^T K_b^T \Omega K_b K_\omega^\perp
=
\begin{bmatrix}
433\sqrt{3}w_1/24 & -97w_1/8 \\
-97w_1/8 & 17\sqrt{3}w_1/8
\end{bmatrix}.
\]
This matrix is not positive definite, since its determinant $-32w_1^2$ is not positive. Hence this case can be ruled out.
\item Suppose that $w_3=w_1=0$. In this case, the right null space of
$\Omega K_b$ within $\mathcal{K}_b$ is spanned by
\[
K_\omega =
\begin{bmatrix}
\sqrt{3}/2 \\
1 \\
0
\end{bmatrix}.
\]
A basis matrix for its orthogonal complement is given by
\[
K_\omega^\perp =
\begin{bmatrix}
1 & 0 \\
-\sqrt{3}/2 & 0 \\
0 & 1
\end{bmatrix}.
\]
It remains to check whether the restriction of $K_b^T\Omega K_b$ to this
orthogonal complement is positive definite. However,
\[
(K_\omega^\perp)^T K_b^T \Omega K_b K_\omega^\perp
=
\begin{bmatrix}
-49\sqrt{3}w_2/12 & 7\sqrt{3}w_2/4 \\
7\sqrt{3}w_2/4 & -\sqrt{3}w_2/2
\end{bmatrix}.
\]
This matrix is not positive definite, since its determinant $-49w_2^2/16$ is
not positive. Hence this case can be ruled out.
\end{enumerate}

In conclusion, there exists no self-stress $\boldsymbol{\omega}$ for which the
polyhedral surface is (first-order) indeterminate. The polyhedral surface shown in Figure \ref{fig: 2nd rigid and saddle} cannot be stablized by any prestress and is a saddle point of any prestressed energy.

Next we check if the polyhedral surface is second-order rigid. Let 
\begin{equation*}
	\boldsymbol{\rho}' = K_b a = K_b \begin{bmatrix}
		a_1 & a_2 & a_3
	\end{bmatrix}^\mathrm{T}, \quad a_1, ~a_2, ~a_3 \in \mathbb{R}.
\end{equation*}
From \eqref{eq:second-order-rigidity-test-2b}, calculate the quadratic form:
\begin{equation}
\begin{aligned}
\boldsymbol{\Omega} \cdot \left(\boldsymbol{\rho}' \otimes \boldsymbol{\rho}'\right) ={}& \quad
\left(
a_1\sqrt{3}\left(17w_1-4w_2-w_3\right)/3
+a_2\left(-11w_1+2w_2+w_3\right)
+a_3\sqrt{3}\left(-3w_1+w_2\right)
\right)a_1 \notag\\
&+
\left(
a_1\left(-11w_1+2w_2+w_3\right)
+a_2\sqrt{3}\left(7w_1-w_2-w_3/2\right)
+a_3\left(6w_1-3w_2/2\right)
\right)a_2 \notag\\
&+
\left(
a_1\sqrt{3}\left(-3w_1+w_2\right)
+a_2\left(6w_1-3w_2/2\right)
+a_3\sqrt{3}\left(3w_1-w_2\right)/2
\right)a_3.
	\end{aligned}
\end{equation}
It is a linear function in $w_1, ~w_2, ~w_3$. In order for the quadratic form to be zero for all $\boldsymbol{\omega}$, we have:
\begin{equation} \label{eq:quadratic-form}
\begin{cases}
\left(34a_1^2-36a_1a_3+42a_2^2+9a_3^2\right)\sqrt{3}/6
-22\left(a_1-6a_3/11\right)a_2 = 0, \\[6pt]
\dfrac{\sqrt{3}}{2}
\left(-8a_1^2+12a_1a_3-6a_2^2-3a_3^2\right)\sqrt{3}/6
+4a_2\left(a_1-3a_3/4\right) = 0, \\[6pt]
\left(-2a_1^2-3a_2^2\right)\sqrt{3}/6+2a_1a_2 = 0.
\end{cases}
\end{equation}
The Gröbner basis of \eqref{eq:quadratic-form} is:
\begin{equation*}
\begin{bmatrix}
2\sqrt{3}a_2a_3+3a_2^2 &
-4\sqrt{3}a_1a_3+\sqrt{3}a_3^2+8a_1a_2+10a_2a_3 &
6\sqrt{3}a_2a_3+4a_1^2-12a_1a_3+3a_3^2 \\
7\sqrt{3}a_3^3+34a_2a_3^2 &
68a_1a_3^2+a_3^3 &
a_3^4
\end{bmatrix}.
\end{equation*}
The only solution is $a_1 = a_2 = a_3 = 0$, hence there is no extendable $\boldsymbol{\rho}'$, and the polyhedral surface is second-order rigid.

From the above analysis, we conclude that the polyhedral surface shown in Figure \ref{fig: 2nd rigid and saddle} is unstable and  second-order rigid.

\subsection{Unstable and third-order flexible}

\begin{figure}
	\noindent \centering	\includegraphics[width=1\linewidth]{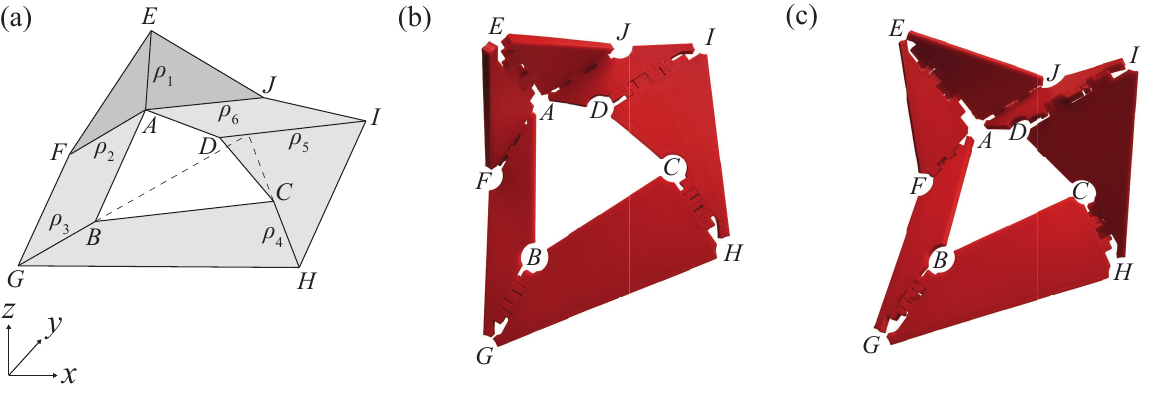}
	\caption{(a) A three-dimensional polyhedral surface with a hole that is rigid in $\mathbb{R}^3$, unstable, and  third-order flexible. (b) A 3D-printed model of the polyhedral surface and (c) a visualization of one direction of an extendable first-order flex. The pip-and-socket hinges were independently reverse-engineered from a commercial Polydron element \citep{harvey_planar_1975}.}
	\label{fig: 3rd order flexible and saddle}
\end{figure}

The three-dimensional polyhedral surface shown in Figure \ref{fig: 3rd order flexible and saddle} is rigid, unstable, and third-order flexible. This polyhedral surface is rigid in $\mathbb{R}^3$ because its six non-concurrent hinges impose six independent constraints. No additional folding motion can arise from any subset of the creases, since at most three creases are concurrent at any point. The higher-order shakiness arises from the separate concurrence of $\rho_1,\rho_2,\rho_6$ and $\rho_3,\rho_4,\rho_5$; however, these two concurrent hinge systems cannot generate a relative motion along the line joining their points of intersection.

The coordinates for vertices $A$ to $J$ are in the following matrix
\begin{equation*}
\begin{bmatrix}
0 & 0 & 3 & 3/2 & 0 & -1 & -1 & 4 & 4 & 2 \\
0 & -2\sqrt{3} & -\sqrt{3} & -\sqrt{3}/6 & 0 & -\sqrt{3} & -3\sqrt{3} & -2\sqrt{3} & 2\sqrt{3}/3 & 2\sqrt{3}/3 \\
0 & 0 & 0 & 0 & 2 & 0 & 0 & 0 & 0 & 0
\end{bmatrix}.
\end{equation*}

The rigidity matrix is:
\begin{equation*}
\dfrac{\dif \boldsymbol{f}}{\dif \boldsymbol{\rho}} = \begin{bmatrix}
0 & -1/2 & -1/2 & 1/2 & \sqrt{3}/2 & \sqrt{3}/2 \\
0 & -\sqrt{3}/2 & -\sqrt{3}/2 & -\sqrt{3}/2 & 1/2 & 1/2 \\
1 & 0 & 0 & 0 & 0 & 0 \\
0 & 0 & 0 & 0 & 0 & 0 \\
0 & 0 & 0 & 0 & 0 & 0 \\
0 & 0 & -\sqrt{3} & -\sqrt{3} & 1 & 0
\end{bmatrix}.
\end{equation*}

The right null space $\mathcal{K}_b$ of the rigidity matrix is two-dimensional. We denote by $K_b$ a matrix whose columns form a basis of $\mathcal{K}_b$:
\begin{equation*}
	K_b = \begin{bmatrix}
0 & \sqrt{3}/3 & \sqrt{3}/3 & -\sqrt{3}/3 & 0 & 1 \\
0 & 0 & 2\sqrt{3}/3 & -\sqrt{3}/3 & 1 & 0
\end{bmatrix}^\mathrm{T}.
\end{equation*}
The self-stress $\boldsymbol{\omega}$ is:
\begin{equation*}
	\boldsymbol{\omega} = \begin{bmatrix}
		0 & 0 &
	0 & w_1 & w_2 & 0 
	\end{bmatrix}^{\mathrm{T}}, \quad w_1, ~w_2 \in \mathbb{R}.
\end{equation*}
The stress matrix $\boldsymbol{\Omega}= \boldsymbol{\omega} \cdot \dif^2 \boldsymbol{f}/\dif \boldsymbol{\rho}^2$ is:
\begin{equation*}
	\boldsymbol{\Omega} = \begin{bmatrix}
0 & 0 & 0 & 0 & 0 & 0 \\
0 & 0 &
3w_1/2-\sqrt{3}w_2/2 &
3w_1/2-\sqrt{3}w_2/2 &
-\sqrt{3}w_1/2+w_2/2 &
0 \\
0 &
3w_1/2-\sqrt{3}w_2/2 &
3w_1/2-\sqrt{3}w_2/2 &
3w_1/2-\sqrt{3}w_2/2 &
-\sqrt{3}w_1/2+w_2/2 &
0 \\
0 &
3w_1/2-\sqrt{3}w_2/2 &
3w_1/2-\sqrt{3}w_2/2 &
3w_1/2+\sqrt{3}w_2/2 &
-\sqrt{3}w_1/2-w_2/2 &
0 \\
0 &
-\sqrt{3}w_1/2+w_2/2 &
-\sqrt{3}w_1/2+w_2/2 &
-\sqrt{3}w_1/2-w_2/2 &
w_1/2-\sqrt{3}w_2/2 &
0 \\
0 & 0 & 0 & 0 & 0 & 0
\end{bmatrix}.
\end{equation*}
When checking prestress stability, the positive definiteness of $\Omega$ on
$\mathcal{K}_b$ is then equivalent to the positive definiteness
of the reduced matrix $K_b^T \Omega K_b$:
\begin{equation*}
K_b^T \Omega K_b =
\begin{bmatrix}
\sqrt{3}w_2/3 & 2\sqrt{3}w_2/3 \\
2\sqrt{3}w_2/3 & 2\sqrt{3}w_2/3
\end{bmatrix}.
\end{equation*}

There is no self-stress for which $K_b^T\Omega K_b$ is positive definite, since its determinant is $-2w_2^2/3\leq 0$. Hence, the configuration $\boldsymbol{\rho}$ is not (first-order) prestress stable.

The next step is to determine whether the configuration $\boldsymbol{\rho}$
can be stabilized at higher order. We begin by analyzing the dimension of the
right null space of the stress matrix $\mathcal{K}_\omega$, which is obtained from the rank of the
reduced stress matrix $\Omega K_b$:
\begin{equation*}
    \Omega K_b =
\begin{bmatrix}
0 & 0 & (\sqrt{3}w_1-w_2)/2 & (\sqrt{3}w_1-3w_2)/2 & (-w_1+\sqrt{3}w_2)/2 & 0 \\
0 & 0 & 0 & -2w_2 & 0 & 0
\end{bmatrix}^\mathrm{T}.
\end{equation*}
The Gröbner basis of all the $2\times 2$ minors of $\Omega K_b$ is:
\begin{equation*}
    \left[
w_2^2,\;
w_1w_2
\right],
\end{equation*}
and they vanish
simultaneously when
$w_2=0$. In this case, the right null space of
$\Omega K_b$ within $\mathcal{K}_b$ is spanned by
\[
K_\omega =
\begin{bmatrix}
0\\
1
\end{bmatrix}.
\]
A basis matrix for its orthogonal complement is given by
\[
K_\omega^\perp =
\begin{bmatrix}
1 \\
0 \\
\end{bmatrix}.
\]
The restriction of $K_b^T\Omega K_b$ to this
orthogonal complement 
$(K_\omega^\perp)^T K_b^T \Omega K_b K_\omega^\perp
= 0$. It is not positive, hence there exists no self-stress $\boldsymbol{\omega}$ for which the
polyhedral surface is (first-order) indeterminate. The polyhedral surface shown in Figure \ref{fig: 3rd order flexible and saddle} is therefore unstable.

Next we check if the polyhedral surface is second-order rigid. Let 
\begin{equation*}
	\boldsymbol{\rho}' = K_b a = K_b \begin{bmatrix}
		a_1 & a_2 
	\end{bmatrix}^\mathrm{T}, \quad a_1,~a_2 \in \mathbb{R}. 
\end{equation*}
From \eqref{eq:second-order-rigidity-test-2b}, calculate the quadratic form:
\begin{equation}
\begin{aligned}
\boldsymbol{\Omega} \cdot \left(\boldsymbol{\rho}' \otimes \boldsymbol{\rho}'\right) =\left(a_1w_2\sqrt{3}/3+2a_2w_2\sqrt{3}/3\right)a_1
+\left(2a_1w_2\sqrt{3}/3+2a_2w_2\sqrt{3}/3\right)a_2.
	\end{aligned}
\end{equation}
It is a linear function in $w_2$. In order for the quadratic form to be zero for all $\boldsymbol{\omega}$, we have:
\begin{equation} 
\sqrt{3}\left(a_1^2+4a_1a_2+2a_2^2\right)/3 = 0.
\end{equation}
The nonzero solution, i.e. an extendable $\boldsymbol{\rho}'$ is given by:
\begin{equation} \label{eq: extendable rho'}
 a_1=(-2 \pm \sqrt{2})a_2.
\end{equation}
We can then compute $\boldsymbol{\rho}''$ by solving \eqref{eq:second-order-rigidity-test-1b}; the solution set forms a two-dimensional linear space:
\begin{equation*}
    \boldsymbol{\rho}' = K_b \begin{bmatrix}
		-2 \pm \sqrt{2} & 1 
	\end{bmatrix}^\mathrm{T} = \begin{bmatrix}
0 &
\sqrt{3}\left(-2 \pm \sqrt{2}\right)/3 &
\pm \sqrt{6}/3 &
\sqrt{3}\left(1 \pm \sqrt{2}\right)/3 &
1 &
-2 \pm \sqrt{2}
\end{bmatrix}^\mathrm{T},
\end{equation*}
\begin{equation*}	
\frac{\dif^2 \boldsymbol{f}}{\dif \boldsymbol{\rho}^2}
 \cdot \left(\boldsymbol{\rho}' \otimes \boldsymbol{\rho}'\right)
= \begin{bmatrix}
0 & 0 & \sqrt{3}\left(1 \mp 2\sqrt{2}/3\right) & 0 & 0 & 0
\end{bmatrix}^\mathrm{T},
\end{equation*}
\begin{equation*}
    \boldsymbol{\rho}'' = \begin{bmatrix}
\sqrt{3}\left(1 \mp 2\sqrt{2}/3\right) &
-b_1-2b_2 &
b_1 &
b_2 &
\sqrt{3}b_1+\sqrt{3}b_2 &
-\sqrt{3}b_1-2\sqrt{3}b_2
\end{bmatrix}^\mathrm{T}, \quad b_1, ~b_2 \in \mathbb{R}.
\end{equation*}

Since $\mathrm{dim}(\mathcal{K}_b)>1$, we cannot determine whether or not the polyhedral surface is exactly third-order rigid using \eqref{eq:third-order-rigidity-test-2}. However, we can test if there is a $(1, 3)$-flex. From \eqref{eq:third-order-rigidity-test-1c}, we then compute,
\begin{equation*}
-\dfrac{\dif^3 \boldsymbol{f}}{\dif \boldsymbol{\rho}^3} 
\cdot (\boldsymbol{\rho}')^{\otimes 3} = 
\begin{bmatrix}
\sqrt{3}\left(7 \mp 5\sqrt{2}\right)/2 &
\left(-35 \pm 25\sqrt{2}\right)/6 &
0 &
0 &
0 &
5/3 \mp \sqrt{2}
\end{bmatrix}^\mathrm{T},
\end{equation*}
\begin{equation*}
-3\,\dfrac{\dif^2 \boldsymbol{f}}{\dif \boldsymbol{\rho}^2} 
\cdot \left(\boldsymbol{\rho}'' \otimes \boldsymbol{\rho}'\right) = \begin{bmatrix}
0 &
0 &
3b_2\left(1 \mp \sqrt{2}\right) &
0 &
6b_1+6b_2 \mp 3\sqrt{2}b_1 &
0
\end{bmatrix}^\mathrm{T}.
\end{equation*}

The right-hand side vector of \eqref{eq:third-order-rigidity-test-1c} is the sum of the above two vectors:
\begin{equation*}
    \begin{bmatrix}
\sqrt{3}\left(7 \mp 5\sqrt{2}\right)/2 &
\left(-35 \pm 25\sqrt{2}\right)/6 &
3b_2\left(1 \mp \sqrt{2}\right) &
0 &
6b_1+6b_2 \mp 3\sqrt{2}b_1 &
5/3 \mp \sqrt{2}
\end{bmatrix}^\mathrm{T}.
\end{equation*}
For \eqref{eq:third-order-rigidity-test-1c} to be solvable, its right-hand side must be orthogonal to every self-stress $\boldsymbol{\omega}$. This compatibility condition yields
$b_1=(-2\mp\sqrt{2})b_2$.
Therefore, choosing $\boldsymbol{\rho}'$ as in \eqref{eq: extendable rho'} and $\boldsymbol{\rho}''$ such that $b_1=(-2\mp\sqrt{2})b_2$ produces a $(1,3)$-flex. Hence, the polyhedral surface shown in Figure \ref{fig: 3rd order flexible and saddle} is third-order flexible.

From the above analysis, we conclude that the polyhedral surface shown in Figure \ref{fig: 3rd order flexible and saddle} is rigid, unstable, and third-order flexible.

\subsection{Indeterminate, not second-order prestress stable, and third-order flexible}

\begin{figure}
	\noindent \centering	\includegraphics[width=1\linewidth]{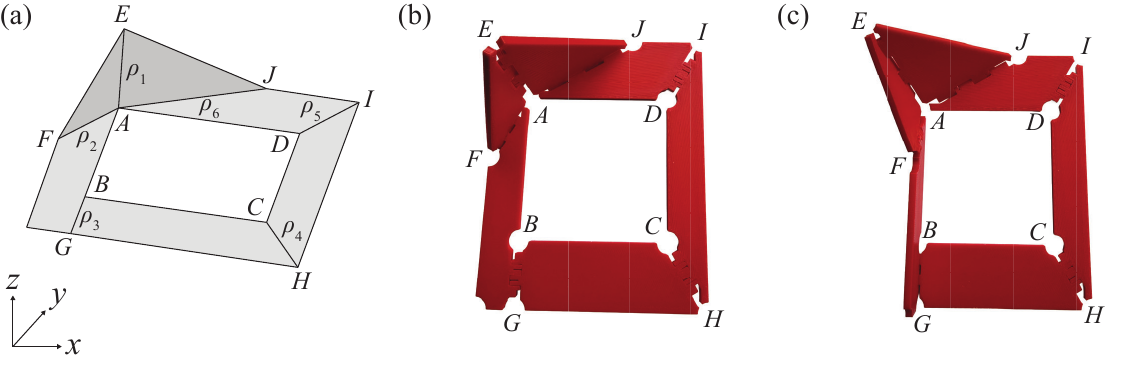}
	\caption{(a) A three-dimensional polyhedral surface with a hole. Rigidity test shows that it is indeterminate, not second-order prestress stable, and third-order flexible. (b) A 3D-printed model of the polyhedral surface. The pip-and-socket hinges were independently reverse-engineered from a commercial Polydron element \citep{harvey_planar_1975}. (c) A folded state of (a). The folding motion is analogous to that of a degree-4 rigid-origami vertex centred at $A$, while the panels $BGHC$, $CHID$, and $ADIJ$ remain fixed. }
	\label{fig: 3rd order flexible and indeterminate}
\end{figure}

The three-dimensional polyhedral surface in Figure \ref{fig: 3rd order flexible and indeterminate} is flexible in $\mathbb{R}^3$. The folding motion is analogous to that of a degree-4 rigid-origami vertex centred at $A$, while the panels $BGHC$, $CHID$, and $ADIJ$ remain fixed. The rigidity test shows that this polyhedral surface is  indeterminate, not second-order prestress stable, and third-order flexible.  The coordinates for vertices $A$ to $J$ are in the following matrix
\begin{equation*}
\begin{bmatrix}
0 & 0 & 2 & 2 & 0 & -1/2 & 0 & 5/2 & 5/2 & 3/2 \\
0 & -2 & -2 & 0 & 0 & -\sqrt{3}/2 & -\sqrt{3}/2 - 2 & -\sqrt{3}/2 - 2 & \sqrt{3}/2 & \sqrt{3}/2 \\
0 & 0 & 0 & 0 & 1 & 0 & 0 & 0 & 0 & 0
\end{bmatrix}
\end{equation*}.

The rigidity matrix is:
\begin{equation*}
\dfrac{\dif \boldsymbol{f}}{\dif \boldsymbol{\rho}} = \begin{bmatrix}
0 & -1/2 & 0 & 1/2 & 1/2 & \sqrt{3}/2 \\
0 & -\sqrt{3}/2 & -1 & -\sqrt{3}/2 & \sqrt{3}/2 & 1/2 \\
1 & 0 & 0 & 0 & 0 & 0 \\
0 & 0 & 0 & 0 & 0 & 0 \\
0 & 0 & 0 & 0 & 0 & 0 \\
0 & 0 & 0 & -\sqrt{3}+1 & \sqrt{3} & 0
\end{bmatrix}.
\end{equation*}

The right null space $\mathcal{K}_b$ of the rigidity matrix is two-dimensional. We denote by $K_b$ a matrix whose columns form a basis of $\mathcal{K}_b$:
\begin{equation*}
	K_b = \begin{bmatrix} 0 & \sqrt{3} & -1 & 0 & 0 & 1 \\ 0 & 5/2+\sqrt{3}/2 & -3/2-3\sqrt{3}/2 & 3/2+\sqrt{3}/2 & 1 & 0 \end{bmatrix}^\mathrm{T}.
\end{equation*}
The self-stress $\boldsymbol{\omega}$ is:
\begin{equation*}
	\boldsymbol{\omega} = \begin{bmatrix}
		0 & 0 & 0 & w_1 & w_2 & 0
	\end{bmatrix}^{\mathrm{T}}, \quad w_1, ~w_2 \in \mathbb{R}.
\end{equation*}
The stress matrix $\boldsymbol{\Omega}= \boldsymbol{\omega} \cdot \dif^2 \boldsymbol{f}/\dif \boldsymbol{\rho}^2$ is:
\begin{equation*}
	\boldsymbol{\Omega} = 
\begin{bmatrix}
0 & 0 & 0 & 0 & 0 & 0 \\
0 & 0 & 0 & ((3-\sqrt{3})w_1+(1-\sqrt{3})w_2)/2 & (-3w_1+\sqrt{3}w_2)/2 & 0 \\
0 & 0 & 0 & (\sqrt{3}-1)w_1 & -\sqrt{3}w_1 & 0 \\
0 & ((3-\sqrt{3})w_1+(1-\sqrt{3})w_2)/2 & (\sqrt{3}-1)w_1 & ((3-\sqrt{3})w_1+(\sqrt{3}-1)w_2)/2 & (-3w_1-\sqrt{3}w_2)/2 & 0 \\
0 & (-3w_1+\sqrt{3}w_2)/2 & -\sqrt{3}w_1 & (-3w_1-\sqrt{3}w_2)/2 & (3w_1-\sqrt{3}w_2)/2 & 0 \\
0 & 0 & 0 & 0 & 0 & 0
\end{bmatrix}.
\end{equation*}
When checking prestress stability, the positive definiteness of $\Omega$ on
$\mathcal{K}_b$ is then equivalent to the positive definiteness
of the reduced matrix $K_b^T \Omega K_b$:
\begin{equation*}
K_b^T \Omega K_b =
\begin{bmatrix} 0 & 0 \\ 0 & (\sqrt{3}w_2-3w_1-6w_2)/(2\sqrt{3}-2) \end{bmatrix}.
\end{equation*}

As the two leading principal minors of $K_b^T \Omega K_b$ are all zero, we conclude that there exists no
self-stress for which $K_b^T\Omega K_b$ is positive definite, and the
configuration $\boldsymbol{\rho}$ is not (first-order) prestress stable.

The next step is to determine whether the configuration $\boldsymbol{\rho}$
can be stabilized at higher order. We begin by analyzing the dimension of the
right null space $\mathcal{K}_\omega$, which is obtained from the rank of the
reduced stress matrix $\Omega K_b$:
\begin{equation*}
    \Omega K_b =
\begin{bmatrix}
0 & 0 & 0 & ((\sqrt{3}-1)w_1+(\sqrt{3}-3)w_2)/2 & (-\sqrt{3}w_1+3w_2)/2 & 0 \\
0 & 0 & 0 & (-\sqrt{3}w_1+(1-2\sqrt{3})w_2)/2 & 0 & 0
\end{bmatrix}^\mathrm{T}.
\end{equation*}
There is only one nonzero $2\times 2$ minor
\begin{equation*}
    \left[
3w_1^2 + (6-4\sqrt{3})w_1w_2 + (3-6\sqrt{3})w_2^2
\right],
\end{equation*}
and it vanishes when
$w_1 = \left(\sqrt{3}/3 - 2\right)w_2$ or $w_1 = \sqrt{3}w_2$. Let us do a case-by-case discussion:
\begin{enumerate} [label={[\arabic*]}]
\item Suppose that $w_1 = \left(\sqrt{3}/3 - 2\right)w_2$. In this case, the right null space of
$\Omega K_b$ within $\mathcal{K}_b$ is spanned by
\[
K_\omega =
\begin{bmatrix}
0\\
1
\end{bmatrix}.
\]
A basis matrix for its orthogonal complement is given by
\[
K_\omega^\perp =
\begin{bmatrix}
1 \\
0 \\
\end{bmatrix}.
\]
It remains to check whether the restriction of $K_b^T\Omega K_b$ to this
orthogonal complement is positive definite. However,
\[
(K_\omega^\perp)^T K_b^T \Omega K_b K_\omega^\perp
= 0.
\]
Hence this case can be ruled out.
\item Suppose that $w_1 = \sqrt{3}w_2$. In this case, the right null space of
$\Omega K_b$ within $\mathcal{K}_b$ is spanned by
\[
K_\omega =
\begin{bmatrix}
1 \\
0
\end{bmatrix}.
\]
A basis matrix for its orthogonal complement is given by
\[
K_\omega^\perp =
\begin{bmatrix}
0 \\
1 \\
\end{bmatrix}.
\]
It remains to check whether the restriction of $K_b^T\Omega K_b$ to this
orthogonal complement is positive definite. As,
\[
(K_\omega^\perp)^T K_b^T \Omega K_b K_\omega^\perp
=
-(3+2\sqrt{3})w_2,
\]
we conclude that the choice $w_2<0$ will make the polyhedral surface indeterminate, i.e., possible to be stabilized at higher-order.
\end{enumerate}

Next we check if the polyhedral surface is second-order rigid. Let 
\begin{equation*}
	\boldsymbol{\rho}' = K_b a = K_b \begin{bmatrix}
		a_1 & a_2 
	\end{bmatrix}^\mathrm{T}, \quad a_1,~a_2 \in \mathbb{R}.
\end{equation*}
Calculate the quadratic form:
\begin{equation}
\begin{aligned}
\boldsymbol{\Omega} \cdot \left(\boldsymbol{\rho}' \otimes \boldsymbol{\rho}'\right) =-a_2^2\left(3(\sqrt{3}+1)w_1+(3+5\sqrt{3})w_2\right)/4.
	\end{aligned}
\end{equation}
For the quadratic form to vanish for every self-stress $\boldsymbol{\omega}$, we must have $a_2=0$. Hence the extendable  $\boldsymbol{\rho}'$ is obtained by imposing $a_2=0$. We can then compute $\boldsymbol{\rho}''$ by solving \eqref{eq:second-order-rigidity-test-1b}; the solution set forms a two-dimensional linear space:
\begin{equation*}
    \boldsymbol{\rho}' = K_b \begin{bmatrix}
		1 & 0 
	\end{bmatrix}^\mathrm{T} = \begin{bmatrix} 0 & \sqrt{3} & -1 & 0 & 0 & 1  \end{bmatrix}^\mathrm{T},
\end{equation*}
\begin{equation*}	
\frac{\dif^2 \boldsymbol{f}}{\dif \boldsymbol{\rho}^2}
 \cdot \left(\boldsymbol{\rho}' \otimes \boldsymbol{\rho}'\right)
= \begin{bmatrix} 0 & 0 & \sqrt{3}/2 & 0 & 0 & 0 \end{bmatrix}^\mathrm{T},
\end{equation*}
\begin{equation*}
    \boldsymbol{\rho}'' = \begin{bmatrix}
\sqrt{3}/2 &
\left(2-\sqrt{3}/3\right)b_1+\sqrt{3}b_2 &
-\sqrt{3}b_1-b_2 &
b_1 &
\left(1-\sqrt{3}/3\right)b_1 &
b_2
\end{bmatrix}, \quad b_1, ~b_2 \in \mathbb{R}.
\end{equation*}

Now we are ready to test the second-order prestress stability. Recall that the self-stress that makes the polyhedral surface indeterminate is:
\begin{equation*}
	\boldsymbol{\omega}_{\mathrm{i}} = \begin{bmatrix}
		0 & 0 & 0 &  \sqrt{3}w_2 & w_2 & 0
	\end{bmatrix}^{\mathrm{T}}, \quad w_2 < 0.
\end{equation*}
First of all we have:
\begin{equation*}
   {\boldsymbol{\omega}}_{\mathrm{i}}  \cdot \dfrac{ \dif^3 \boldsymbol{f}}{\dif \boldsymbol{\rho}^3} = \boldsymbol{0}. 
\end{equation*}
Second, consider  \eqref{eq:second-order-prestress-stability-test-d} for all $(1, 2)$-flexes $(\boldsymbol{\rho}', ~\boldsymbol{\rho}'')$:  \begin{equation*}
    3\,\Omega_\mathrm{i} \cdot (\boldsymbol{\rho}'' \otimes \boldsymbol{\rho}'')  = -2\sqrt{3}w_2b_1^2,
\end{equation*}
\begin{equation*}
    \dfrac{ \dif^4 \boldsymbol{f}}{\dif \boldsymbol{\rho}^4}  \cdot (\boldsymbol{\rho}')^{\otimes 4} = \begin{bmatrix} 0 & 0 & 11\sqrt{3}/4 & 0 & 0 & 0 \end{bmatrix}^\mathrm{T},
\end{equation*}
\begin{equation*}
\Omega_{\mathrm{i}}^{\mathrm{II}} \cdot (\boldsymbol{\rho}')^{\otimes 4} =  \boldsymbol{\omega}_{\mathrm{i}} \cdot \dfrac{ \dif^4 \boldsymbol{f}}{\dif \boldsymbol{\rho}^4}  \cdot (\boldsymbol{\rho}')^{\otimes 4} = 0,
\end{equation*}
hence,
\begin{equation} \label{eq:second-order-determinant}
    3\,\Omega_{\mathrm{i}} \cdot (\boldsymbol{\rho}'' \otimes \boldsymbol{\rho}'') + \Omega_{\mathrm{i}}^{\mathrm{II}} \cdot (\boldsymbol{\rho}')^{\otimes 4} = -2\sqrt{3}w_2b_1^2.
\end{equation}

We cannot let \eqref{eq:second-order-determinant} be positive for all $(1, 2)$-flexes $(\boldsymbol{\rho}', \boldsymbol{\rho}'')$ since it equals to zero if $b_1 = 0$. Consequently, the polyhedral surface shown in Figure \ref{fig: 3rd order flexible and indeterminate} is not second-order prestress stable.

It remains to see if the polyhedral surface is third-order rigid. Since $\mathrm{dim}(\mathcal{K}_b)>1$, we cannot determine whether or not the polyhedral surface is exactly third-order rigid using \eqref{eq:third-order-rigidity-test-2}. However, we can test if there is a $(1, 3)$-flex. From \eqref{eq:third-order-rigidity-test-1c}, we then compute,
\begin{equation*}
-\dfrac{\dif^3 \boldsymbol{f}}{\dif \boldsymbol{\rho}^3} 
\cdot (\boldsymbol{\rho}')^{\otimes 3}
= -\dfrac{\dif^3 \boldsymbol{f}}{\dif \boldsymbol{\rho}^3} 
\cdot \left(K_b \begin{bmatrix}
1 \\ 0 
\end{bmatrix}\right)^{\otimes 3} = 
\begin{bmatrix}
-\sqrt{3} & 3/2 & 0 & 0 & 0 & 0
\end{bmatrix}^\mathrm{T},
\end{equation*}
\begin{equation*}
-3\,\dfrac{\dif^2 \boldsymbol{f}}{\dif \boldsymbol{\rho}^2} 
\cdot \left(\boldsymbol{\rho}'' \otimes \boldsymbol{\rho}'\right) =
\begin{bmatrix}
0 & 0 & (6-\sqrt{3})b_1/2 + 3\sqrt{3}b_2/2 & 0 & 0 & 0
\end{bmatrix}.
\end{equation*}

\eqref{eq:third-order-rigidity-test-1c} admits a solution, since its right-hand side is orthogonal to every self-stress $\boldsymbol{\omega}$ for any $b_1, ~b_2$. Thus, every such pair $(\boldsymbol{\rho}',\boldsymbol{\rho}'')$ can be extended to a $(1,3)$-flex. It follows that the polyhedral surface described in Figure \ref{fig: 3rd order flexible and indeterminate} is third-order flexible.

Combining the above results, we conclude that the flexible polyhedral surface shown in Figure \ref{fig: 3rd order flexible and indeterminate} is indeterminate, not second-order prestress stable, and third-order flexible from the rigidity test.

\section{Discussion}

\subsection{Critical order for flexibility}

As mentioned in Section~\ref{section: introduction}, for geometric constraint systems there exists a sufficiently large integer $n^*$ such that $n^*$-th order flexibility is equivalent to actual flexibility. For polyhedral surfaces, the case $n^* = 1$ corresponds to regular realizations of generically flexible polyhedral surfaces \citep{he_rigid_2019}, while $n^* = 2$ corresponds to certain generically rigid polyhedral surfaces satisfying additional geometric conditions. An example is given by developable and flat-foldable quad-mesh origami \citep{demaine_zero-area_2016}. Note that, the energy-based definition used in this article is consistent with the kinematic definition up to second-order rigidity in previous literature. We conjecture that $n^* = 2$ also holds for discrete Voss surfaces, or V-hedra \citep{voss_diejenigen_1888}, since their folding angles satisfy algebraic relations closely analogous to those arising in developable and flat-foldable quad-mesh origami. It would therefore be particularly promising to investigate the value of $n^*$ within known classes of quad-mesh rigid origami. This question is tractable because such quad-mesh origami have only one degree of freedom, i.e., $\dim (\mathcal{K}_b) = 1$, so the existence of a $(1,n)$-flex is sufficient to establish $n$-th order flexibility.

\subsection{Applicability to general geometric constraint systems}

The computational framework developed in Section~\ref{section: second and third} and summarized in Figure~\ref{fig: flowchart}, as illustrated by the examples in Section~\ref{section: examples}, is directly applicable to geometric constraint systems involving other types of constraints. Indeed, the formulation assumes only a general energy function $E(\boldsymbol{f})$, defined in terms of a tuple of measurements $\boldsymbol{f}(\boldsymbol{\rho})$, where $\boldsymbol{\rho}$ denotes a set of generalized coordinates. The principal physical assumption is that the material stiffness matrix 
$
\dif^2 E / \dif \boldsymbol{f}^2
$
is positive definite. This assumption is appropriate for ordinary passive materials/physical systems at a stable state, as discussed in Section~\ref{section: introduction}. However, some materials and mechanical systems possess neutral deformation modes along which the material stiffness vanishes, so that the corresponding stiffness matrix is only positive semidefinite. Extending the theory of higher-order prestress stability to such systems would therefore constitute a natural and potentially important direction for future research.

\subsection{Challenges in extending the theory}

There are three major challenges in extending the present rigidity framework. The first is to characterize second-order indeterminacy. This notion is a natural extension of first-order indeterminacy and describes systems that are not stabilized at second order but may still be stabilized by third- or higher-order effects. For first-order indeterminacy, the existing theory provides the condition~\eqref{eq:semi-definite-test}, formulated in terms of the right null space $\mathcal{K}_b$ of the rigidity matrix and the stress matrix $\Omega$, under which the total stiffness matrix
$ \dif^2 E / \dif \boldsymbol{\rho}^2 $
in~\eqref{eq: total stiffness} is positive semidefinite under the small-prestress assumption. In the second-order indeterminate case, however, the relevant energy expansion contains an additional quartic contribution $
\Omega^{\mathrm{II}} \cdot (\boldsymbol{\rho}')^{\otimes 4}$.
An analogous criterion would therefore need to ensure the nonnegativity of the quartic form associated with
$
\dif^4 E/\dif \boldsymbol{\rho}^4
$
on the space of $(1, 2)$-flexes. At present, it is unclear how to derive a tractable condition analogous to~\eqref{eq:semi-definite-test}.

The second challenge is to determine third- and higher-order rigidity when $\dim(\mathcal{K}_b)>1$. To establish $n$-th order rigidity in this setting, one must rule out the existence of a $(j,nj)$-flex for every positive integer $j$. In general, it is not sufficient to exclude only $(1,n)$-flexes. The problem therefore involves, in principle, an infinite family of higher-order extension tests, which makes the determination of higher-order rigidity substantially more difficult. Nevertheless, for special frameworks, it may be possible to treat all values of $j$ simultaneously. A bar-joint framework example is provided by \citet{stachel_proposal_2007} (see also \citet[Example~3]{nawratil_global_2025}), where, for every positive integer $j$, a $(j,3j-1)$-flex exists but no $(j,3j)$-flex exists.

The third challenge concerns fractional orders of rigidity and prestress stability, which may arise beyond third order. Establishing the existence of such fractional orders in different types of geometric constraint systems would provide a promising direction for further research. For example, a bar-joint framework presented in \citet[Example~7]{nawratil_global_2025} appears to exhibit rigidity of order $7/2$.

\subsection{Feasibility of computation}

The computational workflow described in Figure~\ref{fig: flowchart} involves three main types of symbolic calculations.

\begin{enumerate} [label={[\arabic*]}]
\item \textbf{Polynomial inequality problems.}
In the test for prestress stability, one seeks a self-stress for which the reduced stress matrix is positive definite on the infinitesimal flex space $\mathcal{K}_b$. Let
$
d=\dim(\mathcal{K}_b)
$
and let $s$ denote the dimension of the left null space of the rigidity matrix, equivalently, the dimension of the self-stress space. The reduced stress matrix is a symmetric $d\times d$ matrix whose entries depend linearly on the $s$ coefficients of the self-stress vector $\boldsymbol{\omega}$. By Sylvester's criterion, positive definiteness can be tested by requiring its $d$ leading principal minors to be positive. These minors are polynomials in the self-stress coefficients, with degrees bounded respectively by $1,\ldots,d$. The prestress-stability test therefore reduces to the feasibility of a system of polynomial inequalities in $s$ variables.

The test for second-order prestress stability leads to a more difficult quantified polynomial inequality. In this case, one seeks an indeterminate self-stress $\boldsymbol{\omega}_{\mathrm{i}}$ such that the relevant higher-order energy term is positive for every admissible $(1,2)$-flex. After parametrizing the space of $(1,2)$-flexes $\mathcal{B}$ with free parametres $\boldsymbol{b}$, this condition typically takes the form
$$
\exists\,\boldsymbol{\omega}_{\mathrm{i}}
\quad
\forall\,\boldsymbol{b}\in\mathcal{B}:
\quad
p(\boldsymbol{\omega}_{\mathrm{i}},\boldsymbol{b})>0,
$$
The polynomial $p$ is affine in $\boldsymbol{\omega}_{\mathrm{i}}$ and has degree at most four in $\boldsymbol{b}$. Although the dependence on the self-stress coefficients is linear, verifying positivity for all admissible $\boldsymbol{b}$ remains a global quartic-positivity problem. 

Exact symbolic treatment of the above two tests requires triangular decomposition or quantifier-elimination methods, whose practical complexity can grow rapidly with the number of stress and flex parameters.

\item \textbf{Polynomial equation systems and Gr\"obner bases.}
Polynomial equation systems arise both in the determination of indeterminate self-stresses $\boldsymbol{\omega}_{\mathrm{i}}$ and in the test for second-order rigidity. For indeterminate self-stresses, the equations are given by setting all minors of the reduced stress matrix equal to zero. Since the matrix entries are linear in the self-stress coefficients, an $r\times r$ minor has degree at most $r$, and hence at most $d$, in the $s$ stress variables. In the second-order rigidity test, the compatibility equations are quadratic in the variables parametrizing $\mathcal{K}_b$.

Although the input equations therefore have relatively low degree, their solution set is frequently positive-dimensional. In such cases, the objective is not merely to compute a finite list of isolated solutions, but to obtain global algebraic information about the entire solution variety. Gr\"obner basis computations may then generate intermediate polynomials of substantially higher degree and with much larger coefficients than those appearing in the original system. Moreover, the number of monomials can increase rapidly with both the number of variables and the degrees reached during elimination. These effects make Gr\"obner basis computations one of the principal bottlenecks of the workflow.

\item \textbf{Linear-algebraic feasibility problems.}
The computation of second- and third-order flexes is substantially more tractable. They reduces to determining whether a prescribed linear combination of vectors lies in the column space of the rigidity matrix. Equivalently, this vector must be orthogonal to every self-stress. The tests only have polynomial complexity and are computationally inexpensive compared with polynomial-system solving and quantified polynomial inequalities.

\end{enumerate}

Overall, the exact symbolic workflow is most practical when both the left and right nullities of the rigidity matrix are small. Based on the authors' computational experience, calculations are typically manageable when both nullities are below approximately five, although this threshold should be regarded only as an empirical guideline rather than a universal bound. When the nullities become larger, exact symbolic computations may become prohibitively expensive, and numerical polynomial solving may provide more practical alternatives. Such numerical approaches can substantially improve scalability, but they generally sacrifice exactness and may require additional certification when rigorous conclusions about rigidity or stability are needed.

\section{Conclusion}

This article translates a recently developed energy-based theory of higher-order rigidity into a practical framework for analysing polyhedral surfaces. Building on the general theory, we introduced criteria for second-order prestress stability and third-order rigidity and demonstrated their use through concrete examples exhibiting different levels of higher-order shakiness. The presentation has emphasized geometric and mechanical interpretation, following the computational workflow step by step to make the underlying theory accessible to broad audience. Although developed here for polyhedral surfaces, the workflow is applicable more broadly to a wide range of geometric constraint systems and physical near-mechanisms.

\appendix

\section*{Appendix}

\section{$n$-th order derivatives of the closure constraint} \label{app: derivatives}

This appendix provides a detailed derivation of the $n$th-order derivative of the closure constraint $\boldsymbol{f}$ for any positive integer $n$, complementing the discussion in Section~\ref{section: n-th order derivatives}.

Starting with \eqref{eq: closure 1}, let 
	\begin{align} 
R_i=
		\left[ \begin{array}{ccc}
			\cos \alpha_{i} & -\sin \alpha_{i} & 0\\
			\sin \alpha_{i} &  \cos \alpha_{i} & 0\\
			0                &                 0 & 1
		\end{array} \right]
		\left[ \begin{array}{ccc}
			1                &                 0 & 0             \\
			0                &     \cos \rho_{i} & -\sin \rho_{i}\\
			0                &     \sin \rho_{i} &  \cos \rho_{i}
		\end{array} \right],
\end{align}
and denote:
\begin{equation}
	R= \overset {\curvearrowright} {\prod_{i=1}^{n}} R_i=\left[\begin{array}{ccc}
			1                &                 0 & 0             \\
			0                &     1 &  0 \\
			0                &     0 &  1
		\end{array} \right].
\end{equation}
Note again that $\overset {\curvearrowright}{}$ means the product is evaluated by multiplying matrices in an order where those with higher indices appear on the right.  $\overset {\curvearrowleft}{}$ means the product is evaluated by multiplying matrices in an order where those with higher indices appear on the left.
Consider the coordinate frames attached to panel $i$ ($1 \le i \le n$), as shown in Figure~\ref{fig: vertex and cycle}(a). 
The direction vectors of the local axes $x_i$, $y_i$, and $z_i$ are written as the columns of the following matrix in the selected \textbf{global coordinate system}:
\begin{equation}
X_i=
\begin{bmatrix}
x_{i1} & y_{i1} & z_{i1}\\
x_{i2} & y_{i2} & z_{i2}\\
x_{i3} & y_{i3} & z_{i3}
\end{bmatrix},
\end{equation}
and we have:
\begin{equation} \label{eq: R and X_n}
	R= \overset {\curvearrowright}
	{\prod_{i=1}^{n}} R_i=\begin{bmatrix}
			1                &                 0 & 0             \\
			0                &     1 &  0 \\
			0                &     0 &  1
		\end{bmatrix}=X_nX_n^\mathrm{T}=X_nRX_n^\mathrm{T}.
\end{equation}
Next, we compute the first-order derivatives of $R$:
\begin{equation}
\begin{aligned} \label{eq: first-order derivative of R}
	\dfrac{\partial R}{\partial \rho_j} & =X_n\dfrac{\partial R}{\partial \rho_j}X_n^\mathrm{T}=X_n\left(\overset {\curvearrowright}{\prod_{i=1}^{j-1}}R_i\right)\dfrac{\mathrm{d} R_j}{\mathrm{d} \rho_j}\left(\overset {\curvearrowright} {\prod_{i=j+1}^{n}}R_i\right)X_n^\mathrm{T} 
	=X_n\left(\overset {\curvearrowright} {\prod_{i=1}^{j}}R_i\right)R_j^\mathrm{T}\dfrac{\mathrm{d} R_j}{\mathrm{d} \rho_j}\left(\overset {\curvearrowleft}{\prod_{i=1}^{j}}R_i^\mathrm{T}\right) X_n^\mathrm{T} \\ & =\left(X_n\overset {\curvearrowright}{\prod_{i=1}^{j}}R_i\right)
	\left[ \begin{array}{ccc}
		0 & 0 & 0\\
		0 & 0 & -1\\
		0 & 1 & 0
	\end{array} \right]
	\left(X_n\overset {\curvearrowright}{\prod_{i=1}^{j}}R_i\right)^\mathrm{T} =X_j
	\left[ \begin{array}{ccc}
		0 & 0 & 0\\
		0 &  0 & -1\\
		0                &                 1 & 0
	\end{array} \right]
	X_j^\mathrm{T} =
\left[
\begin{array}{ccc}
0 & -x_{j3} & x_{j2} \\
x_{j3} & 0 & -x_{j1} \\
-x_{j2} & x_{j1} & 0
\end{array}
\right].
\end{aligned}
\end{equation}
This is exactly
\begin{equation}
	\dfrac{\partial \boldsymbol{f}^v}{\partial \rho_j} = x_j.
\end{equation}

Similarly, starting from \eqref{eq: closure 2}, let
\begin{align}
\left[
\begin{array}{cc}
R_i & v_i \\
0 & 1
\end{array}
\right]
=
\left[
\begin{array}{cccc}
\cos \beta_i & -\sin \beta_i & 0 & l_i \cos \gamma_i \\
\sin \beta_i & \cos \beta_i  & 0 & l_i \sin \gamma_i \\
0 & 0 & 1 & 0 \\
0 & 0 & 0 & 1
\end{array}
\right]
\left[
\begin{array}{cccc}
1 & 0 & 0 & 0 \\
0 & \cos \rho_i & -\sin \rho_i & 0 \\
0 & \sin \rho_i & \cos \rho_i  & 0 \\
0 & 0 & 0 & 1
\end{array}
\right],
\end{align}
where $v_i=[l_i \cos \gamma_i;~l_i \sin \gamma_i;~0]$ is the position of origin $O_i$ measured in the local coordinate system built on panel $i-1$. Consequently, we have an additional constraint for all $v_i$:
\begin{equation} \label{eq: cycle constraint}
	\sum_{i=1}^{n} \overset {\curvearrowright} {\prod_{k=1}^{i-1}}R_kv_i=\begin{bmatrix}
	0 \\ 0 \\ 0
	\end{bmatrix}=X_n\sum_{i=1}^{n} \overset {\curvearrowright}{\prod_{k=1}^{i-1}}R_kv_i.
\end{equation}
The first-order derivative of the left hand side is
\begin{equation} 
	\begin{aligned}
	\dfrac{\partial}{\partial \rho_j}\left(\sum_{i=1}^{n} \overset {\curvearrowright}{\prod_{k=1}^{i-1}}R_kv_i\right) & 
	=\dfrac{\partial}{\partial \rho_j}\left(X_n\sum_{i=1}^{n} \overset {\curvearrowright}{\prod_{k=1}^{i-1}}R_kv_i\right) 
	=\sum_{i=j+1}^{n} \left(X_n \overset {\curvearrowright} {\prod_{k=1}^{j-1}} R_k\right) \dfrac{\mathrm{d} R_j}{\mathrm{d} \rho_j} \overset {\curvearrowright} {\prod_{k=j+1}^{i-1}} R_k v_i  \\ & = \left(X_n \overset {\curvearrowright} {\prod_{k=1}^{j-1}} R_k\right) \dfrac{\mathrm{d} R_j}{\mathrm{d} \rho_j} \sum_{i=j+1}^{n} \overset {\curvearrowright} {\prod_{k=j+1}^{i-1}} R_k v_i.
\end{aligned}
\end{equation}
For the last term in the above equation, we have
\begin{equation} \label{eq: relation for cycle} 
	\begin{gathered}
	\sum_{i=j+1}^{n} \overset {\curvearrowright} {\prod_{k=j+1}^{i-1}} R_k v_i= \left(\overset {\curvearrowleft} {\prod_{k=1}^{j}} R_k^\mathrm{T} \right) \sum_{i=j+1}^{n} \overset {\curvearrowright}{\prod_{k=1}^{i-1}}R_k v_i=- \left(\overset {\curvearrowleft}{\prod_{k=1}^{j}}R_k^\mathrm{T} \right) \sum_{i=1}^{j} \overset {\curvearrowright}{\prod_{k=1}^{i-1}}R_kv_i,
\end{gathered}
\end{equation}
hence 
\begin{equation}
\begin{aligned}
	\dfrac{\partial}{\partial \rho_j}\left(\sum_{i=1}^{n} \overset {\curvearrowright} {\prod_{k=1}^{i-1}}R_kv_i\right)& =-\left(X_n\overset {\curvearrowright}{\prod_{k=1}^{j-1}}R_k\right)\dfrac{\mathrm{d} R_j}{\mathrm{d} \rho_j}\left(\overset {\curvearrowleft}{\prod_{k=1}^{j}}R_k^\mathrm{T}\right)\sum_{i=1}^{j} \overset {\curvearrowright}{\prod_{k=1}^{i-1}}R_kv_i \\
	& =-\left(X_n\overset {\curvearrowright}{\prod_{k=1}^{j-1}}R_k\right)\dfrac{\mathrm{d} R_j}{\mathrm{d} \rho_j}\left(\overset {\curvearrowleft}{\prod_{k=1}^{j}}R_k^\mathrm{T}X_n^\mathrm{T}\right)\left(X_n\sum_{i=1}^{j} \overset {\curvearrowright}{\prod_{k=1}^{i-1}}R_kv_i\right)=-\dfrac{\partial R}{\partial \rho_j} \left(X_n\sum_{i=1}^{j} \overset {\curvearrowright}{\prod_{k=1}^{i-1}} R_kv_i\right) \\
	& =-\left[
\begin{array}{ccc}
0 & -x_{j3} & x_{j2} \\
x_{j3} & 0 & -x_{j1} \\
-x_{j2} & x_{j1} & 0
\end{array}
\right] \left(X_n\sum_{i=1}^{j} \overset {\curvearrowright} {\prod_{k=1}^{i-1}} R_kv_i\right) = -x_j \times O_j,
\end{aligned}
\end{equation}
where 
\begin{equation}
	O_j=X_n\sum_{i=1}^{j} \overset {\curvearrowright} {\prod_{k=1}^{i-1}} R_kv_i
\end{equation} 
is exactly the position of origin $O_j$ measured in the global coordinate system. We conclude that for every cycle:
\begin{equation}
	\dfrac{\partial \boldsymbol{f}^v}{\partial \rho_j} = - x_j \times O_j.
\end{equation}

We now calculate the second-order derivatives. For $1 \le k < j \le n$,
\begin{equation} 
	\begin{gathered}
	\dfrac{\partial^2 R}{\partial \rho_k \partial \rho_j}=X_n\dfrac{\partial^2 R}{\partial \rho_k \partial \rho_j}X_n^\mathrm{T}=X_n\left(\overset {\curvearrowright}{\prod_{i=1}^{k-1}}R_i\right) \dfrac{\mathrm{d} R_k}{\mathrm{d} \rho_k} \left(\overset {\curvearrowright}{\prod_{i=k+1}^{j-1}}R_i\right) \dfrac{\mathrm{d} R_j}{\mathrm{d} \rho_j} \left(\overset {\curvearrowright}{\prod_{i=j+1}^{n}}R_i \right) X_n^\mathrm{T}
	\end{gathered},
\end{equation}
since
\begin{equation} 
	\begin{gathered}
		\overset {\curvearrowright}{\prod_{i=k+1}^{j-1}}R_i=\overset {\curvearrowleft}{\prod_{i=1}^{k}}R_i^\mathrm{T} \overset {\curvearrowleft}{\prod_{i=j}^{n}}R_i^\mathrm{T}=\overset {\curvearrowleft}{\prod_{i=1}^{k}}R_i^\mathrm{T} \overset {\curvearrowright}{\prod_{i=1}^{j-1}}R_i
	\end{gathered},
\end{equation}
we have
\begin{equation} 
	\begin{aligned}
	\dfrac{\partial^2 R}{\partial \rho_k \partial \rho_j} & =X_n\left(\overset {\curvearrowright}{\prod_{i=1}^{k}}R_i\right) R_k^\mathrm{T}\dfrac{\mathrm{d} R_k}{\mathrm{d} \rho_k} \left(\overset {\curvearrowleft}{\prod_{i=1}^{k}}R_i^\mathrm{T} \right)X_n^\mathrm{T}X_n \left(\overset {\curvearrowright}{\prod_{i=1}^{j}}R_i\right)R_j^\mathrm{T}\dfrac{\mathrm{d} R_j}{\mathrm{d} \rho_j}\left(\overset {\curvearrowleft}{\prod_{i=1}^{j}}R_i^\mathrm{T}\right)X_n^\mathrm{T} \\
		& =\left(X_n\overset {\curvearrowright}{\prod_{i=1}^{k}}R_i\right) \left[ \begin{array}{ccc}
			0 & 0 & 0\\
			0 & 0 & -1\\
			0 & 1 & 0
		\end{array} \right] \left(X_n\overset {\curvearrowright}{\prod_{i=1}^{k}}R_i\right)^\mathrm{T}  \left(X_n\overset {\curvearrowright}{\prod_{i=1}^{j}}R_i\right)\left[ \begin{array}{ccc}
			0 & 0 & 0\\
			0 & 0 & -1\\
			0 & 1 & 0
		\end{array} \right]\left(X_n\overset {\curvearrowright}{\prod_{i=1}^{j}}R_i\right)^\mathrm{T} \\
&=
\left[
\begin{array}{ccc}
0 & -x_{k3} & x_{k2} \\
x_{k3} & 0 & -x_{k1} \\
-x_{k2} & x_{k1} & 0
\end{array}
\right]
\left[
\begin{array}{ccc}
0 & -x_{j3} & x_{j2} \\
x_{j3} & 0 & -x_{j1} \\
-x_{j2} & x_{j1} & 0
\end{array}
\right] \\
&=
\left[
\begin{array}{ccc}
-x_{k2}x_{j2}-x_{k3}x_{j3} & x_{k2}x_{j1} & x_{k3}x_{j1} \\
x_{k1}x_{j2} & -x_{k3}x_{j3}-x_{k1}x_{j1} & x_{k3}x_{j2} \\
x_{k1}x_{j3} & x_{k2}x_{j3} & -x_{k1}x_{j1}-x_{k2}x_{j2}
\end{array}
\right].
	\end{aligned}
\end{equation}

When $k=j$, 
\begin{equation} 
	\begin{aligned}
	\dfrac{\partial^2 R}{\partial \rho_j^2} & =X_n\dfrac{\partial^2 R}{\partial \rho_j^2}X_n^\mathrm{T}=X_n\left(\overset {\curvearrowright}{\prod_{i=1}^{j-1}}R_i\right)\dfrac{\mathrm{d}^2 R_j}{\mathrm{d} \rho_j^2}\left(\overset {\curvearrowright}{\prod_{i=j+1}^{n}}R_i\right)X_n^\mathrm{T}=X_n\left(\overset {\curvearrowright}{\prod_{i=1}^{j-1}}R_i\right)\dfrac{\mathrm{d}^2 R_j}{\mathrm{d} \rho_j^2}\left(\overset {\curvearrowleft}{\prod_{i=1}^{j}}R_i^\mathrm{T}\right)X_n^\mathrm{T} \\
	& =X_n\left(\overset {\curvearrowright}{\prod_{i=1}^{j}}R_i\right)R_j^\mathrm{T}\dfrac{\mathrm{d}^2 R_j}{\mathrm{d} \rho_j^2}\left(\overset {\curvearrowleft}{\prod_{i=1}^{j}}R_i^\mathrm{T}\right)X_n^\mathrm{T} =\left(X_n\overset {\curvearrowright}{\prod_{i=1}^{j}}R_i\right)
		\left[ \begin{array}{ccc}
			0 & 0 & 0\\
			0 & -1 & 0\\
			0 & 0 & -1
		\end{array} \right]
		\left(X_n\overset {\curvearrowright}{\prod_{i=1}^{j}}R_i\right)^\mathrm{T}
		\\ & =
X_j
\left[
\begin{array}{ccc}
0 & 0 & 0\\
0 & -1 & 0\\
0 & 0 & -1
\end{array}
\right]
X_j^\mathrm{T} =
\left[
\begin{array}{ccc}
-x_{j2}^2-x_{j3}^2 & x_{j2}x_{j1} & x_{j3}x_{j1} \\
x_{j1}x_{j2} & -x_{j3}^2-x_{j1}^2 & x_{j3}x_{j2} \\
x_{j1}x_{j3} & x_{j2}x_{j3} & -x_{j1}^2-x_{j2}^2
\end{array}
\right],
	\end{aligned}
\end{equation}
and when $1 \le j < k \le n$ just switch the subscripts. This goes to:
\begin{equation}
	\dfrac{\partial^2 \boldsymbol{f}^v}{\partial \rho_k \partial \rho_j} = \begin{bmatrix}	x_{j3}x_{k2} \\[4pt]
	x_{j1}x_{k3} \\[4pt]
	x_{j2}x_{k1}
	\end{bmatrix},\mathrm{~for~} k \le j.
\end{equation}

Next, when $1 \le k < j \le n$, the second-order derivative of the left hand side of \eqref{eq: cycle constraint} is
\begin{equation}
	\begin{aligned}
		\dfrac{\partial^2}{\partial \rho_k \partial \rho_j}\left(\sum_{i=1}^{n} \overset {\curvearrowright}{\prod_{l=1}^{i-1}}R_lv_i\right)& =\dfrac{\partial^2}{\partial \rho_k \partial \rho_j}\left(X_n\sum_{i=1}^{n} \overset {\curvearrowright}{\prod_{l=1}^{i-1}}R_lv_i\right) =X_n\dfrac{\partial}{\partial \rho_k}\left[\sum_{i=j+1}^{n} \left(\overset {\curvearrowright}{\prod_{l=1}^{j-1}} R_l \right) \dfrac{\mathrm{d} R_j}{\mathrm{d} \rho_j} \overset {\curvearrowright}{\prod_{l=j+1}^{i-1}} R_l v_i\right]\\
		& =X_n\sum_{i=j+1}^{n} \left(\overset {\curvearrowright}{\prod_{l=1}^{k-1}} R_l\right) \dfrac{\mathrm{d} R_k}{\mathrm{d} \rho_k} \left(\overset {\curvearrowright}{\prod_{l=k+1}^{j-1}} R_l\right) \dfrac{\mathrm{d} R_j}{\mathrm{d} \rho_j} \overset {\curvearrowright}{\prod_{l=j+1}^{i-1}}R_l v_i \\	& = X_n\left(\overset {\curvearrowright}{\prod_{l=1}^{k-1}} R_l\right) \dfrac{\mathrm{d} R_k}{\mathrm{d} \rho_k} \left(\overset {\curvearrowright} {\prod_{l=k+1}^{j-1}} R_l \right)\dfrac{\mathrm{d} R_j}{\mathrm{d} \rho_j} \sum_{i=j+1}^{n} \overset {\curvearrowright}{\prod_{l=j+1}^{i-1}} R_l v_i.
	\end{aligned}
\end{equation}
Apply \eqref{eq: relation for cycle} we obtain:  
\begin{equation}
	\begin{aligned}
		\dfrac{\partial^2}{\partial \rho_k \partial \rho_j}\left(\sum_{i=1}^{n} \overset {\curvearrowright}{\prod_{l=1}^{i-1}}R_lv_i\right)
		& =-X_n\left(\overset {\curvearrowright}{\prod_{l=1}^{k-1}} R_l\right) \dfrac{\mathrm{d} R_k}{\mathrm{d} \rho_k} \left(\overset {\curvearrowright}{\prod_{l=k+1}^{j-1}} R_l \right) \dfrac{\mathrm{d} R_j}{\mathrm{d} \rho_j} \left(\overset {\curvearrowleft}{\prod_{l=1}^{j}} R_l^\mathrm{T} \right) \sum_{i=1}^{j} \prod_{l=1}^{i-1}R_lv_i \\
		& =-\left[ \begin{array}{cccc}
			0 & -x_{3k} & x_{2k} \\
			x_{3k} & 0 & -x_{1k} \\
			-x_{2k} & x_{1k} & 0 
		\end{array} \right]\left[ \begin{array}{cccc}
			0 & -x_{3j} & x_{2j} \\
			x_{3j} & 0 & -x_{1j} \\
			-x_{2j} & x_{1j} & 0 
		\end{array} \right] \left(X_n \sum_{i=1}^{j} \overset {\curvearrowright} {\prod_{l=1}^{i-1}} R_lv_i\right) \\
		& = -x_k \times (x_j \times O_j). 
	\end{aligned}
\end{equation}
That is to say:
\begin{equation}
	\dfrac{\partial^2 \boldsymbol{f}^c}{\partial \rho_k \partial \rho_j} = - x_k \times \left(x_j \times O_j\right) \mathrm{~for~} k \le j.
\end{equation}
Similarly, when $k=j$, the result will be $-x_j \times (x_j \times O_j)$, and when $1 \le j < k \le n$ just switch the subscripts. 

The above derivation clearly shows the pattern to calculate $n$-th order derivatives. Suppose $j_1, ~j_2, ~ \cdots, ~j_n$ are the dummy variables, for $1 \le j_n \le j_{n-1} \le \cdots \le j_1 \le n$,
\begin{equation}
    \dfrac{\partial^n R}{\partial \rho_{j_n} \partial \rho_{j_{n-1}} \cdots \partial \rho_{j_1}} = x_{j_n}^\times x_{j_{n-1}}^\times \cdots x_{j_1}^\times,
\end{equation}
\begin{equation}
\dfrac{\partial^n}{\partial \rho_{j_n} \partial \rho_{j_{n-1}} \cdots \partial \rho_{j_1}}\left(\sum_{i=1}^{n} \overset {\curvearrowright}{\prod_{l=1}^{i-1}}R_lv_i\right) = - x_{j_n} \times \cdots \times (x_{j_2} \times (x_{j_1} \times O_{j_1})).
\end{equation}

\section{Explicit expressions from the Faà di Bruno's formula} \label{app: Faa}

In this appendix, we derive the higher-order derivatives of the closure constraint $\dif^i \boldsymbol{f} / \dif t^i$ and energy functional $\dif^i E / \dif t^i$ by applying the \textbf{multivariant Faà di Bruno's formula}:
\begin{equation*}
	\dfrac{\dif^i \boldsymbol{f}}{\dif t^i} = \sum \limits_{k =1}^i \sum \limits_{\substack{\mathrm{partitions~of~}i \\ \mathrm{~into~} k \mathrm{~parts}}} \dfrac{i !}{m_1 ! m_2 ! \cdots m_i!} \dfrac{\dif^k \boldsymbol{f}}{\dif \boldsymbol{\rho}^k}  \cdot  \bigotimes_{j=1}^i \left( \frac{\boldsymbol{\rho}^{(j)}}{j!} \right)^{\otimes m_j}
\end{equation*}
where the sum is over all distinct $i$-tuples of non-negative integers $(m_1, ~m_2, ~\cdots, ~m_i)$ satisfying $m_1 + 2m_2 + \cdots + im_i = i$ and $m_1 + m_2 + \cdots + m_i = k$, $\dif^k \boldsymbol{f}/\dif \boldsymbol{\rho}^k$ is a dimension $k+1$ tensor, $\boldsymbol{\rho}^{(j)} = \dif^j \boldsymbol{\rho} / \dif t^j$ at $t = 0$, $\otimes$ is the tensor product, and $\cdot$ denotes tensor contraction over the last $k$ indices of $\dif^k \boldsymbol{f}/\dif \boldsymbol{\rho}^k$ with the order-$k$ tensor $\otimes_{j=1}^i \left( p^{(j)}/j! \right)^{\otimes m_j}$.

To further clarify the calculation we write the result of the first three terms:
\begin{equation} \label{eq: flex 1}
	\dfrac{\dif \boldsymbol{f}}{\dif t} = \dfrac{\dif \boldsymbol{f}}{\dif \boldsymbol{\rho}} \cdot \boldsymbol{\rho}'
\end{equation}
\begin{equation} \label{eq: flex 2}
	\dfrac{\dif^2 \boldsymbol{f}}{\dif t^2} = \dfrac{\dif \boldsymbol{f}}{\dif \boldsymbol{\rho}} \cdot \boldsymbol{\rho}'' + \dfrac{\dif^2 \boldsymbol{f}}{\dif \boldsymbol{\rho}^2} \cdot \left(\boldsymbol{\rho}' \otimes \boldsymbol{\rho}' \right)
\end{equation}
\begin{equation} \label{eq: flex 3}
	\dfrac{\dif^3 \boldsymbol{f}}{\dif t^3} = \dfrac{\dif \boldsymbol{f}}{\dif \boldsymbol{\rho}} \cdot \boldsymbol{\rho}''' + 3\dfrac{\dif^2 \boldsymbol{f}}{\dif \boldsymbol{\rho}^2} \cdot \left(\boldsymbol{\rho}'' \otimes \boldsymbol{\rho}'\right) + \dfrac{\dif^3 \boldsymbol{f}}{\dif \boldsymbol{\rho}^3} \cdot (\boldsymbol{\rho}')^{\otimes 3} \\ 
\end{equation}
In index notation, for instance,
\begin{equation*}
	\dfrac{\dif^3 \boldsymbol{f}}{\dif \boldsymbol{\rho}^3} \cdot (\boldsymbol{\rho}')^{\otimes 3} = \sum \limits_{j_1, j_2, j_3} \dfrac{\partial^3 f_i}{\partial \rho_{j_1} \partial \rho_{j_2} \partial \rho_{j_3}} \rho_{j_1}'\rho_{j_2}'\rho_{j_3}'
\end{equation*}
Furthermore, in terms of the energy derivatives $\dif^i E / \dif t^i$:
\begin{equation} \label{eq: energy 1}
	\dfrac{\dif E}{\dif t} = \dfrac{\dif E}{\dif \boldsymbol{\rho}} \cdot \boldsymbol{\rho}' = 0 \mathrm{~~since~~}  \dfrac{\dif E}{\dif \boldsymbol{\rho}} = \boldsymbol{0} \mathrm{~~from~the~equilibrium~\eqref{eq: equilibrium}}
\end{equation} 
\begin{equation} \label{eq: energy 2}
	\begin{aligned}
		\dfrac{\dif^2 E}{\dif t^2} & = \dfrac{\dif E}{\dif \boldsymbol{\rho}} \cdot \boldsymbol{\rho}'' + \dfrac{\dif^2 E}{\dif \boldsymbol{\rho}^2} \cdot \left(\boldsymbol{\rho}'\otimes\boldsymbol{\rho}' \right) = \dfrac{\dif^2 E}{\dif \boldsymbol{\rho}^2} \cdot \left(\boldsymbol{\rho}'\otimes\boldsymbol{\rho}'\right) \\
		& = \left( \boldsymbol{\omega} \cdot \dfrac{\dif^2 \boldsymbol{f}}{\dif \boldsymbol{\rho}^2} + \dfrac{\dif^2 E}{\dif \boldsymbol{f}^2} \cdot \left(\dfrac{\dif \boldsymbol{f}}{\dif \boldsymbol{\rho}} \otimes \dfrac{\dif \boldsymbol{f}}{\dif \boldsymbol{\rho}} \right) \right) \cdot \left(\boldsymbol{\rho}' \otimes \boldsymbol{\rho}'\right) (\mathrm{the~Hessian~of~energy})
	\end{aligned}
\end{equation}
\begin{equation*} 
	\dfrac{\dif^3 E}{\dif t^3} = 3 \dfrac{\dif^2 E}{\dif \boldsymbol{\rho}^2} \cdot \left(\boldsymbol{\rho}'' \otimes \boldsymbol{\rho}'\right) + \dfrac{\dif^3 E}{\dif \boldsymbol{\rho}^3} \cdot (\boldsymbol{\rho}')^{\otimes 3} \\
\end{equation*}
The expressions for $\dif^i E / \dif \boldsymbol{\rho}^i~(i \ge 3, ~i\in \mathbb{Z}_{+})$ include the symmetrized tensor product, which is obtained by averaging over all permutations of the indices. 
\begin{equation*} 
	\dfrac{\dif^3 E}{\dif \boldsymbol{\rho}^3} =  \boldsymbol{\omega} \cdot \dfrac{\dif^3 \boldsymbol{f}}{\dif \boldsymbol{\rho}^3} + 3 \dfrac{\dif^2 E}{\dif \boldsymbol{f}^2} \cdot \left(\dfrac{\dif^2 \boldsymbol{f}}{\dif \boldsymbol{\rho}^2} \odot \dfrac{\dif \boldsymbol{f}}{\dif \boldsymbol{\rho}} \right) + \dfrac{\dif^3 E}{\dif \boldsymbol{f}^3} \cdot \left(\dfrac{\dif \boldsymbol{f}}{\dif \boldsymbol{\rho}}\right)^{\otimes 3}   
\end{equation*}
\begin{equation*}
	\dfrac{\dif^2 \boldsymbol{f}}{\dif \boldsymbol{\rho}^2} \odot \dfrac{\dif \boldsymbol{f}}{\dif \boldsymbol{\rho}} = \dfrac{1}{3} \left(\dfrac{\partial^2 f_{i_1}}{\partial \rho_{j_1}\rho_{j_2}}\dfrac{\partial f_{i_2}}{\partial \rho_{j_3}}+\dfrac{\partial^2 f_{i_1}}{\partial \rho_{j_3}\rho_{j_1}}\dfrac{\partial f_{i_2}}{\partial \rho_{j_2}}+\dfrac{\partial^2 f_{i_1}}{\partial \rho_{j_2}\rho_{j_3}}\dfrac{\partial f_{i_2}}{\partial \rho_{j_1}}\right)
\end{equation*}
and
\begin{equation} \label{eq: energy 3}
	\begin{aligned}
		\dfrac{\dif^3 E}{\dif t^3} & = 3 \left(\boldsymbol{\omega} \cdot \dfrac{\dif^2 \boldsymbol{f}}{\dif \boldsymbol{\rho}^2} + \dfrac{\dif^2 E}{\dif \boldsymbol{f}^2} \cdot \left(\dfrac{\dif \boldsymbol{f}}{\dif \boldsymbol{\rho}} \otimes \dfrac{\dif \boldsymbol{f}}{\dif \boldsymbol{\rho}} \right) \right) \cdot \left(\boldsymbol{\rho}'' \otimes \boldsymbol{\rho}'\right) \\ & \quad + \left( \boldsymbol{\omega} \cdot \dfrac{\dif^3 \boldsymbol{f}}{\dif \boldsymbol{\rho}^3} + 3 \dfrac{\dif^2 E}{\dif \boldsymbol{f}^2} \cdot \left(\dfrac{\dif^2 \boldsymbol{f}}{\dif \boldsymbol{\rho}^2} \odot \dfrac{\dif \boldsymbol{f}}{\dif \boldsymbol{\rho}} \right) + \dfrac{\dif^3 E}{\dif \boldsymbol{f}^3} \cdot \left(\dfrac{\dif \boldsymbol{f}}{\dif \boldsymbol{\rho}}\right)^{\otimes 3}  \right) \cdot (\boldsymbol{\rho}')^{\otimes 3} \\
	\end{aligned}
\end{equation}

High-order derivatives of the constraint $\boldsymbol{f}$:

\begin{equation} 
	\begin{aligned}
		\dfrac{\dif^4 \boldsymbol{f}}{\dif t^4}  & = \dfrac{ \dif \boldsymbol{f}}{\dif \boldsymbol{\rho}} \cdot \boldsymbol{\rho}^{(4)} + 4\dfrac{ \dif^2 \boldsymbol{f}}{\dif \boldsymbol{\rho}^2} \cdot \left(\boldsymbol{\rho}''' \otimes \boldsymbol{\rho}'\right) + 3\dfrac{ \dif^2 \boldsymbol{f}}{\dif \boldsymbol{\rho}^2} \cdot \left(\boldsymbol{\rho}'' \otimes \boldsymbol{\rho}''\right) \\
		& \quad + 6\dfrac{ \dif^3 \boldsymbol{f}}{\dif \boldsymbol{\rho}^3}\cdot \left(\boldsymbol{\rho}'' \otimes \boldsymbol{\rho}' \otimes \boldsymbol{\rho}' \right) + \dfrac{ \dif^4 \boldsymbol{f}}{\dif \boldsymbol{\rho}^4} \cdot \left(\boldsymbol{\rho}'\right)^{\otimes 4} \\ 
	\end{aligned}
\end{equation}
\begin{equation} 
	\begin{aligned}
		\dfrac{\dif^5 \boldsymbol{f}}{\dif t^5}  & = \dfrac{ \dif \boldsymbol{f}}{\dif \boldsymbol{\rho}}\cdot \boldsymbol{\rho}^{(5)} + 5\dfrac{ \dif^2 \boldsymbol{f}}{\dif \boldsymbol{\rho}^2} \cdot \left(\boldsymbol{\rho}^{(4)}\otimes\boldsymbol{\rho}'\right) + 10\dfrac{ \dif^2 \boldsymbol{f}}{\dif \boldsymbol{\rho}^2} \cdot \left(\boldsymbol{\rho}'''\otimes\boldsymbol{\rho}''\right) \\
		& \quad + 10 \dfrac{ \dif^3 \boldsymbol{f}}{\dif \boldsymbol{\rho}^3} \cdot \left(\boldsymbol{\rho}'''\otimes\boldsymbol{\rho}'\otimes\boldsymbol{\rho}'\right)
		+ 15 \dfrac{ \dif^3 \boldsymbol{f}}{\dif \boldsymbol{\rho}^3} \cdot \left(\boldsymbol{\rho}''\otimes\boldsymbol{\rho}''\otimes\boldsymbol{\rho}'\right) \\
		& \quad + 10 \dfrac{ \dif^4 \boldsymbol{f}}{\dif \boldsymbol{\rho}^4} \cdot \left(\boldsymbol{\rho}''\otimes (\boldsymbol{\rho}')^{\otimes 3} \right) + \dfrac{ \dif^5 \boldsymbol{f}}{\dif \boldsymbol{\rho}^5} \cdot (\boldsymbol{\rho}')^{\otimes 5} \\ 
	\end{aligned}
\end{equation}
\begin{equation} 
	\begin{aligned}
		\dfrac{\dif^6 \boldsymbol{f}}{\dif t^6}  & = \dfrac{ \dif \boldsymbol{f}}{\dif \boldsymbol{\rho}} \cdot \boldsymbol{\rho}^{(6)} + 6 \dfrac{ \dif^2 \boldsymbol{f}}{\dif \boldsymbol{\rho}^2} \cdot \left(\boldsymbol{\rho}^{(5)} \otimes \boldsymbol{\rho}'\right) + 15 \dfrac{ \dif^2 \boldsymbol{f}}{\dif \boldsymbol{\rho}^2} \cdot \left(\boldsymbol{\rho}^{(4)}\otimes\boldsymbol{\rho}''\right) + 10 \dfrac{ \dif^2 \boldsymbol{f}}{\dif \boldsymbol{\rho}^2} \cdot \left(\boldsymbol{\rho}'''\otimes\boldsymbol{\rho}'''\right) \\ 
		& + 15\dfrac{ \dif^3 \boldsymbol{f}}{\dif \boldsymbol{\rho}^3} \cdot \left(\boldsymbol{\rho}^{(4)}\otimes\boldsymbol{\rho}'\otimes\boldsymbol{\rho}'\right) + 60\dfrac{ \dif^3 \boldsymbol{f}}{\dif \boldsymbol{\rho}^3} \cdot \left(\boldsymbol{\rho}'''\otimes\boldsymbol{\rho}''\otimes\boldsymbol{\rho}'\right) + 15\dfrac{ \dif^3 \boldsymbol{f}}{\dif \boldsymbol{\rho}^3} \cdot \left(\boldsymbol{\rho}''\right)^{\otimes 3} \\
		& + 20 \dfrac{ \dif^4 \boldsymbol{f}}{\dif \boldsymbol{\rho}^4} \cdot \left(\boldsymbol{\rho}'''\otimes (\boldsymbol{\rho}')^{\otimes 3}\right) + 45 \dfrac{ \dif^4 \boldsymbol{f}}{\dif \boldsymbol{\rho}^4} \cdot (\boldsymbol{\rho}''\otimes\boldsymbol{\rho}''\otimes\boldsymbol{\rho}'\otimes\boldsymbol{\rho}') \\
		& + 15 \dfrac{ \dif^5 \boldsymbol{f}}{\dif \boldsymbol{\rho}^5} \cdot \left(\boldsymbol{\rho}''\otimes (\boldsymbol{\rho}')^{\otimes 4} \right) + \dfrac{ \dif^6 \boldsymbol{f}}{\dif \boldsymbol{\rho}^6} \cdot (\boldsymbol{\rho}')^{\otimes 6} \\ 
	\end{aligned}
\end{equation} 
High-order derivatives of the energy $E$:
\begin{equation}
	\begin{aligned} \label{eq: energy 4}
		\dfrac{\dif^4 E}{\dif t^4} & = 4\dfrac{\dif^2 E}{\dif \boldsymbol{\rho}^2} \cdot \left(\boldsymbol{\rho}'''\otimes\boldsymbol{\rho}'\right) + 3\dfrac{\dif^2 E}{\dif \boldsymbol{\rho}^2} \cdot \left(\boldsymbol{\rho}'' \otimes \boldsymbol{\rho}''\right) \\
		& \quad + 6\dfrac{\dif^3 E}{\dif \boldsymbol{\rho}^3} \cdot \left(\boldsymbol{\rho}''\otimes\boldsymbol{\rho}'\otimes\boldsymbol{\rho}'\right) + \dfrac{\dif^4 E}{\dif \boldsymbol{\rho}^4} \cdot (\boldsymbol{\rho}')^{\otimes 4}
	\end{aligned}
\end{equation}
\begin{equation*} 
	\begin{aligned}
		& = \left(\boldsymbol{\omega} \cdot \dfrac{\dif^2 \boldsymbol{f}}{\dif \boldsymbol{\rho}^2} + \dfrac{ \dif^2 E}{\dif \boldsymbol{f}^2} \cdot \left(\dfrac{ \dif \boldsymbol{f}}{\dif \boldsymbol{\rho}} \otimes \dfrac{ \dif \boldsymbol{f}}{\dif \boldsymbol{\rho}} \right)\right)\cdot \left(4\boldsymbol{\rho}'''\otimes\boldsymbol{\rho}' + 3 \boldsymbol{\rho}'' \otimes \boldsymbol{\rho}'' \right)\\
		& \quad + \left( \boldsymbol{\omega} \cdot \dfrac{ \dif^3 \boldsymbol{f}}{\dif \boldsymbol{\rho}^3} + 3 \dfrac{ \dif^2 E}{\dif \boldsymbol{f}^2} \cdot \left(\dfrac{ \dif^2 \boldsymbol{f}}{\dif \boldsymbol{\rho}^2} \odot \dfrac{ \dif \boldsymbol{f}}{\dif \boldsymbol{\rho}} \right) + \dfrac{ \dif^3 E}{\dif \boldsymbol{f}^3} \cdot \left(\dfrac{ \dif \boldsymbol{f}}{\dif \boldsymbol{\rho}}\right)^{\otimes 3}  \right) \cdot \left(6\boldsymbol{\rho}''\otimes\boldsymbol{\rho}'\otimes\boldsymbol{\rho}'\right) \\ 
	\end{aligned}
\end{equation*}
\begin{equation*} 
	\begin{aligned}
		& + \left(\boldsymbol{\omega} \cdot \dfrac{ \dif^4 \boldsymbol{f}}{\dif \boldsymbol{\rho}^4} + 4 \dfrac{ \dif^2 E}{\dif \boldsymbol{f}^2} \cdot \left(\dfrac{ \dif^3 \boldsymbol{f}}{\dif \boldsymbol{\rho}^3} \odot \dfrac{ \dif \boldsymbol{f}}{\dif \boldsymbol{\rho}} \right) + 3 \dfrac{ \dif^2 E}{\dif \boldsymbol{f}^2} \cdot \left(\dfrac{ \dif^2 \boldsymbol{f}}{\dif \boldsymbol{\rho}^2} \otimes \dfrac{ \dif^2 \boldsymbol{f}}{\dif \boldsymbol{\rho}^2}\right)\right.  \\
		& \quad + \left. 6 \dfrac{ \dif^3 E}{\dif \boldsymbol{f}^3} \cdot \left(\dfrac{ \dif^2 \boldsymbol{f}}{\dif \boldsymbol{\rho}^2} \odot \dfrac{ \dif \boldsymbol{f}}{\dif \boldsymbol{\rho}} \odot \dfrac{ \dif \boldsymbol{f}}{\dif \boldsymbol{\rho}} \right) + \dfrac{ \dif^4 E}{\dif \boldsymbol{f}^4} \cdot \left(\dfrac{ \dif \boldsymbol{f}}{\dif \boldsymbol{\rho}}\right)^{\otimes 4} \right) \cdot (\boldsymbol{\rho}')^{\otimes 4} 
	\end{aligned}
\end{equation*}
\begin{equation*} 
	\begin{aligned}
		\dfrac{\dif^5 E}{\dif t^5} & = 5\dfrac{\dif^2 E}{\dif \boldsymbol{\rho}^2} \cdot \left(\boldsymbol{\rho}^{(4)} \otimes \boldsymbol{\rho}'\right) + 10\dfrac{\dif^2 E}{\dif \boldsymbol{\rho}^2} \cdot \left(\boldsymbol{\rho}'''\otimes\boldsymbol{\rho}''\right) + 10\dfrac{\dif^3 E}{\dif \boldsymbol{\rho}^3} \cdot \left(\boldsymbol{\rho}'''\otimes\boldsymbol{\rho}'\otimes\boldsymbol{\rho}'\right) \\
		& \quad + 15\dfrac{\dif^3 E}{\dif \boldsymbol{\rho}^3} \cdot \left(\boldsymbol{\rho}''\otimes\boldsymbol{\rho}''\otimes\boldsymbol{\rho}'\right) + 10 \dfrac{\dif^4 E}{\dif \boldsymbol{\rho}^4} \cdot \left(\boldsymbol{\rho}''\otimes (\boldsymbol{\rho}')^{\otimes 3} \right) + \dfrac{\dif^5 E}{\dif \boldsymbol{\rho}^5} \cdot (\boldsymbol{\rho}')^{\otimes 5} \\ 
	\end{aligned}
\end{equation*}
\begin{equation*}
	\begin{aligned}
		& = \left(\boldsymbol{\omega} \cdot \dfrac{\dif^2 \boldsymbol{f}}{\dif \boldsymbol{\rho}^2} + \dfrac{ \dif^2 E}{\dif \boldsymbol{f}^2} \cdot \left(\dfrac{ \dif \boldsymbol{f}}{\dif \boldsymbol{\rho}} \otimes \dfrac{ \dif \boldsymbol{f}}{\dif \boldsymbol{\rho}} \right)\right)\cdot \left(5\boldsymbol{\rho}^{(4)} \otimes \boldsymbol{\rho}' + 10\boldsymbol{\rho}'''\otimes\boldsymbol{\rho}'' \right)\\
		& \quad + \left( \boldsymbol{\omega} \cdot \dfrac{ \dif^3 \boldsymbol{f}}{\dif \boldsymbol{\rho}^3} + 3 \dfrac{ \dif^2 E}{\dif \boldsymbol{f}^2} \cdot \left(\dfrac{ \dif^2 \boldsymbol{f}}{\dif \boldsymbol{\rho}^2} \odot \dfrac{ \dif \boldsymbol{f}}{\dif \boldsymbol{\rho}} \right) + \dfrac{ \dif^3 E}{\dif \boldsymbol{f}^3} \cdot \left(\dfrac{ \dif \boldsymbol{f}}{\dif \boldsymbol{\rho}}\right)^{\otimes 3}  \right) \cdot \\ 
		& \quad \quad \left(10\boldsymbol{\rho}'''\otimes\boldsymbol{\rho}'\otimes\boldsymbol{\rho}'+15\boldsymbol{\rho}''\otimes\boldsymbol{\rho}''\otimes\boldsymbol{\rho}'\right) \\
	\end{aligned}
\end{equation*}
\begin{equation*} 
	\begin{aligned}
		& + \left(\boldsymbol{\omega} \cdot \dfrac{ \dif^4 \boldsymbol{f}}{\dif \boldsymbol{\rho}^4} + 4 \dfrac{ \dif^2 E}{\dif \boldsymbol{f}^2} \cdot \left(\dfrac{ \dif^3 \boldsymbol{f}}{\dif \boldsymbol{\rho}^3} \odot \dfrac{ \dif \boldsymbol{f}}{\dif \boldsymbol{\rho}} \right) + 3 \dfrac{ \dif^2 E}{\dif \boldsymbol{f}^2} \cdot \left(\dfrac{ \dif^2 \boldsymbol{f}}{\dif \boldsymbol{\rho}^2} \otimes \dfrac{ \dif^2 \boldsymbol{f}}{\dif \boldsymbol{\rho}^2}\right)\right.  \\
		& \quad + \left. 6 \dfrac{ \dif^3 E}{\dif \boldsymbol{f}^3} \cdot \left(\dfrac{ \dif^2 \boldsymbol{f}}{\dif \boldsymbol{\rho}^2} \odot \dfrac{ \dif \boldsymbol{f}}{\dif \boldsymbol{\rho}} \odot \dfrac{ \dif \boldsymbol{f}}{\dif \boldsymbol{\rho}} \right) + \dfrac{ \dif^4 E}{\dif \boldsymbol{f}^4} \cdot \left(\dfrac{ \dif \boldsymbol{f}}{\dif \boldsymbol{\rho}}\right)^{\otimes 4} \right) \cdot \left(10\boldsymbol{\rho}''\otimes (\boldsymbol{\rho}')^{\otimes 3} \right)
	\end{aligned}
\end{equation*}
\begin{equation*} 
	\begin{aligned}
		& + \left( \boldsymbol{\omega} \cdot \dfrac{ \dif^5 \boldsymbol{f}}{\dif \boldsymbol{\rho}^5} + 5\dfrac{ \dif^2 E}{\dif \boldsymbol{f}^2} \cdot \left(\dfrac{ \dif^4 \boldsymbol{f}}{\dif \boldsymbol{\rho}^4} \odot \dfrac{ \dif \boldsymbol{f}}{\dif \boldsymbol{\rho}} \right) + 10\dfrac{ \dif^2 E}{\dif \boldsymbol{f}^2} \cdot \left(\dfrac{ \dif^3 \boldsymbol{f}}{\dif \boldsymbol{\rho}^3} \odot \dfrac{ \dif^2 \boldsymbol{f}}{\dif \boldsymbol{\rho}^2} \right) \right. \\
		& \quad + 10 \dfrac{ \dif^3 E}{\dif \boldsymbol{f}^3} \cdot \left(\dfrac{ \dif^3 \boldsymbol{f}}{\dif \boldsymbol{\rho}^3} \odot \dfrac{ \dif \boldsymbol{f}}{\dif \boldsymbol{\rho}} \odot \dfrac{ \dif \boldsymbol{f}}{\dif \boldsymbol{\rho}} \right) + 15\dfrac{ \dif^3 E}{\dif \boldsymbol{f}^3} \cdot \left(\dfrac{ \dif^2 \boldsymbol{f}}{\dif \boldsymbol{\rho}^2} \odot \dfrac{ \dif^2 \boldsymbol{f}}{\dif \boldsymbol{\rho}^2} \odot \dfrac{ \dif \boldsymbol{f}}{\dif \boldsymbol{\rho}} \right)  \\
		& \left. \quad + 10 \dfrac{ \dif^4 E}{\dif \boldsymbol{f}^4} \cdot \left(\dfrac{ \dif^2 \boldsymbol{f}}{\dif \boldsymbol{\rho}^2} \odot \left(\dfrac{ \dif \boldsymbol{f}}{\dif \boldsymbol{\rho}}\right)^{\odot 3} \right) + \dfrac{ \dif^5 E}{\dif \boldsymbol{f}^5} \cdot \left(\dfrac{ \dif \boldsymbol{f}}{\dif \boldsymbol{\rho}}\right)^{\otimes 5} \right) \cdot (\boldsymbol{\rho}')^{\otimes 5} \\ 
	\end{aligned}
\end{equation*}
\begin{equation*}
	\begin{aligned}
		\dfrac{\dif^6 E}{\dif t^6}  & = 6\dfrac{\dif^2 E}{\dif \boldsymbol{\rho}^2} \cdot \left(\boldsymbol{\rho}^{(5)} \otimes \boldsymbol{\rho}'\right) + 15\dfrac{\dif^2 E}{\dif \boldsymbol{\rho}^2} \cdot \left(\boldsymbol{\rho}^{(4)}\otimes\boldsymbol{\rho}''\right) + 10\dfrac{\dif^2 E}{\dif \boldsymbol{\rho}^2} \cdot \left(\boldsymbol{\rho}'''\otimes\boldsymbol{\rho}'''\right) \\ 
		& \quad + 15\dfrac{\dif^3 E}{\dif \boldsymbol{\rho}^3} \cdot \left(\boldsymbol{\rho}^{(4)}\otimes\boldsymbol{\rho}'\otimes \boldsymbol{\rho}'\right) + 60\dfrac{\dif^3 E}{\dif \boldsymbol{\rho}^3} \cdot \left(\boldsymbol{\rho}'''\otimes\boldsymbol{\rho}''\otimes\boldsymbol{\rho}'\right) \\
		& \quad + 15\dfrac{\dif^3 E}{\dif \boldsymbol{\rho}^3} \cdot \left(\boldsymbol{\rho}''\right)^{\otimes 3} + 20 \dfrac{\dif^4 E}{\dif \boldsymbol{\rho}^4} \cdot \left(\boldsymbol{\rho}'''\otimes (\boldsymbol{\rho}')^{\otimes 3} \right) + 45 \dfrac{\dif^4 E}{\dif \boldsymbol{\rho}^4} \cdot (\boldsymbol{\rho}''\otimes\boldsymbol{\rho}''\otimes\boldsymbol{\rho}'\otimes\boldsymbol{\rho}') \\
		& \quad + 15 \dfrac{\dif^5 E}{\dif \boldsymbol{\rho}^5} \cdot \left(\boldsymbol{\rho}''\otimes (\boldsymbol{\rho}')^{\otimes 4} \right) + \dfrac{\dif^6 E}{\dif \boldsymbol{\rho}^6} \cdot (\boldsymbol{\rho}')^{\otimes 6} \\ 
	\end{aligned}
\end{equation*}
\begin{equation*}
	\begin{aligned}
		& = \left(\boldsymbol{\omega} \cdot \dfrac{\dif^2 \boldsymbol{f}}{\dif \boldsymbol{\rho}^2} + \dfrac{ \dif^2 E}{\dif \boldsymbol{f}^2} \cdot \left(\dfrac{ \dif \boldsymbol{f}}{\dif \boldsymbol{\rho}} \otimes \dfrac{ \dif \boldsymbol{f}}{\dif \boldsymbol{\rho}} \right)\right)\cdot \left(6\boldsymbol{\rho}^{(5)} \otimes \boldsymbol{\rho}' + 15\boldsymbol{\rho}^{(4)}\otimes\boldsymbol{\rho}''+10\boldsymbol{\rho}'''\otimes\boldsymbol{\rho}''' \right)\\
		& \quad + \left( \boldsymbol{\omega} \cdot \dfrac{ \dif^3 \boldsymbol{f}}{\dif \boldsymbol{\rho}^3} + 3 \dfrac{ \dif^2 E}{\dif \boldsymbol{f}^2} \cdot \left(\dfrac{ \dif^2 \boldsymbol{f}}{\dif \boldsymbol{\rho}^2} \odot \dfrac{ \dif \boldsymbol{f}}{\dif \boldsymbol{\rho}} \right) + \dfrac{ \dif^3 E}{\dif \boldsymbol{f}^3} \cdot \left(\dfrac{ \dif \boldsymbol{f}}{\dif \boldsymbol{\rho}}\right)^{\otimes 3}  \right) \cdot \\ 
		& \quad \quad \left(15\boldsymbol{\rho}^{(4)}\otimes\boldsymbol{\rho}'\otimes \boldsymbol{\rho}'+60\boldsymbol{\rho}'''\otimes\boldsymbol{\rho}''\otimes\boldsymbol{\rho}'+15\left(\boldsymbol{\rho}''\right)^{\otimes 3}\right) \\
	\end{aligned}
\end{equation*}
\begin{equation*} 
	\begin{aligned}
		& + \left(\boldsymbol{\omega} \cdot \dfrac{ \dif^4 \boldsymbol{f}}{\dif \boldsymbol{\rho}^4} + 4 \dfrac{ \dif^2 E}{\dif \boldsymbol{f}^2} \cdot \left(\dfrac{ \dif^3 \boldsymbol{f}}{\dif \boldsymbol{\rho}^3} \odot \dfrac{ \dif \boldsymbol{f}}{\dif \boldsymbol{\rho}} \right) + 3 \dfrac{ \dif^2 E}{\dif \boldsymbol{f}^2} \cdot \left(\dfrac{ \dif^2 \boldsymbol{f}}{\dif \boldsymbol{\rho}^2} \otimes \dfrac{ \dif^2 \boldsymbol{f}}{\dif \boldsymbol{\rho}^2}\right) + 6 \dfrac{ \dif^3 E}{\dif \boldsymbol{f}^3} \cdot \left(\dfrac{ \dif^2 \boldsymbol{f}}{\dif \boldsymbol{\rho}^2} \odot \dfrac{ \dif \boldsymbol{f}}{\dif \boldsymbol{\rho}} \odot \dfrac{ \dif \boldsymbol{f}}{\dif \boldsymbol{\rho}} \right) \right.  \\
		& \quad + \left.  \dfrac{ \dif^4 E}{\dif \boldsymbol{f}^4} \cdot \left(\dfrac{ \dif \boldsymbol{f}}{\dif \boldsymbol{\rho}}\right)^{\otimes 4} \right) \cdot \left(20\boldsymbol{\rho}''\otimes (\boldsymbol{\rho}')^{\otimes 3} + 45\boldsymbol{\rho}''\otimes\boldsymbol{\rho}''\otimes\boldsymbol{\rho}'\otimes\boldsymbol{\rho}'\right)
	\end{aligned}
\end{equation*}
\begin{equation*} 
	\begin{aligned}
		& + \left( \boldsymbol{\omega} \cdot \dfrac{ \dif^5 \boldsymbol{f}}{\dif \boldsymbol{\rho}^5} + 5\dfrac{ \dif^2 E}{\dif \boldsymbol{f}^2} \cdot \left(\dfrac{ \dif^4 \boldsymbol{f}}{\dif \boldsymbol{\rho}^4} \odot \dfrac{ \dif \boldsymbol{f}}{\dif \boldsymbol{\rho}} \right) + 10\dfrac{ \dif^2 E}{\dif \boldsymbol{f}^2} \cdot \left(\dfrac{ \dif^3 \boldsymbol{f}}{\dif \boldsymbol{\rho}^3} \odot \dfrac{ \dif^2 \boldsymbol{f}}{\dif \boldsymbol{\rho}^2} \right) \right. \\
		& \quad + 10 \dfrac{ \dif^3 E}{\dif \boldsymbol{f}^3} \cdot \left(\dfrac{ \dif^3 \boldsymbol{f}}{\dif \boldsymbol{\rho}^3} \odot \dfrac{ \dif \boldsymbol{f}}{\dif \boldsymbol{\rho}} \odot \dfrac{ \dif \boldsymbol{f}}{\dif \boldsymbol{\rho}} \right) + 15\dfrac{ \dif^3 E}{\dif \boldsymbol{f}^3} \cdot \left(\dfrac{ \dif^2 \boldsymbol{f}}{\dif \boldsymbol{\rho}^2} \odot \dfrac{ \dif^2 \boldsymbol{f}}{\dif \boldsymbol{\rho}^2} \odot \dfrac{ \dif \boldsymbol{f}}{\dif \boldsymbol{\rho}} \right)  \\
		& \left. \quad + 10 \dfrac{ \dif^4 E}{\dif \boldsymbol{f}^4} \cdot \left(\dfrac{ \dif^2 \boldsymbol{f}}{\dif \boldsymbol{\rho}^2} \odot \left(\dfrac{ \dif \boldsymbol{f}}{\dif \boldsymbol{\rho}}\right)^{\odot 3} \right) + \dfrac{ \dif^5 E}{\dif \boldsymbol{f}^5} \cdot \left(\dfrac{ \dif \boldsymbol{f}}{\dif \boldsymbol{\rho}}\right)^{\otimes 5} \right) \cdot \left(15\boldsymbol{\rho}''\otimes (\boldsymbol{\rho}')^{\otimes 4} \right) \\ 
	\end{aligned}
\end{equation*}
\begin{equation*} 
	\begin{aligned}
		& + \left( \boldsymbol{\omega} \cdot \dfrac{ \dif^6 \boldsymbol{f}}{\dif \boldsymbol{\rho}^6} + 6\dfrac{ \dif^2 E}{\dif \boldsymbol{f}^2} \cdot \left(\dfrac{ \dif^5 \boldsymbol{f}}{\dif \boldsymbol{\rho}^5} \odot \dfrac{ \dif \boldsymbol{f}}{\dif \boldsymbol{\rho}} \right) + 15\dfrac{ \dif^2 E}{\dif \boldsymbol{f}^2} \cdot \left(\dfrac{ \dif^4 \boldsymbol{f}}{\dif \boldsymbol{\rho}^4} \odot \dfrac{ \dif^2 \boldsymbol{f}}{\dif \boldsymbol{\rho}^2} \right) + 10\dfrac{ \dif^2 E}{\dif \boldsymbol{f}^2} \cdot \left(\dfrac{ \dif^3 \boldsymbol{f}}{\dif \boldsymbol{\rho}^3} \otimes \dfrac{ \dif^3 \boldsymbol{f}}{\dif \boldsymbol{\rho}^3} \right) \right. \\ 
		& \quad + 15\dfrac{ \dif^3 E}{\dif \boldsymbol{f}^3} \cdot \left(\dfrac{ \dif^4 \boldsymbol{f}}{\dif \boldsymbol{\rho}^4} \odot \dfrac{ \dif \boldsymbol{f}}{\dif \boldsymbol{\rho}} \odot \dfrac{ \dif \boldsymbol{f}}{\dif \boldsymbol{\rho}} \right) + 60 \dfrac{ \dif^3 E}{\dif \boldsymbol{f}^3} \cdot \left(\dfrac{\dif^3 \boldsymbol{f}}{\dif \boldsymbol{\rho}^3} \odot \dfrac{\dif^2 \boldsymbol{f}}{\dif \boldsymbol{\rho}^2} \odot \dfrac{ \dif \boldsymbol{f}}{\dif \boldsymbol{\rho}} \right) + 15 \dfrac{ \dif^3 E}{\dif \boldsymbol{f}^3} \cdot \left(\dfrac{\dif^2 \boldsymbol{f}}{\dif \boldsymbol{\rho}^2}\right)^{\otimes 3} \\
		& \quad + 20 \dfrac{ \dif^4 E}{\dif \boldsymbol{f}^4} \cdot \left(\dfrac{\dif^3 \boldsymbol{f}}{\dif \boldsymbol{\rho}^3} \odot \left(\dfrac{ \dif \boldsymbol{f}}{\dif \boldsymbol{\rho}}\right)^{\odot 3} \right) + 45 \dfrac{ \dif^4 E}{\dif \boldsymbol{f}^4} \cdot \left(\dfrac{\dif^2 \boldsymbol{f}}{\dif \boldsymbol{\rho}^2} \odot \dfrac{\dif^2 \boldsymbol{f}}{\dif \boldsymbol{\rho}^2} \odot \dfrac{ \dif \boldsymbol{f}}{\dif \boldsymbol{\rho}} \odot \dfrac{ \dif \boldsymbol{f}}{\dif \boldsymbol{\rho}} \right) \\
		& \quad + \left. 15 \dfrac{ \dif^5 E}{\dif \boldsymbol{f}^5} \cdot \left(\dfrac{\dif^2 \boldsymbol{f}}{\dif \boldsymbol{\rho}^2} \odot \left(\dfrac{ \dif \boldsymbol{f}}{\dif \boldsymbol{\rho}}\right)^{\odot 4} \right) + \dfrac{ \dif^6 E}{\dif \boldsymbol{f}^6} \cdot \left(\dfrac{ \dif \boldsymbol{f}}{\dif \boldsymbol{\rho}}\right)^{\otimes 6} \right) \cdot (\boldsymbol{\rho}')^{\otimes 6}  \\
	\end{aligned}
\end{equation*}

\section{Proof of the weak zero lemma}
\label{app: weak zero}

\begin{lem}
\citep{gortler_higher_2025} If $\Delta \boldsymbol{\rho} \in \mathcal{K}_b \setminus \{\boldsymbol{0}\}$ is a `weak zero' of a stress matrix $\Omega$:
\begin{equation}
\Omega \cdot (\Delta \boldsymbol{\rho} \otimes \Delta \boldsymbol{\rho}) = 0
\quad \text{but} \quad
\Omega \cdot \Delta \boldsymbol{\rho} \neq \boldsymbol{0},
\end{equation}
then $\dif^2 E/\dif \boldsymbol{\rho}^2$ becomes indefinite.
\end{lem}

\begin{proof}
	Let $\boldsymbol{\phi}$ be a vector not orthogonal to $\boldsymbol{\Omega} \cdot \Delta \boldsymbol{\rho}$, i.e., $\boldsymbol{\Omega} \cdot \left(\boldsymbol{\phi} \otimes \Delta \boldsymbol{\rho}\right) \neq 0$, then we can change the sign of $\boldsymbol{\phi}$ so that $	\boldsymbol{\Omega} \cdot \left(\boldsymbol{\phi} \otimes \Delta \boldsymbol{\rho}\right) < 0$. Next, calculate the energy growth under perturbation $\Delta \boldsymbol{\rho} + s\boldsymbol{\phi}$, $s > 0$:
	\begin{equation*}
		\begin{aligned}
			\frac{\dif^2 E}{\dif \boldsymbol{\rho}^2} \cdot ((\Delta \boldsymbol{\rho} + s\boldsymbol{\phi} )\otimes (\Delta \boldsymbol{\rho} + s\boldsymbol{\phi})) &  =  \left( \boldsymbol{\Omega} + \dfrac{ \dif^2 E}{\dif \boldsymbol{f}^2} \cdot \left(\dfrac{ \dif \boldsymbol{f}}{\dif \boldsymbol{\rho}} \otimes \dfrac{ \dif \boldsymbol{f}}{\dif \boldsymbol{\rho}} \right) \right) \cdot \left(\left(\Delta \boldsymbol{\rho}+s\boldsymbol{\phi}\right) \otimes \left(\Delta \boldsymbol{\rho}+s\boldsymbol{\phi}\right)\right) \\
			& =  2s \boldsymbol{\Omega} \cdot \left(\boldsymbol{\phi} \otimes \Delta \boldsymbol{\rho}\right) + s^2 \frac{\dif^2 E}{\dif \boldsymbol{\rho}^2} \cdot \left(\boldsymbol{\phi} \otimes \boldsymbol{\phi} \right) \\
		\end{aligned}
	\end{equation*}
	Hence,
	\begin{equation*}
		\frac{\dif^2 E}{\dif \boldsymbol{\rho}^2} \cdot ((\Delta \boldsymbol{\rho} + s\boldsymbol{\phi} )\otimes (\Delta \boldsymbol{\rho} + s\boldsymbol{\phi})) < 0 \mathrm{~~if~~} s < -  \dfrac{\dif^2 E/\dif \boldsymbol{\rho}^2 \cdot \left(\boldsymbol{\phi} \otimes \boldsymbol{\phi} \right) }{ 2 \boldsymbol{\Omega} \cdot \left(\boldsymbol{\phi} \otimes \Delta \boldsymbol{\rho}\right)}
	\end{equation*}
	It means that $\dif^2 E/\dif \boldsymbol{\rho}^2$ has negative eigenvalues. 
\end{proof}


\section*{Acknowledgement}

This work was supported by the National Science Foundation under Grant No.~DMS-1929284 during Z.H.'s participation in the ICERM semester program ``Geometry of Materials, Packings and Rigid Frameworks,'' and by the RWTH Port to Europe Postdoc Fellowship, funded under the Excellence Strategy of the Federal Government and the Länder. D.R. acknowledges support from the Deutsche Forschungsgemeinschaft (DFG) under grant CRC/TRR 280 (project ID 417002380).

Z.H. is deeply grateful to Steven Gortler for his insightful discussions and ideas, which substantially influenced this work. The authors also thank Louis Theran, Miranda Holmes-Cerfon, Lisa Manning, and Christian Santangelo for helpful and stimulating discussions during the ``Geometry and Materials'' workshop and the semester program.

\bibliographystyle{cas-model2-names}

\bibliography{rigid_origami}

@article{tachi_rigid_2015,
	title = {Rigid folding of periodic origami tessellations},
	volume = {6},
	doi = {10.1090/mbk/095.1/10},
	journal = {Origami},
	author = {Tachi, Tomohiro},
	year = {2015},
	pages = {97--108},
}

@article{miura_method_1985,
	title = {Method of packaging and deployment of large membranes in space},
	volume = {618},
	journal = {Title The Institute of Space and Astronautical Science Report},
	author = {Miura, Koryo},
	year = {1985},
	pages = {1},
}

@inproceedings{demaine_zero-area_2016,
	title = {Zero-{Area} {Reciprocal} {Diagram} of {Origami}},
	volume = {2016},
	booktitle = {Proceedings of {IASS} {Annual} {Symposia}},
	publisher = {International Association for Shell and Spatial Structures (IASS)},
	author = {Demaine, Erik D. and Demaine, Martin L. and Huffman, David A. and Hull, Thomas C. and Koschitz, Duks and Tachi, Tomohiro},
	year = {2016},
	pages = {1--10},
}

@article{tachi_design_2012,
	title = {Design of infinitesimally and finitely flexible origami based on reciprocal figures},
	volume = {16},
	number = {2},
	journal = {Journal of Geometry and Graphics},
	author = {Tachi, Tomohiro},
	year = {2012},
	pages = {223--234},
}

@article{he_rigid_2019,
	title = {On rigid origami {I}: piecewise-planar paper with straight-line creases},
	volume = {475},
	copyright = {All rights reserved},
	shorttitle = {On rigid origami {I}},
	doi = {10.1098/rspa.2019.0215},
	number = {2232},
	urldate = {2020-01-05},
	journal = {Proceedings of the Royal Society A: Mathematical, Physical and Engineering Sciences},
	author = {He, Zeyuan and Guest, Simon D.},
	month = dec,
	year = {2019},
	pages = {20190215},
}

@article{connelly_second-order_1996,
	title = {Second-order rigidity and prestress stability for tensegrity frameworks},
	volume = {9},
	doi = {10.1137/s0895480192229236},
	number = {3},
	journal = {SIAM Journal on Discrete Mathematics},
	author = {Connelly, Robert and Whiteley, Walter},
	year = {1996},
	pages = {453--491},
}

@article{connelly_higher-order_1994,
	title = {Higher-order rigidity -- {What} is the proper definition?},
	volume = {11},
	doi = {10.1007/bf02574003},
	number = {2},
	journal = {Discrete \& Computational Geometry},
	author = {Connelly, Robert and Servatius, Herman},
	year = {1994},
	pages = {193--200},
}

@article{alexandrov_sufficient_1998,
	title = {Sufficient conditions for the extendibility of an n-th order flex of polyhedra},
	volume = {39},
	number = {2},
	journal = {Beitr. Algebra Geom},
	author = {Alexandrov, Victor},
	year = {1998},
	pages = {367--378},
}

@inproceedings{stachel_proposal_2007,
	title = {A proposal for a proper definition of higher-order rigidity},
	volume = {2007},
	booktitle = {Tensegrity {Workshop}},
	author = {Stachel, Hellmuth},
	year = {2007},
}

@incollection{whiteley_rigidity_2002,
	title = {Rigidity of molecular structures: generic and geometric analysis},
	shorttitle = {Rigidity of molecular structures},
	booktitle = {Rigidity theory and applications},
	publisher = {Springer},
	author = {Whiteley, Walter},
	year = {2002},
	pages = {21--46},
}

@article{crapo_statics_1982,
	title = {Statics of frameworks and motions of panel structures: a projective geometric introduction},
	shorttitle = {Statics of frameworks and motions of panel structures},
	journal = {Structural Topology, 1982, núm. 6},
	publisher = {Université du Québec à Montréal},
	author = {Crapo, Henry and Whiteley, Walter},
	year = {1982},
}

@book{connelly_frameworks_2022,
	title = {Frameworks, {Tensegrities}, and {Symmetry}},
	isbn = {978-0-521-87910-1},
	doi = {10.1017/9780511843297},
	language = {en},
	publisher = {Cambridge University Press},
	author = {Connelly, Robert and Guest, Simon D.},
	month = jan,
	year = {2022},
}

@article{watanabe_method_2009,
	title = {The method for judging rigid foldability},
	volume = {4},
	journal = {Origami},
	publisher = {AK Peters Natick, MA},
	author = {Watanabe, Naohiko and Kawaguchi, Ken-ichi},
	year = {2009},
	pages = {165--174},
}

@article{meloni_engineering_2021,
	title = {Engineering {Origami}: {A} {Comprehensive} {Review} of {Recent} {Applications}, {Design} {Methods}, and {Tools}},
	volume = {8},
	issn = {2198-3844},
	shorttitle = {Engineering {Origami}},
	doi = {10.1002/advs.202000636},
	language = {en},
	number = {13},
	urldate = {2023-07-24},
	journal = {Advanced Science},
	author = {Meloni, Marco and Cai, Jianguo and Zhang, Qian and Sang-Hoon Lee, Daniel and Li, Meng and Ma, Ruijun and Parashkevov, Teo Emilov and Feng, Jian},
	year = {2021},
	pages = {2000636},
}

@book{voss_diejenigen_1888,
	title = {Diejenigen {Flächen}, auf denen zwei {Schaaren} geodätischer {Linien} ein conjugirtes {System} bilden},
	publisher = {Verlagd. K. Akad.},
	author = {Voss, Aurel},
	year = {1888},
}

@article{misseroni_origami_2024,
	title = {Origami engineering},
	volume = {4},
	copyright = {2024 Springer Nature Limited},
	issn = {2662-8449},
	doi = {10.1038/s43586-024-00313-7},
	language = {en},
	number = {1},
	urldate = {2024-08-05},
	journal = {Nature Reviews Methods Primers},
	publisher = {Nature Publishing Group},
	author = {Misseroni, Diego and Pratapa, Phanisri P. and Liu, Ke and Kresling, Biruta and Chen, Yan and Daraio, Chiara and Paulino, Glaucio H.},
	month = jun,
	year = {2024},
	pages = {1--19},
}

@incollection{sabitov_local_1992,
	address = {Berlin, Heidelberg},
	title = {Local {Theory} of {Bendings} of {Surfaces}},
	isbn = {978-3-662-02751-6},
	doi = {10.1007/978-3-662-02751-6_3},
	language = {en},
	urldate = {2024-11-19},
	booktitle = {Geometry {III}: {Theory} of {Surfaces}},
	publisher = {Springer},
	author = {Sabitov, I. Kh.},
	editor = {Burago, Yu. D. and Zalgaller, V. A.},
	year = {1992},
	pages = {179--250},
}

@article{rembs_verbiegungen_1933,
	title = {Verbiegungen höherer {Ordnung} und ebene {Flächenrinnen}},
	volume = {36},
	issn = {1432-1823},
	doi = {10.1007/BF01188611},
	language = {de},
	number = {1},
	urldate = {2024-11-20},
	journal = {Mathematische Zeitschrift},
	author = {Rembs, Eduard},
	month = dec,
	year = {1933},
	pages = {110--121},
}

@article{efimov_theorems_1952,
	title = {Some theorems about rigidity and non-bendability (in {Russian}).},
	volume = {7},
	number = {5},
	journal = {Usp. Mat. Nauk},
	author = {Efimov, N. V.},
	year = {1952},
	pages = {215--224},
}

@article{tarnai_higher-order_1989,
	title = {Higher-order infinitesimal mechanisms},
	volume = {102},
	number = {3-4},
	journal = {Acta Technica Acad. Sci. Hung},
	author = {Tarnai, T.},
	year = {1989},
	pages = {363--378},
}

@article{salerno_how_1992,
	title = {How to recognize the order of infinitesimal mechanisms: {A} numerical approach},
	volume = {35},
	copyright = {Copyright © 1992 John Wiley \& Sons, Ltd},
	issn = {1097-0207},
	shorttitle = {How to recognize the order of infinitesimal mechanisms},
	doi = {10.1002/nme.1620350702},
	language = {en},
	number = {7},
	urldate = {2024-11-20},
	journal = {International Journal for Numerical Methods in Engineering},
	author = {Salerno, Ginevra},
	year = {1992},
	pages = {1351--1395},
}

@book{kuznetsov_underconstrained_2012,
	title = {Underconstrained {Structural} {Systems}},
	isbn = {978-1-4612-3176-9},
	doi = {10.1007/978-1-4612-3176-9},
	language = {en},
	publisher = {Springer Science \& Business Media},
	author = {Kuznetsov, E. N.},
	month = dec,
	year = {2012},
	note = {Google-Books-ID: UJHhBwAAQBAJ},
}

@article{chen_order_2011,
	title = {The order of local mobility of mechanisms},
	volume = {46},
	issn = {0094-114X},
	doi = {10.1016/j.mechmachtheory.2011.04.007},
	number = {9},
	urldate = {2024-11-20},
	journal = {Mechanism and Machine Theory},
	author = {Chen, C.},
	month = sep,
	year = {2011},
	pages = {1251--1264},
}

@article{gaspar_finite_1994,
	title = {Finite mechanisms have no higher-order rigidity},
	volume = {106},
	journal = {Acta Technica Academiae Scientiarum Hungaricae},
	publisher = {AKADEMIAI KIADO},
	author = {Gáspár, Zs and Tarnai, T.},
	year = {1994},
	pages = {119--126},
}

@article{hayakawa_panel-point_2024,
	title = {Panel-point model for rigidity and flexibility analysis of rigid origami},
	volume = {121},
	copyright = {All rights reserved},
	issn = {0925-7721},
	doi = {10.1016/j.comgeo.2024.102100},
	urldate = {2024-11-21},
	journal = {Computational Geometry},
	author = {Hayakawa, Kentaro and He, Zeyuan and Guest, Simon D.},
	month = aug,
	year = {2024},
	pages = {102100},
}

@misc{tachi_proper_2024,
	title = {A {Proper} {Definition} of {Higher} {Order} {Rigidity}},
	doi = {10.48550/arXiv.2410.15541},
	urldate = {2024-12-10},
	publisher = {arXiv},
	author = {Tachi, Tomohiro},
	month = oct,
	year = {2024},
}

@article{he_rigid_2022,
	title = {On rigid origami {III}: local rigidity analysis},
	volume = {478},
	copyright = {All rights reserved},
	shorttitle = {On rigid origami {III}},
	doi = {10.1098/rspa.2021.0589},
	number = {2258},
	urldate = {2024-12-14},
	journal = {Proceedings of the Royal Society A: Mathematical, Physical and Engineering Sciences},
	publisher = {Royal Society},
	author = {He, Zeyuan and Guest, Simon D.},
	month = feb,
	year = {2022},
	pages = {20210589},
}

@misc{gortler_higher_2025,
	title = {Higher {Order} {Rigidity} and {Energy}},
	doi = {10.48550/arXiv.2506.03108},
	urldate = {2025-06-05},
	publisher = {arXiv},
	author = {Gortler, Steven J. and Holmes-Cerfon, Miranda and Theran, Louis},
	month = jun,
	year = {2025},
	note = {arXiv:2506.03108 [math]},
}

@article{nawratil_global_2025,
	title = {A global approach for the redefinition of higher-order flexibility and rigidity},
	volume = {205},
	issn = {0094-114X},
	doi = {10.1016/j.mechmachtheory.2024.105853},
	urldate = {2025-06-06},
	journal = {Mechanism and Machine Theory},
	author = {Nawratil, Georg},
	month = mar,
	year = {2025},
	pages = {105853},
}

@misc{nawratil_flexes_2025,
	title = {On flexes associated with higher-order flexible bar-joint frameworks},
	doi = {10.48550/arXiv.2502.01124},
	urldate = {2025-06-26},
	publisher = {arXiv},
	author = {Nawratil, Georg},
	month = feb,
	year = {2025},
}

@misc{he_new_2025,
	title = {A new method for generalizing non-self-intersecting flexible polyhedra},
	doi = {10.48550/arXiv.2505.05629},
	urldate = {2025-08-05},
	publisher = {arXiv},
	author = {He, Zeyuan and Guest, Simon D.},
	month = may,
	year = {2025},
}

@article{tachi_rigid-foldable_2012,
	title = {Rigid-{Foldable} {Cylinders} and {Cells}},
	volume = {53},
	number = {4},
	journal = {Journal of the International Association for Shell and Spatial Structures},
	author = {Tachi, Tomohiro and Miura, Koryo},
	month = dec,
	year = {2012},
	pages = {217--226},
}

@article{kroy_glassy_2007,
	title = {The glassy wormlike chain},
	volume = {9},
	issn = {1367-2630},
	doi = {10.1088/1367-2630/9/11/416},
	language = {en},
	number = {11},
	urldate = {2026-04-05},
	journal = {New Journal of Physics},
	author = {Kroy, Klaus and Glaser, Jens},
	month = nov,
	year = {2007},
	pages = {416--416},
}

@article{guest_stiffness_2006,
	title = {The stiffness of prestressed frameworks: {A} unifying approach},
	volume = {43},
	issn = {0020-7683},
	shorttitle = {The stiffness of prestressed frameworks},
	doi = {10.1016/j.ijsolstr.2005.03.008},
	number = {3},
	urldate = {2026-04-05},
	journal = {International Journal of Solids and Structures},
	author = {Guest, Simon},
	month = feb,
	year = {2006},
	pages = {842--854},
}

@article{marchetti_hydrodynamics_2013,
	title = {Hydrodynamics of soft active matter},
	volume = {85},
	doi = {10.1103/RevModPhys.85.1143},
	number = {3},
	urldate = {2026-04-05},
	journal = {Reviews of Modern Physics},
	publisher = {American Physical Society},
	author = {Marchetti, M. C. and Joanny, J. F. and Ramaswamy, S. and Liverpool, T. B. and Prost, J. and Rao, Madan and Simha, R. Aditi},
	month = jul,
	year = {2013},
	pages = {1143--1189},
}

@article{ingber_tensegrity_1997,
	title = {Tensegrity: {The} architectural basis of cellular mechanotransduction},
	volume = {59},
	issn = {0066-4278, 1545-1585},
	shorttitle = {{TENSEGRITY}},
	doi = {10.1146/annurev.physiol.59.1.575},
	language = {en},
	number = {1},
	urldate = {2026-04-05},
	journal = {Annual Review of Physiology},
	author = {Ingber, D. E.},
	month = oct,
	year = {1997},
	pages = {575--599},
}

@patent{harvey_planar_1975,
	title = {Planar element especially for a constructional toy},
	nationality = {GB},
	language = {en},
	assignee = {Individual},
	number = {GB1378942A},
	urldate = {2026-07-16},
	author = {Harvey, Edward Henry},
	month = jan,
	year = {1975},
}

\end{document}